\documentclass[twoside,11pt]{article}
\date{~}
\input xy
\xyoption{all}
\usepackage{makeidx}
\usepackage{amssymb}
\usepackage{bm}
\makeindex
\font\tenfrak=eufm10
\font\sevenfrak=eufm7
\font\fivefrak=eufm5
        \newfam\frakfam \def\frak{\fam\frakfam\tenfrak}
                \textfont\frakfam=\tenfrak
\font\tenDDl=msbm10  
\font\sevenDDl=msbm7 
\font\fiveDDl=msbm5 
        \newfam\DDlfam \def\DDl{\fam\DDlfam\tenDDl} 
                \textfont\DDlfam=\tenDDl
\scriptfont\DDlfam=\sevenDDl
        \scriptscriptfont\DDlfam=\fiveDDl
\scriptfont\frakfam=\sevenfrak
        \scriptscriptfont\frakfam=\fivefrak
\def\comdG{{\bf G}}
\def\comdg{G}
\def\Box{\quad}

\long\def\nodo#1{}
\def\genfd{{\bf k}}
\def\mat22#1#2#3#4{ \left(\begin{array}{cc} #1 & #2\\ #3 & #4
\end{array}\right) } 
\newcommand{\id}{{\rm id}}
\def\dj{d\kern-.30em\raise1.25ex\vbox{\hrule width .3em height .03em}}
\def\Dj{D\rlap{\kern-.70em\raise0.75ex
    \vbox{\hrule width .3em height .03em}}}

\def\lnamedef#1{\expandafter\edef\csname lpt#1 \endcsname}
\def\lnameuse#1{\expandafter\csname lpt#1 \endcsname}
\def\myedef#1{\expandafter\edef\csname #1\endcsname}
\def\muse#1{\expandafter\csname #1\endcsname}
\def\defeq{:=}
\def\titilde#1{\tilde{\tilde{#1}}}
\newcounter{point}
\newcommand{\thispt}{{\bf \arabic{section}.\arabic{point}}}
\newcommand{\thispta}{{\bf \arabic{section}.\arabic{point}a}}
\newcommand{\ppt}{\vskip .01in\addtocounter{point}{1}
{\bf \arabic{section}.\arabic{point} }}
\newcommand{\ppta}{{\bf \arabic{section}.\arabic{point}a }}
\newcommand{\pptb}{{\bf \arabic{section}.\arabic{point}b }}
\newcommand{\pptc}{{\bf \arabic{section}.\arabic{point}c }}
\newcommand{\pptd}{
{\bf \arabic{section}.\arabic{point}d }}

\def\lnamedef#1{\expandafter\edef\csname lpt#1 \endcsname}
\def\lnameuse#1{\expandafter\csname lpt#1 \endcsname}
\newcommand\ldef[1]{\lnamedef{#1}{\arabic{section}.\arabic{point}}}
\newcommand\luse[1]{{\bf\lnameuse{#1}}}

\newcommand\minsection[1]{\section{#1}\setcounter{point}{0}}

\def\thebibliography#1{\vskip32pt plus4pt minus4pt
\centerline{\bf References}
\vglue4pt plus2pt minus2pt
\nobreak\list
 {[\arabic{enumi}]}{\settowidth\labelwidth{[#1]}\leftmargin\labelwidth
\labelsep=3pt \advance\leftmargin\labelsep
 \usecounter{enumi}}
 \def\newblock{\hskip .09em plus .33em minus .06em}
 \sloppy\clubpenalty4000\widowpenalty4000
 \sfcode`\.=1000\relax}

 \def\Fhd#1#2{\smash{\mathop{\hbox to 14mm{\rightarrowfill}}
\limits^{\scriptstyle#1}_{\scriptstyle#2}}}

\def\Fhg#1#2{\smash{\mathop{\hbox to 14mm{\leftarrowfill}}
\limits^{\scriptstyle#1}_{\scriptstyle#2}}}

 \def\fhd#1#2{\smash{\mathop{\hbox to 8mm{\rightarrowfill}}
\limits^{\scriptstyle#1}_{\scriptstyle#2}}}

\def\fhg#1#2{\smash{\mathop{\hbox to 8mm{\leftarrowfill}}
\limits^{\scriptstyle#1}_{\scriptstyle#2}}}

\def\fhnull#1#2{\smash{\mathop{\hbox to 0mm{
}}
\limits^{\scriptstyle#1}_{\scriptstyle#2}}}

\def\fhUp#1#2#3{\smash{\mathop{\hbox to 8mm{
$#2$}}
\limits^{\scriptstyle#1}_{\scriptstyle#3}}}

\def\diagram#1{\def\normalbaselines{\baselineskip=0pt
\lineskip=0pt\lineskiplimit=0pt}   \matrix{#1}}

\def\mfall{\mbox{ }\forall}
\def\wind#1{{\bf \index{#1} #1}}

\def\filtf{{\cal L}}
\def\catlSp{{\frak lSp}}
\def\catnlSp{{\frak nlSp}}
\def\catRings{{\underline{Rings}}}
\def\catCommRings{{\underline{CommRings}}}
\def\catqcoh{{\frak Qcoh}}
\def\catmonoids{{\frak Mon}}
\def\catfodc{{\frak Fodc}}

\def\abelian{abelian } 
\def\abelianM{abelian} %

\title{Noncommutative localization
\\
in noncommutative geometry}

\author{Zoran \v{S}koda}

\begin{document}
\maketitle
\begin{abstract}
The aim of these notes is to collect 
and motivate the basic localization toolbox for the geometric study
of ``spaces'' locally described by noncommutative rings
and their categories of modules.

We present the basics of {\sc Ore}
localization of rings and modules in great detail.
Common practical techniques are studied as well.
We also describe a counterexample to a folklore test principle
for Ore sets.
Localization in negatively filtered rings arising
in deformation theory is presented.
A new notion of the differential {\sc Ore} condition is introduced
in the study of localization of differential calculi.

To aid the geometrical viewpoint, localization is
studied with emphasis on descent formalism, flatness,
\abelian categories of quasicoherent sheaves and generalizations,
and natural pairs of adjoint functors for sheaf and module categories.
The key motivational theorems from the seminal works of {\sc Gabriel}
on localization, \abelian categories and schemes are quoted without proof,
as well as the related statements
of {\sc Popescu, Eilenberg-Watts, Deligne} and {\sc Rosenberg}.

The {\sc Cohn} universal localization does not have good flatness
properties, but it is determined by the localization map
already at the ring level, like the perfect localizations are.
Cohn localization is here related to the quasideterminants of {\sc
Gelfand} and {\sc Retakh}; and this may help the understanding of both
subjects.
\end{abstract}


\tableofcontents

\minsection{Introduction}{\def\thefootnote{}}

\ppt {\bf Objectives and scope.}
This is an introduction to Ore localizations and generalizations,
with {\it geometric} applications in mind.
\vskip .01in

The existing bibliography in localization theory is relatively
vast, including several monographs. Localizations proliferate
particularly in category theory, with flavours adapted to various
situations, like bicategories, toposes, {\sc Quillen}'s model
categories, triangulated categories etc. A noncommutative
algebraic geometer replaces a space with a ring or more general
'algebra', or with some category whose objects mimic modules over
the 'algebra', or mimic sheaves over a space, or he/she studies a
more general category which is
glued from such local ingredients.
This setup suggests that we may eventually need a similar
toolbox to the one used by the category and homotopy theorists;
however the simplest cases are in the area of ring and module theory.
Even so, we shift the emphasis from purely ring-theoretic questions
to more geometrical ones.
\vskip .02in

A {\bf localized ring} is typically structurally simpler than the
original, but retaining some of its features. A useful tool for a
ring theorist is a controlled but {\it substantial simplification}
of a ring, often as extreme as passing to a local ring or a
quotient skewfield. On the contrary, our main geometrical goal are
{\it those localizations which may play the role of noncommutative
analogues of (rings of regular functions on principal Zariski)
open sets} in an (affine, say) variety. Rings of functions on these
open sets may be slightly simpler than the rings of {\it global}
functions, but not as radically as when, say, passing to a local
ring. We start with the very basics of localization procedures.
The geometric notion of a {\it cover} by localizations is studied
in the noncommutative context as well. Only recent geometrically
minded works include some elements of such study, as some key
features of covers, the {\it globalization lemma} in particular,
were recognized only in mid-eighties. \vskip .02in

We use an {\it elementary method} to prove the existence and
simple properties of the Ore localized rings, in line with the
original 1931 paper of {\sc O.~Ore}~\cite{ore:31} (who however
assumed no zero divisors). Modern treatments often neglect details
when sketching the elementary method. Another modern method (to
prove existence etc. (\cite{McConnell})), following {\sc Asano}
(\cite{Asano}), is cheap, but does not give an equivalent experience
with the Ore method. Calculations similar to the ones in our
proofs appear in concrete examples when using or checking Ore
conditions {\it in practice}. We also use this method to examine
when there is an induced localization of a first order
differential calculus over a noncommutative ring, and come to a
condition, not previously observed, which we call the
``differential Ore condition''. The elementary method has the
advantage of being parallel to the calculus of (left) fractions in
general categories, which has an aura of being a difficult
subject, but is more transparent
after learning Ore localization the elementary way. \vskip .015in

Our next expositional goal is to obtain some practical criteria
for finding and dealing with Ore localizations in practice.
Folklore strategies to prove that some set is Ore are examined in
detail. In a section on Ore localization in 'negatively' filtered
rings, we explore similar methods.

\vskip .02in While Ore localization is treated in a comprehensive
introductory style, more general localizations are sketched in a
survey style. For advanced topics there are often in place good
references, assuming that the reader knows the motivation, and at
least the Ore case. Both requisites may be fulfilled by reading
the present notes. We emphasize facts important in geometry, often
omitted, or which are only folklore. In order to clear up some
sources of confusion, we sketch and compare several levels of
generality; mention competing terminologies; and examine
difficulties of geometrical interpretation and usage.

\vskip .02in We focus on localizations of the category $R-{\rm
Mod}$ of all {\it left} modules over a fixed ring $R$.
Localization in other specific categories, e.g. central bimodules
(symmetric localization, cf.~\cite{JaraVerschorenVidal}),
bimodules, and the standard approach via {\it injective hulls}
\index{injective hulls} are omitted. One of the reasons is
that often there are too few central bimodules over a noncommutative ring
and 2-sided ideals in particular. Bimodules in general
are interpreted as generalizing the maps of noncommutative
rings, as explained in the text. Generalities on localization in
arbitrary categories, and \abelian in particular, are outlined for
convenience.

\vskip .02in
As {\sc Cohn} localization can be found in other works in this volume,
we include only a short introduction
with two goals: putting it in our context and,
more importantly, relating it to the recent
subject of quasideterminants. Anybody aware of both subjects
is aware of some connection, and we try to spell it out
as precisely as we can.
\vskip .02in

\ppt {\bf Prerequisites on algebraic structures.}
Basic notions on rings and modules are freely used:
unital ring, left ideal, center, left module,
bimodule, domain (ring with no zero divisors),
skewfield (division ring), graded rings and modules,
and operations with those.

\ppt {\bf Prerequisites on categories.} The language of categories
and functors is assumed, including a (universal) initial and
terminal object, \index{initial object}\index{terminal object}
(projective=inverse) limit, colimit (= inductive/direct limit),
(co)products, adjoint functors, Yoneda lemma, and the categorical duality
(inverting arrows, dual category, dual statements about categories).

Appendix~A in~\cite{Weibel} suffices for our purposes;
for more see~\cite{BD:cats,Borceux,MacLane}.

\ppt A morphism $f : B \rightarrow C$ is epi (mono) if
for any pair $e,e'$ of morphisms
from $C$ to any $D$ (resp. from any $A$ to $B$), equality
$ef = e'f$ (resp. $fe = fe'$) implies $e = e'$.

A {\bf subobject} \index{subobject} \ldef{subobject}
is an equivalence class of monomorphisms.
A pair $(F,{\rm in})$ consisting of a functor $F$,
and a natural transformation of functors ${\rm in} : F\hookrightarrow G$,
is a {\bf subfunctor} of functor $G$, if all
${\rm in}_M : F(M)\hookrightarrow G(M)$ are monomorphisms.
Explicitly, the naturality of ${\rm in}$
reads $\forall f: M \rightarrow N$,
${\rm in}_N \circ F(f) =  G(f)\circ {\rm in}_M : F(M)\rightarrow G(N)$.
Clearly, if $F$ is a subfunctor of an additive ($\genfd$-linear)
functor $G$, between additive ($\genfd$-linear)
categories, then $F$ is additive ($\genfd$-linear) as well.

\ppt A (small) diagram $d$ in the category ${\cal C}$
will be viewed as a functor
from some (small) category $D$ into ${\cal C}$.
For fixed $D$, such functors
and natural transformations make a category ${\cal C}^D$.
Every object $c$ in ${\cal C}$ gives rise to
a constant diagram $c^D$ sending each $X$ in $D$ into $c$.
This rule extends to a functor $()^D : {\cal C}\to {\cal C}^D$.
A {\bf cone over diagram}
$d : D \to {\cal C}$ is a natural transformation
$I_c : c^D \Rightarrow d$ for some $c \in {\cal C}$.
A morphism $I_c \to I_{c'}$ of cones over $d$ is a morphism
$\phi : c'\to c$ such that $I_{c'} = I_c\circ \phi^D$.
A terminal among the cones over $D$
will be called a {\bf limiting cone} over $D$.
A {\bf colimiting cone} in ${\cal C}$ is
a limiting cone in opposite category ${\cal C}^{\rm op}$.
Consider a 'parallel' family of morphisms
$\{f_\gamma : A \to B\}_{\gamma\in\Gamma}$
as a diagram $D$ with 2 objects and $|\Gamma|$ arrows in an obvious way.
In this case, a cone over $D$ is given by
a single map $I : c \rightarrow A$.
We call the diagram
$c \stackrel{I}\to A \stackrel{f_\gamma}\Rightarrow B$
a \wind{fork diagram}.
It is called an {\bf equalizer (diagram)}
\index{equalizer}
if $I :c \rightarrow A$ is in addition a limiting cone; by abuse
of language one often calls $I$, or even $c$ an equalizer.
Equalizers in ${\cal C}^{\rm op}$ are referred to
as {\bf coequalizers}.\index{coequalizer}
A morphism $I \to A$ of the cone of
an equalizer diagram with $\Gamma = \{1,2\}$
is also called a \wind{kernel} of parallel
pair $f_1,f_2 : A \to B$. Cokernels are defined dually.

A {\bf zero} object  \index{zero object} $0$
is an object which is simultaneously initial and terminal.
In that case, composition $X\to 0\to X$ is
also denoted by $0 : X \to X$.
A (co)kernel of a {\it single} morphism $f: A \to B$
in a category with $0$
is by definition a (co)kernel of pair
$f,0 : A \to B$.

\ppt \ldef{presrefl} A functor $F : {\cal C}\to {\cal C}'$ induces
a (pullback) functor for diagrams $F\circ\, : {\cal C}^D\to ({\cal
C}')^D$. It is defined by $d \mapsto F\circ d$ for every diagram
$d: D \to {\cal C}$, and $\alpha \mapsto F(\alpha)$ where
$(F(\alpha))_M := F(\alpha_M)$ for $\alpha : d \Rightarrow d'$.

$F$ {\bf preserves} limits \index{preserving limits}
of some class ${\cal P}$ of
diagrams in ${\cal A}$ if it sends any limiting cone
$p_0 \rightarrow p$ over any
$p \in {\cal P}$ in ${\cal A}$ into a limiting
cone in ${\cal A}'$. $F$ {\bf reflects} limits \index{reflecting limits}
if any cone $p_0 \rightarrow p$ over any
$p \in {\cal P}$ in ${\cal A}$ must be a limiting cone if
$F$ sends it to a limiting cone in ${\cal A}'$.
The same holds whenever the word 'limit' is replaced by 'colimit',
and cone $p_0 \rightarrow p$ by a cocone $p \rightarrow p_0$. 

\ppt An ${\bf Ab}$-category (or {\bf preadditive category})
is a category ${\cal A}$ with an \abelian group operation $+$
on each set ${\cal A}(X,Y)$,
such that each composition map $\circ : {\cal A}(X,Y)\times {\cal A}(Y,Z)
\to{\cal A}(Y,Z)$ is bilinear. An ${\bf Ab}$-category is
{\bf additive} if it contains a zero object and pairwise,
hence all finite, products. Automatically then, finite products
agree with finite coproducts.
Recall that an additive category ${\cal A}$ is {\bf \abelian} if
each morphism $f$ in ${\cal A}$ has a kernel and a cokernel morphism, and
the kernel object of a cokernel equals the cokernel object of a kernel.
We assume that the reader is comfortable
with elementary notions on \abelian groups like
exact sequences and left (right) exact functors
in the greater generality of \abelian categories.

\ppt {\bf Gabriel-Mitchell-Popescu embedding theorem.}
{\it Every small \abelian category
is equivalent as an \abelian category
to a subcategory of category of left modules over a
certain ring $R$.} {\it Proof:} ~\cite{Popescu}.

\ppt {\bf Prerequisites on spaces of modern geometry.}
We expect familiarity with the notions of a {\it presheaf,
separated presheaf} and {\it sheaf},
and with some examples describing a geometry
via a topological space with a structure sheaf on it;
as well as the idea of gluing from some sort of local models,
behind the concepts of (super)manifolds, analytic spaces and schemes.
Earlier exposure to commutative algebraic
(or analytic) varieties and schemes
is assumed in style of some sections, but no specific facts are required;
an abstract sketch of the main features of scheme-like
theories is supplied in the text below.

\ppt {\bf Conventions.} The word {\it map} means a set-theoretic map
unless it is accompanied with a specification of
the ambient category in question
when it means a (homo)morphism, e.g. a {\it map of rings}
means a ring (homo)morphism.
The word {\it noncommutative} means
``not necessarily commutative''.
Though for many constructions this is not necessary,
we mostly deal with unital rings and modules, unless said otherwise.
{\it Ideal} without a modifier always means '1-sided (usually left) ideal'.

The symbol for inclusion $\subset$ may include the case of equality.
The unadorned tensor symbol is over ${\DDl Z}$, except for elements
in given tensor products,
like $a\otimes b := a\otimes_R b \in A\otimes_R B$.
For algebras and modules over a commutative ring,
this ring is usually denoted $\genfd$.
These conventions may be locally overwritten by contextual remarks.

\minsection{Noncommutative geometry}

{\sc Descartes} introduced the method of coordinates, which
amounted to the identification of real vector spaces with the spaces
described by the axioms of {\sc Euclid}.
{\sc Lagrange} considered more general curvilinear
coordinates in analytic mechanics to obtain
exhaustive treatments of space.
Topological spaces do not have distinguished coordinate functions,
but the whole algebra of functions suffices. The {\sc
Gelfand-Neimark} theorem (e.g.~\cite{Landsman:bookCQM})
states that the category of compact
Hausdorff topological spaces is equivalent to the opposite of
category of commutative unital $C^*$-algebras. This is
accomplished by assigning to a compact $X$ the Banach
$\star$-algebra $C(X)$ of all continuous ${\DDl C}$-valued
functions (with the supremum-norm and involution $f^*(x) =
\overline{f(x)}$). In the other direction one (re)constructs $X$
as a Gelfand spectrum of the algebra $A$, which is a space whose
points are continuous characters $\chi : A \to {\DDl C}^*$ endowed
with spectral topology. These characters appear as the
evaluation functionals $\chi_x$ on $A$ at points $x \in X$, where
$\chi_x(f) = f(x)$. Each annihilator ${\rm Ann}\,\chi = \{a \in
A\,|\, \chi(a) = 0\}$ is a maximal ideal of $A$ and all maximal
ideals are of that form.

\ppt {\bf Towards noncommutative algebraic geometry.} For any
commutative ring $R$, {\sc Grothendieck} replaced maximal ideals
from the theory of affine varieties and from the Gelfand-Neimark
picture, by arbitrary prime ideals, which behave better
functorially, and he endowed the resulting spectrum ${\rm
Spec}\,R$ with a non-Hausdorff Zariski topology and a structure
sheaf, defined with the help of commutative localization. This
amounts to a fully faithful contravariant functor ${\rm Spec}$
from $\underline{CommRings}$ to the category ${\catlSp}$ of
locally ringed spaces. In other words, the essential image of this
functor, the category of geometric affine schemes ${\rm Aff} =
{\rm Spec}(\underline{CommRings})$ is equivalent to the formal
category of affine schemes which is by the definition
$\underline{CommRings}^{\rm op}$. Notions of points, open
subspaces and sheaves are used to define ${\catlSp}$ and the
functor ${\rm Spec}$. The functor ${\rm Spec}$ takes values in a
category described in local geometrical terms, translating
algebraic concepts into geometric ones. The functor enables the
transfer of intuition and methods between algebra and geometry.
This interplay is to a large extent the basic {\it raison
d'\^{e}tre} for the subject of algebraic geometry. The spaces in
${\catlSp}$ may be glued via topologies and sheaves, and certain
limit constructions may be performed there what gives a great
flexibility in usage of range of other subcategories in $\catlSp$
e.g. schemes, algebraic spaces, formal schemes, almost schemes
etc. Useful constructions like blow-ups, quotients by actions of
groups, flat families, infinitesimal neighborhoods etc. take us
often out of the realm of affine schemes.

The dictionary between the geometric properties and abstract algebra
may be partially extended to include noncommutative algebra.
Noncommutative geometry means exploring the idea
of faithfully extending the $\rm Spec$ functor (or, analogously, the
Gelfand's functor above) to noncommutative algebras as a domain and
some geometrical universe ${\catnlSp}$ generalizing the category $\catlSp$
as a target, and thinking geometrically
of the consequences of such construction(s).
The category $\catnlSp$ should ideally contain more
general noncommutative schemes, extending the fact that spaces in
$\catlSp$ feature topology, enabling us to glue
affine schemes over such subsets. Useful examples of noncommutative
spaces are often studied in contexts which are much more restricted
then what we require for $\catnlSp$; for instance one gives
away topology or points, or one works only
with Noetherian algebras close to commutative (say of finite GK-dimension,
or finite dimensional over the center, when the latter is a field),
or the category is not big enough to include the whole of $\catlSp$.
For example, {\sc van Oystaeyen} and his school~
(\cite{vOyst:assalg, vOyst:refl,vOystWill:Groth})
consider a certain class of graded rings for which
they can use localizations to define their version of
noncommutative ${\bf Proj}$-functor.
A more restricted class of graded rings supplying examples
very close in behaviour to commutative projective varieties
is studied by {\sc Artin, Zhang, Stafford}
and others (see~\cite{ArtinZhang94,Staf:ICM} and refs. therein).
{\sc Y.~Soibelman}~(\cite{Soib:motivational})
advocates examples of natural compactifications
of moduli spaces of commutative spaces with
noncommutative spaces as points on the boundary.

Thus we often restrict ourselves either to smaller geometric
realms ${\rm Ns} \subset {\catnlSp}$ than $\catnlSp$ containing
for example only projective ``noncommutative varieties'' of some sort,
or to give up points, topological spaces in
ordinary sense and work with a more intrinsic embedding
of the category of affine schemes (now $\underline{Rings}^{\rm op}$
into some category of (pre)sheaves over $\underline{Rings}^{\rm op}$
using Yoneda lemma, Grothendieck topologies
and related concepts~(\cite{KontsRos:MP,orlov:qcsh,Ros:NcSS}).
In the commutative case, while both the spectral and functorial approaches
are interchangeably used (EGA prefers spectral, while SGA
and {\sc Demazure}-{\sc Gabriel}~\cite{Dem-Gab-eng} choose functorial;
the latter motivated with niceties in the treatment of group schemes),
the more difficult foundational constructions of
theoretical nature are done using the Yoneda embedding approach
(algebraic geometry over model categories~(\cite{Toen:higher,ToenVez:mixt});
{\sc Thomason}'s work~(\cite{TThomason}) on K-theory and derived categories
(cf. also~\cite{Balmer:psh});
${\DDl A}^1$-homotopy theory of schemes~(\cite{Morel})).

\ppt One often stops consciously half way toward the construction
of the functor $\underline{Rings}^{\rm op}\to {\catnlSp}$; e.g.
start with rings and do nothing except for introducing
a small class of open sets, e.g. commutative localizations,
ignoring other natural candidates, because {\it it is difficult
to work with them}.  Unlike in the above discussed case of restricting
the class of spaces ${\rm Ns}\subset {\catnlSp}$, we are conservative
in the details of the spectral description, thus landing in some
intermediate ``semilocal'' category ${\bf slSp}$ by means of a fully
faithful embedding $\underline{Rings}^{\rm op}\to {\bf slSp}$.

\ppta An example of ${\bf slSp}$ is as follows.
Consider the center $Z(R)$ of a ring $R$
and construct the commutative ring ${\rm Spec}\,Z(R)$. Then for each
principal open set $U$ in $Z(R)$, one localizes $R$ at $U$
(a commutative localization) and this essentially
gives the structure sheaf $U \to R_U = {\cal O}(U)$. The problem
is, what if the center is small, hence ${\rm Spec}\,Z(R)$ is small
as well,
hence all the information on $R$ is kept in a few, maybe one,
rings $R_U$, and we did not get far.
In some cases the base space  ${\rm Spec}\,Z'$ is big enough
and we may glue such spectra to interesting
more general ``schemes'' (\cite{semiquantum}).
Taking the center is not functorial,
so we have to modify the categories a bit,
to allow for pairs $(R,C)$, $C\subset Z(R)$ from the start,
and construct from them some ``space'' ${\rm Spec}_2(R,C)$
(see more in Ch.~10).
It is argued in~\cite{semiquantum} that, when $R$ is small
relative to $Z'$, this construction is a satisfactory geometrization and the
standard tools from cohomology theory may be used. They call such a
situation {\it semiquantum geometry}.\index{semiquantum geometry}

A fruitful method is to add a {\it limited class} $\{Q_i\}_i$ of
other localizations on $R-{\rm Mod}$, and {\it think} of $Q_i R$
as the structure ring $R_{U_i}$ over open subset $U_i$. However,
now $U_i$ is not really a subset in ${\rm Spec}\,Z(R)$, but rather
a ``geometric'' label for $Q_i$ viewed as certain open set on
hypothetical noncommutative ${\rm Spec}\,R$. Of course the latter
point of view is central~(\cite{KontsRos:MP,JaraVerschorenVidal,
Ros:Sploc,Versch:96,vOyst:assalg,vOystWill:Groth}) for our subject
far beyond the idea of a small enrichment of the Zariski topology
on the spectrum of the center.

In summary, restricting sharply to a small class of localizations and/or
working with small spectra, projects a coarser local description
${\bf slSp}(R)$ than often desired.

\ppt Alternatively, one may loose some information, for instance
considering the points of spectra but not the sheaves,
or types of spectra with insufficiently many points
for the reconstruction of rings.
We may think of such correspondences as nonfaithful functors
from $\bf slSp$ into some partial geometric realms
${\rm Feature}_\alpha({\catnlSp})$.

\ppt {\sc Manin} has suggested (\cite{Manin:TNG}) a functor from
graded rings into \abelian categories: to a Noetherian ring $R$
assign the quotient of the category of finitely generated graded
$R$-modules by the subcategory of the finitely generated graded
$R$-modules of finite length. In the commutative case, by a
theorem of {\sc Serre}, this category is equivalent to the
category of {\it coherent sheaves} over ${\bf Proj}\,R$. This is
one of the candidates for {\it projective noncommutative geometry}
and we view it as an example of functor of type ${\bf slSp}$.
Manin here actualizes the {\sc Grothendieck}'s advice that {\it to do
geometry one does not need the space itself but only the category
of sheaves on that would-be space}. In this spirit, Grothendieck
defined {\bf topos} (\cite{Johnstone} and \cite{Borceux}, vol.3)
as an abstract category satisfying a list of axioms, whose
consequence is that it is equivalent to the category of sheaves
${\rm Fas}\,{\cal C}$ \index{${\rm Fas}\,{\cal C}$} over some site
${\cal C}$ (a category with a Grothendieck topology). Two
different sites may give rise to the same topos, but their
cohomological behaviour will be the same. Thus they are thought of
as the same generalized space. Likewise, in algebraic geometry, we have
examples for the same heuristics, where
\abelian categories of quasicoherent sheaves of ${\cal
O}$-modules, are in the place of the topos of {\it all} sheaves of sets.
The suitable notion of a morphism between the topoi
is recognized to be 'geometrical morphism' what is also
an adjoint pair of functors with certain additional properties.
In topos theory and applications, Grothendieck
actually utilizes an interplay ('yoga') of 6 {\it standard} functors
attached to the same 'morphism'.
We shall discuss the basic pairs of adjoint functors for the categories of
modules and sheaves. They appear in the disguise of
maps of (noncommutative) rings (affine maps and their abstract
version), as bimodules for two rings, as direct and inverse maps
for ${\cal O}_X$-modules, and as localization functors.

\ppt \ldef{grothcat}
{\bf Grothendieck \index{Grothendieck category}
categories} (G.c.)~\cite{JaraVerschorenVidal,Smith:nag}.
A Grothendieck category is a cocomplete
(having all small limits) \abelian category, having
enough injectives and a small generator.
The category of left $R$-modules, and the category of all sheaves
of left $R$-modules over a fixed topological space, are G.c.'s.
Given a coalgebra $C$, the category of $C$-comodules
is G.c. Given a bialgebra $B$ and a $B$-comodule algebra $E$,
the category of relative $(E,B)$-Hopf modules is a G. c.~\cite{wisbSurvey}.
\vskip .01in

{\bf Theorem.} ({\sc P.~Gabriel},~\cite{Gab:catab}) {\it
The category $\catqcoh_X$ of quasicoherent
sheaves of ${\cal O}_X$-modules over a quasicompact quasiseparated
scheme $X$ is a Grothendieck category.
}
\vskip .01in

It is {\it not} known if $\catqcoh_X$ where $X$ is a general
scheme is cocomplete nor if it has enough injectives. This fact is
behind our decision not to strictly require our \abelian
categories of noncommutative geometry to be a G.c. (which is
fashionable). {\sc Rosenberg}~\cite{Ros:NcSch} requires the weaker
{\it property (sup)} (= {\it categories with exact limits})
introduced by {\sc Gabriel} (\cite{Gab:catab}): for any object $M$
and any ascending chain of subobjects there is a supremum
subobject, and taking such supremums commutes with taking the join
(minimum) with a fixed subobject $N\subset M$. This holds for
$R$-mod, ${\rm Fas}\,{\cal C}$ (for a small site ${\cal C}$) and
$\catqcoh_X$ (for any scheme $X$). \vskip .01in

\ppt {\bf Theorem.}
({\sc P.~Gabriel} for noetherian schemes~(\cite{Gab:catab}, Ch. VI);
{\sc A.L.~Rosenberg} in quasicompact case~(\cite{Ros:recon});
and in general case~(\cite{Ros:USpSch}))
\indent{\it Every scheme $X$ can be reconstructed
from the \abelian category $\catqcoh_X$ uniquely up to an isomorphism
of schemes.} \vskip .02in

This motivates the promotion to a ``space'' of any member of a class
of \abelian categories, usually required to obey some additional
axioms, permitting for (some variant) of $R$-mod
($R$ possibly noncommutative) and $\catqcoh_X$ as prime examples.
A distinguished object ${\cal O}$
in ${\cal A}$, corresponding to the structure sheaf is often
useful part of a data, even at the abstract level, hence the
spaces could be actually pairs $({\cal A}, {\cal O})$. The
study of functors for the categories of modules and categories
of sheaves shows that there is a special role for functors having
various exactness properties (\cite{Ros:NcSch}), e.g. having a
right adjoint, hence such properties are often required in general.
Gluing categories over localizations, and variants
thereof, should be interpreted in good cases as gluings of spaces
from local models. In noncommutative geometry, a local model is usually the
full category of modules over a noncommutative ring. \vskip .02in

\ppt The so-called {\it derived algebraic geometry},
treating in more natural terms and globalizing the infinitesimal
picture of moduli spaces governed by the deformation theory,
appeared recently (\cite{Behrend:dgschemes,Toen:higher}).
Its cousin, {\it homotopical algebraic geometry} appeared promising
in the study of homotopy theories for algebraic varieties,
and also in using the reasoning of algebraic geometry for
ring spectra of homotopy theory and for their globalization.
In such generalizations of algebraic geometry the basic gadgets are
{\it higher categories} (e.g. simplicially enriched, DG, Segal, $A_\infty$,
cf.~\cite{Drinf:DG,Keller:introAinf,Porter:Scat,Toen:higher,ToenVez:mixt}).
The lack of smoothness in some examples of moduli spaces is now
explained as an artifact of the truncation process replacing the
natural and smooth 'derived moduli spaces' by ordinary moduli spaces
('{\it hidden smoothness principle}' due to {\sc Bondal, Deligne, Drinfeld,
Kapranov, Hinich, Kontsevich}...). \vskip .02in

Part of the relevant structure here may be already expressed by
replacing rings by differential graded algebras (dga-s)
(\cite{Behrend:dgschemes}), or, more generally, by introducing
sequences of higher (e.g. '{\sc Massey}') products, as in the
theory of $A_\infty$ (strongly homotopy associative) and
$L_\infty$ (strongly homotopy Lie) algebras. Such generalizations
and special requirements needed to do localization in such
enriched settings, are beyond the scope of the present article. A
noncommutative algebraic geometry framework designed by {\sc
O.A.~Laudal}~(\cite{Laudal:MPI}), with emphasis on the problem of
noncommutative deformation of moduli, implicitly includes the
higher Massey products as well. In the viewpoint put forward by
{\sc Kontsevich} and {\sc Fukaya}, some of the 'dualities' of
modern mathematical physics, e.g. the {\it homological mirror
symmetry}, involve $A_\infty$-categories defined in terms
of geometric data~
(\cite{SoibKonts:mirrorToric,Manin:mmm,Soib:motivational}). The
so-called {\bf quantization}\index{quantization}
~(\cite{Connes:book,Konts:defqalgvar,Landsman:bookCQM,Yek:def}) in
its many versions is generally of deformational and noncommutative
nature. Thus it is not surprising that the formalisms combining
the noncommutative and homological (or even homotopical)
structures benefit from the geometrically sound models of quantum
physics. {\sc Manin} suggested that a more systematic content of a
similar nature exists, programmatically named {\it quantized
algebraic geometry}, which may shed light on hidden aspects
of the geometry of (commutative) varieties, including the deep
subject of motives. \vskip .02in

An interesting interplay of derived categories
of coherent sheaves on varieties
and their close analogues among other triangulated categories,
motivates some 'noncommutative' geometry of the latter
(\cite{BonVdB:gen,BonOrl:ICM}).
Triangulated categories are also only
a ``truncation'' of some other higher categories.\vskip .02in

One should also mention that some
important classes of rings in {\it quantum algebra},
for example quantum groups, may be constructed using categories
of (perverse) sheaves over certain commutative configuration spaces
(\cite{Lusztig:qgrbk}).
Thus the structure of various sheaf categories is an ever repeating theme
which relates the commutative and noncommutative world.
See the essay~\cite{Cart:madday} for further motivation.

\minsection{Abstract localization}

\label{sec:absloc}

We discuss localization of
1. algebraic structures; 2. categories.
These two types are related: typically a localization of a ring $R$
induces a localization of the category $R-{\rm Mod}$
of left modules over $R$.

A recipe ${\cal G}$ for a localization takes as input
a structure $R$ (monoid, lattice, ring), or a category ${\cal A}$,
and a distinguished data $\Sigma$ in $R$ (or ${\cal A}$).
The localizing data $\Sigma$ is selected from some class ${\cal U}(R)$
of structural data, for example elements, endomorphisms or ideals of $R$;
similarly ${\cal U}({\cal A})$ could be a class of subcategories or
collection of morphisms in ${\cal A}$.
Usually not all obvious subclasses of
${\cal U}(R)$ may serve as distinguished data for ${\cal G}$,
and some 'localizability' conditions apply.

A localization procedure  ${\cal G}(R,\Sigma)$
should replace $R$ by another object
$Y$ and a map $i : R \to Y$, which induces,
for given ${\cal G}$, some canonical correspondence
${\cal G}(i) :\Sigma\rightsquigarrow\Sigma_*$
between the localization data $\Sigma$
and some other data $\Sigma_*$ chosen from ${\cal U}(Y)$.
The subclass $\Sigma_*$ should satisfy some natural
requirement, for example that it consist of invertible
elements. Pair $(i,Y)$ should be in some sense smallest,
or universal among all candidates satisfying the
given requirements. For given requirements
only certain collections $\Sigma$ built from
elements in ${\cal U}(R)$ give rise to a universal $(i,Y)$.
Such $\Sigma$ are generically called {\it localizable}
and the resulting $Y$ is denoted $\Sigma^{-1}R$.

In the case of a category ${\cal C}$, a map $i$ is
replaced by a localization functor
$Q^* : {\cal C}\rightarrow \Sigma^{-1}{\cal C}$.
In this article, a localization of a category
will be equivalent to an abstract localization with respect to
a class of morphisms $\Sigma$ in ${\cal C}$,
often using some other equivalent data (e.g. 'localizing subcategory').
Following~\cite{GZ}, we sketch the general case of
a localization at a class of morphisms $\Sigma$, cf. also~\cite{Borceux}.

\ppt An {\bf abstract 1-diagram} ${\cal E}$ is a structure weaker
then a category: it consists of a class  ${\rm Ob}\,{\cal E}$ of
objects and a class ${\rm Mor}\,{\cal E}$ of morphisms equipped
with a source and a target maps $
\mbox{\sc s},\mbox{\sc t}:{\rm Mor}\,{\cal E}\to {\rm Ob}\,{\cal E}$.
No composition, or identity morphisms are supplied.
As usual, for two objects $A,B$ by ${\cal E}(A,B)$
we denote class of morphisms $f$ with $\mbox{\sc s}(f) = A$ and
$\mbox{\sc t}(f) = B$. If each ${\cal E}(A,B)$ is a set,
one may use word (multiple-edge) graph instead. If ${\cal E},{\cal C}$ are
diagrams, an ${\cal E}$-{\bf diagram} in ${\cal C}$, or a morphism
$d : {\cal E}\to {\cal C}$, is any pair of maps ${\rm Ob}\,{\cal
E}\to {\rm Ob}\,{\cal C}$ and ${\rm Mor}\,{\cal E}\to {\rm
Mor}\,{\cal C}$ which commute with source and target maps. Small
abstract 1-diagrams and their morphisms form a category ${\frak
Diagr}_1$. To each category one assigns its underlying abstract
diagram. This correspondence induces a forgetful functor from the
category ${\frak Cat}$ of small categories to ${\frak Diagr}_1$.
The construction of a category of paths below provides the left
adjoint to this functor.

If $n\geq 0$ is an integer,
a {\bf path} of length $n$ from $A$ to $B$ in an abstract diagram ${\cal E}$
is a tuple $(A,f_1,f_2,\ldots,f_n,B)$, where $A$ is an object
and $f_i$ are morphisms in  ${\cal E}$, such that
$\mbox{\sc s}(f_{i+1}) = \mbox{\sc t}(f_i)$ for $i = 1,\ldots, n-1$,
and $\mbox{\sc s}(f_{1}) = A$, $\mbox{\sc t}(f_n) = B$ if $n>0$,
and $A = B$ if $n = 0$.
For any abstract 1-diagram ${\cal E}$ define
a {\bf category} ${\rm Pa}\,{\cal E}$~\index{${\rm Pa}\,{\cal E}$}
{\bf of paths} in ${\cal E}$\index{category of paths}~as follows.
The class of objects is
$${\rm Ob}\,{\rm Pa}\,{\cal E} := {\rm Ob}\,{\cal E}$$
and the morphism class $({\rm Pa}\,{\cal E})(A,B)$
consists of all paths from $A$ to $B$.
One declares ${\rm Id}_A := (A,A)$,
$\mbox{\sc s}' (A,f_1,\ldots,f_n,B) = A$ and
$\mbox{\sc t}' (A,f_1,\ldots,f_n,B) = B$ to be the identity morphisms,
and the source and target maps for ${\rm Pa}\,{\cal E}(A,B)$.
The composition rule is
\[ (A,f_1,\ldots,f_n,B) \circ (B,g_1,\ldots,g_m,C) =
(A,f_1,\ldots,f_n,g_1,\ldots,g_m,C).\]
If each $({\rm Pa}\,{\cal E})(A,B)$ is small we indeed obtain a category.

Consider the canonical ${\cal E}$-diagram
$i_{\cal E}: {\cal E}\to {\rm Pa}\,{\cal E}$ which is
tautological on objects as well as on paths of length 1.
${\rm Pa}\,{\cal E}$ has the following universal property:
an ${\cal E}$-{\bf diagram} $d$ in any category ${\cal C}$
gives rise to a unique functor $d' : {\rm Pa}\,{\cal E}\to {\cal C}$
such that $d = d' \circ i_{\cal E}$.

\ppt Let $\Sigma$ be a family of morphisms in
category ${\cal C}$. If $J : {\cal C}\rightarrow {\cal D}$
is any functor let $\Sigma_* := J(\Sigma)$ be
the class of all morphisms $J(f)$ where $f \in \Sigma$.
Given ${\cal C}$ and $\Sigma$, consider the diagram scheme
${\cal E} = {\cal E}({\cal C},\Sigma)$ with
${\rm Ob}\,{\cal E} := {\rm Ob}\,{\cal C}$
and ${\rm Mor}\,{\cal E} := {\rm Mor}\,{\cal C}\coprod\Sigma$,
$\mbox{\sc s}_{\cal E} = \mbox{\sc s}\coprod \mbox{\sc s}|_\Sigma$,
$\mbox{\sc t}_{\cal E} = \mbox{\sc t}\coprod \mbox{\sc t}|_\Sigma$.
One has canonical inclusions
$\mbox{\sc in}_1 : {\rm Mor}\,{\cal C}\to {\rm Mor}\,{\cal E}$,
$\mbox{\sc in}_2 : \Sigma \hookrightarrow  {\rm Mor}\,{\cal E}$.
Let $\sim$ be the smallest equivalence relation
on ${\rm Pa}\,{\cal E}$ such that
\[\begin{array}{c}
(\mbox{\sc in}_1 v) \circ (\mbox{\sc in}_1 u) \sim \mbox{\sc in}_1 (v\circ u)
\mbox{ if }  v\circ u \mbox{ is defined in } {\cal C},\\
\mbox{\sc in}_1(\id_A) \sim (A,A),\,\,\,\,\,A\in {\cal C},\\
\left.\begin{array}{l}
(\mbox{\sc in}_2 f) \circ (\mbox{\sc in}_1 f)
\sim (\mbox{\sc s}(f),\mbox{\sc s}(f)) \\
(\mbox{\sc in}_1 f) \circ (\mbox{\sc in}_2 f)
\sim (\mbox{\sc t}(f),\mbox{\sc t}(f))
\end{array}\right\rbrace f \in \Sigma.
\end{array}\]
It is direct to show that operation $\circ$
induces a composition on classes of morphisms with respect to
this particular equivalence relation. In this way we obtain
a quotient $\Sigma^{-1}{\cal C}$ of
the category ${\rm Pa}\,{\cal E}$ together with
the canonical functor $Q^*_\Sigma : {\cal C}\to{\rm Pa}\,{\cal E}$
which is tautological on objects and equals
$i_{\cal E} \circ \mbox{\sc in}_1$
followed by the projection to the classes of
equivalence on morphisms.

\ppt {\bf Proposition.} {\it
If $f \in \Sigma$ then the
functor $Q^*_\Sigma : {\cal C}\to\Sigma^{-1}{\cal C}$
sends $f$ to an invertible map $Q^*_\Sigma(f)$.
If $T :  {\cal C}\to{\cal D}$ is any functor such
that $T(s)$ is invertible for any $s \in \Sigma$ then there is
a unique functor $H :\Sigma^{-1}{\cal A}\to {\cal B}$
such that $T = H \circ Q^*_\Sigma$.
}

$\Sigma^{-1}{\cal C}$ is
\wind{category of fractions}
of ${\cal C}$ at $\Sigma$. This construction has a defect, in that
there is no general recipe to determine when two morphisms
in ${\rm Pa}\,{\cal E}$ represent the same morphism in $\Sigma^{-1}{\cal C}$.
If $\Sigma$ satisfies the Ore conditions, below,
there is one.

\ppt  \ldef{absloc} {\bf Proposition.}~\cite{GZ} {\it Let
$T^* \dashv T_*$ be an adjoint pair of functors
(this notation means that $T^*$ is left adjoint to $T_*$),
where $T^* : {\cal A}\rightarrow {\cal B}$.
with adjunction counit $\epsilon : T^* T_* \Rightarrow 1_{\cal B}$.
Let $\Sigma = \Sigma(T_*)$ be the class of all morphisms
$f$ in ${\cal A}$ such that $T^*(f)$ is invertible,
and $Q^*_\Sigma : {\cal A} \to \Sigma^{-1}{\cal A}$ the natural
functor.
Then the following are equivalent:

(i) $T_*$ is fully faithful.

(ii) $\epsilon : T^* T_* \Rightarrow 1_{\cal B}$
is an isomorphism of functors.

(iii) The unique functor $H : \Sigma^{-1}{\cal A}\to {\cal B}$
such that $T^* = H\circ Q^*_\Sigma$ is an equivalence;
in particular $Q^*_\Sigma$ has a right adjoint $Q_{\Sigma *}$.

(iv) (If ${\cal A}$ is small)
For each category ${\cal X}$, functor
${\frak Cat}(-,{\cal X}) : {\frak Cat}({\cal B},{\cal X})
\to  {\frak Cat}({\cal A},{\cal X})$ is fully faithful.
}

\vskip .02in Throughout the paper, any functor $T^*$
agreeing with a functor $Q^*_\Sigma : {\cal C}\to \Sigma^{-1}{\cal C}$
as above up to category equivalences will be referred to as
a {\bf localization functor}\index{localization functor}.
A functor $T^*$ satisfying (i)-(iii)
will be referred to as {\bf a continuous localization functor}.\vskip .02in

\ppt (\ldef{compat}\cite{LuntsRosMP})
{\it Suppose $Q^* : {\cal A}\to {\cal B}$ is a localization functor
(cf.~\luse{absloc}), and $F : {\cal A}\to {\cal A}$ an endofunctor.
If there is a functor $G : {\cal B}\to {\cal B}$
and a natural equivalence of functors
$\alpha : Q^*\circ F \Rightarrow G\circ Q^*$
then there is a unique functor
$F_{\cal B} : {\cal B}\to {\cal B}$ such that
$Q^*\circ F = F_{\cal B}\circ Q^*$. In that case,
we say that {\bf $F$ is compatible with $Q^*$}.}

{\it Proof.} Suppose $f : M \to N$ is a morphism in ${\cal A}$.
Suppose that $Q^*(f)$ is invertible. Then $GQ^*(f) : GQ^*(M)\to
GQ^*(N)$ also has some inverse $s$. The naturality of $\alpha$ and
$\alpha^{-1}$ implies
\[\begin{array}{l}
\alpha_M^{-1} \circ s \circ \alpha_N \circ Q^*F(f) =
\alpha_M^{-1}\circ s \circ GQ^*(f) \circ \alpha_M  = \id_M,
\\
Q^*F(f)\circ \alpha_M^{-1}\circ s\circ \alpha_N =
\alpha_N^{-1}\circ GQ^*(f) \circ s\circ \alpha_N = \id_N,
\end{array}\]
hence $\alpha_M^{-1} \circ s \circ \alpha_N : Q^*F(N)\to Q^*F(M)$ is
the inverse of $Q^*F(f)$. The conclusion is that for any $f$ with
$Q^*(f)$ invertible, $Q^*F(f)$ is invertible as well. In other words,
(by the universal property of the localization),
functor $Q^*F$ factors through the quotient category ${\cal B}$, i.e.
$\exists!\,F_{\cal B} : {\cal B}\to {\cal B}$ with
$Q^*\circ F = F_{\cal B}\circ Q^*$. Q.E.D.

\minsection{Ore localization for monoids}

\indent \ppt A \wind{semigroup} is a set $R$ with a binary {\it
associative} operation. A semigroup with unit element $1 \in R$ is
called a \wind{monoid}. By definition, maps of semigroups are set
maps which respect the multiplication, and maps of monoids should
preserve unit element as well. Monoids and maps of monoids form a
category $\catmonoids$, which has arbitrary products. The notion
of a submonoid is the obvious one.

A subset $S$ of a monoid $R$ is {\bf multiplicative}
\index{multiplicative subset} if $1 \in S$ and whenever $s_1,s_2
\in S$ then $s_1 s_2 \in S$. For a set $S_1 \subset R$ there is a
smallest multiplicative subset $S\subset R$ containing $S_1$,
namely the set of all products $s_1 \cdots s_n$ where $s_i \in
S_1$, including the product of the empty set of elements which
equals $1$ by definition. We say that $S$ is
\wind{multiplicatively generated} by $S_1$.

\vskip .02in
\ppt A multiplicative subset $S\subset R$ is a \wind{left Ore set} if
 \begin{itemize}
\item $    (\forall s \in S\,\, \forall r \in R \,\,
\exists s' \in S \,\,\exists r' \in R) \, ( r' s =  s' r )$
(left Ore condition);
\item $(\forall n_1, n_2 \in R\,\, \forall s \in S)
\,(n_1 s = n_2 s) \Rightarrow (\exists s' \in S, \,s'n_1 = s'n_2)$

(left reversibility).
\end{itemize}

\ppt In traditional ring-theoretic terminology, $S$ is a left {\bf
Ore set} if the first condition holds and $S$ is a \wind{left
denominator set} if both conditions hold. We often say ``left Ore
set'' for a left denominator set, as is increasingly common among
geometers, and the notion of satisfying just the left Ore
condition may be said simply ``satisfying left Ore condition''. By
saying (plural:) ``left Ore condition{\it s}'' we subsume both the
left Ore condition and the left reversibility.

\ppt A monoid $R$ can be viewed as a small category ${\rm Cat}(R)$
with a single object $R$.
Left multiplication by an element $a\in R$ is a morphism in
${\rm Cat}R$ denoted by $L_a$.
We compose the morphisms by composing the maps. Any small category having
one single object is clearly equivalent to ${\rm Cat}(R)$ for a suitable $R$.

This suggests a generalization of the notion of a denominator set
(as well as its applications below) by replacing ${\rm Cat}(R)$ by
an arbitrary category. A \wind{multiplicative system} in a
category ${\cal A}$ is a class $\Sigma$ of morphisms in ${\cal A}$
where all identity morphisms $1_A$, where $A \in {\rm Ob}\,{\cal
A}$, are in $\Sigma$, and for any two composable morphism $s, t
\in \Sigma$ (i.e., the target (range) of $t$  matches the source
(domain) of $s$), also $s\circ t \in S$.

A multiplicative system $\Sigma$ satisfies the left Ore
conditions if it satisfies the ordinary left Ore condition with
all quantifiers conditioned on matching
of the source and target maps appropriately.

More precisely, $\Sigma$ satisfies the left Ore condition if
\[ \forall (s : A\to B) \in \Sigma,\,\forall r : A' \to B,\,
\exists (s': D \to A') \in \Sigma,\,\exists r' : D \rightarrow A,\]
so that $r' \circ s =  s' \circ r$.
$S$ satisfies left reversibility ('simplifiability') if
\[\begin{array}{c}
(\forall n_1, n_2 : A \to B,\, \forall (s: C \to A)  \in \Sigma)
\,(n_1 \circ s = n_2 \circ s) \\ \,\,\,\,\,\Rightarrow
(\exists (s': B \to D) \in \Sigma, \,
s'\circ n_1 = s'\circ n_2).\end{array}\]
We may picture the left simplifiability by the diagram
\[
 \diagram{
C \fhd{s}{} A \fhd{\fhd{n_1}{}}{n_2}
B \stackrel{s'}{\dashrightarrow} D
}
\]
We say that $\Sigma$ is a left denominator system, or
equivalently, that the pair $({\cal A},\Sigma)$ forms a {\bf left
calculus of fractions} \index{left calculus of fractions} if the
left Ore and left simplifiability condition hold. The
book~\cite{Pop:cats} has a good graphical treatment of that
subject. See also~\cite{Borceux,Ehresmann,GZ,Popescu}.

\ppt {\bf Lemma.} {\it Let $f: R\rightarrow R'$
be a surjective map of monoids and $S\subset R$ left Ore.
Then $f(S)$ is left Ore in $R'$.
}

\ppt Let ${\cal D}$ be some category of monoids with additional structure,
i.e. a category with a faithful functor
$(-)_{\rm mon} : {\cal D} \rightarrow \catmonoids$
preserving and reflecting finite equalizers.
If $R$ is an object in ${\cal D}$,
a multiplicative set in $R$ is
by definition any multiplicative set $S \subset (R)_{\rm mon}$.

{\bf Definition.}  {\it Given a multiplicative set $S$ in $R \in
{\cal D}$ we introduce category ${\cal C}_{\cal D}(R,S)$
\index{${\cal C}_{\cal D}(R,S)$} as follows. The objects of ${\cal
C}_{\cal D}(R,S)$ are all pairs $(j, Y)$ where $Y \in {\rm
Ob}\,{\cal D}$ and $j : R \rightarrow Y$ is a morphism in ${\cal
D}$ satisfying
\begin{itemize}
\item $(\forall s \in S)$ $(\exists u \in Y)$
$( u j(s) = j(s) u = 1)$ in $(Y)_{\rm mon}$;
\end{itemize}
The morphisms of pairs $\sigma : (j,Y) \rightarrow (j',Y')$ are
precisely those morphisms $\sigma : Y \rightarrow Y'$ in ${\cal
D}$, for which $\sigma \circ j = j'$. }

In plain words, we consider those morphisms which invert all $s \in S$.

Now we would like multiplying
$j(s_1)^{-1}j(r_1) \cdot j(s_2)^{-1}j(r_2)$
to obtain again a 'left fraction' $j(s)^{-1}j(r)$. For this  it is
enough to be able to 'commute' the two middle terms
in the sense $j(r_1)j(s_2)^{-1} = j(s')^{-1}j(r')$
as $j(s_1)^{-1}j(s')^{-1}) = j(s' s_1)^{-1}$ and $j(r')j(r_2) = j(r'r_2)$
and we are done. This reasoning is the origin of the left Ore condition.
Here is a formal statement:

\ppt \ldef{SRclosed}{\bf Proposition.} {\it (i)
For $(j,Y) \in {\rm Ob}\in{\cal C}_\catmonoids(R,S)$,
$\,j(S)$ is left Ore in $j(R)$ iff
\[ j(S)^{-1}j(R) = \{ j(s)^{-1}j(r)\,|\,s\in S, r\in R\} \subset Y\]
is a submonoid of $Y$. In particular, if $S$ is left Ore in $R$,
$j(S)^{-1}j(R)$ is a submonoid of $Y$ for each
$(j,Y) \in {\rm Ob}({\cal C}(R,S))$.
\newline (ii) If the equivalent conditions in (i) hold, then
}
\begin{equation}\label{eq:cat_left_denom}\begin{array}{c}
\forall (s,r) \in S \times R
\\
j(s)^{-1}j(r) = j(s')^{-1}j(r')
\end{array}
\Leftrightarrow \left\lbrace \begin{array}{l}
\exists \tilde s\in S,\,  \exists \tilde{r}\in R,\\
j(\tilde{s})j(s') = j(\tilde{r})j(s)\\
j(\tilde{s})j(r') = j(\tilde{r})j(r).
\end{array}\right. \end{equation}

{\it Proof.} (i) ($\Rightarrow$) Let $s_1,s_2 \in S$ $r_1,r_2 \in R$.
By the left Ore condition $\exists s' \in S$ $\exists r'\in R$
with $j(s')j(r_1) = j(r')j(s_2)$.
Hence the product
$j(s_1)^{-1}j(r_1)\cdot j(s_2)^{-1}j(r_2)= j(s's_1)^{-1}j(r'r_2)$ belongs to
$Y$.
\newline ($\Leftarrow$) If $j(S)^{-1}j(R)$ is a monoid then
$j(r)j(s)^{-1}\in  j(S)^{-1}j(R)$. In other words, $\exists s' \in
S$ $\exists r' \in R$ such that $j(r)j(s)^{-1}= j(s')^{-1}j(r')$.
Thus $j(s_1)j(r) = j(r_1)j(s)$.

 (ii) By multiplying from the left by $j(s')$
one gets $j(s')j(s^{-1})j(r) = j(r')$. As $S$ is left Ore,
$\exists \tilde{s}\in S$ $\exists \tilde{r} \in R$ such that
$\tilde{s} s' = \tilde{r} s$. This implies $j(\tilde{s})j(s') =
j(\tilde{r})j(s)$ and, consequently, $j(s')j(s^{-1}) = j(\tilde{s})^{-1}
j(\tilde{r})$; then $j(\tilde{s})^{-1}j(\tilde{r})j(r) = j(r')$
and, finally, $j(\tilde{r})j(r) = j(\tilde{s})j(r')$.\vskip .02in

\ppt {\bf Proposition.} \ldef{distributivenessOre} {\it Let
$S,R,Y,j$ be as in~\luse{SRclosed}, and let $R, Y$ be each
equipped with a binary operation, in both cases denoted by $+_0$,
such that $\cdot$ is left distributive with respect to $+_0$. If
$j(S)^{-1}j(R)$ is a submonoid of $Y$, then it is closed with
respect to $+_0$ as well. }

{\it Proof.} The following calculation is valid in $Y$:
\[\begin{array}{l} j(s_1)^{-1} j(r_1) +_0 j(s_2)^{-1}j(r_2) =
j(s_1)^{-1} ( j(r_1) +_0  j(s_1)j(s_2)^{-1}j(r_2)) \\
\,\,\,\,\,\,\,\,\,\,\,\,\,\,\,\,\,\,\,\,\,\,\,\,\,\,\,\,\,\,\,\,\,\,\,\,\,\,
= j(s_1)^{-1} ( j(r_1) +_0  j(\tilde{s})^{-1} j(\tilde{r})j(r_2))\\
\,\,\,\,\,\,\,\,\,\,\,\,\,\,\,\,\,\,\,\,\,\,\,\,\,\,\,\,\,\,\,\,\,\,\,\,\,\,
= j(\tilde{s} s_1)^{-1}(j(\tilde{s})j(r_1) +_0 j(\tilde{r})j(r_2))
\in j(S)^{-1}j(R),\end{array}\]
where $j(\tilde{s})j(s_1)=j(\tilde{r})j(s_2)$ for some $\tilde{s}$,
$\tilde{r}$ by the left Ore condition which holds due~\luse{SRclosed}.

\ppt {\bf Remark.} We do not require $j(a +_0 b) = j(a) +_0 j(b)$ here.

\ppta {\bf Exercise.} Generalize this to a family ${\cal F}$ of
$n$-ary left distributive operations in place of $+_0$, i.e., of
operations of the form $F : X^{\times n}\to X$, for various $n
\geq 0$, such that $L_a \circ F = F \circ L_a^{\times n}$.

\ppt From now on we limit to the case where the category ${\cal
D}$ above corresponds to a variety ${\frak D}$ of algebras $(A,
{\cal L}_A)$ (in the sense of universal algebra) of signature
${\cal L} = (\cdot, 1, {\cal F})$, where ${\cal F}_A$ is a family
of left distributive operations on $A$ on a $(A,\cdot_A,1_A)$. The
reader who does not care for this generality (suitable say for
algebras with operators) can consider 3 basic cases: 1) ${\cal F}
= \emptyset$ when ${\cal D} = \catmonoids$; 2) ${\cal F} = \{+\}$
and algebras are unital rings; 3) ${\cal F} = \{+\}$ and algebras
are associative unital $\genfd$-algebras over a commutative ring
$\genfd$.

\ppt Denote by ${\cal C}_{l,{\cal D}}(R,S)$ the full subcategory
of ${\cal C}_{\cal D}(R,S) $ consisting of those objects $(j, Y)$
which satisfy
\begin{itemize}
\item $(\forall r,r' \in R)\, (j(r) = j(r')
\Leftrightarrow \exists s \in S\,(sr = sr'))$.
\item $j(S)^{-1} j(R)$ is a subring of $Y$
\end{itemize}
Similarly, ${\cal C}_{l,{\cal D}}^-(R,S)$ by definition consists
of objects satisfying the first, but not necessarily the second
property. Denote by ${\cal C}_{r,{\cal D}}(R,S)$ the full
subcategory of ${\cal C}_{\cal D}(R,S) $ consisting of those
objects $(j, Y)$ which satisfy the symmetric conditions
\begin{itemize}
\item $(\forall r,r' \in R)\,(j(r) = j(r')
\Leftrightarrow \exists s \in S\,(rs = r's) )$.
\item $ j(R)j(S)^{-1}$ is a subring of $Y$
\end{itemize}
Finally, the objects in ${\cal C}_{r,{\cal D}}^-(R,S)$, by
definition, satisfy the first, but not necessarily the second
property. If there is a universal initial object in ${\cal
C}_l(R,S)$ (${\cal C}_r(R,S)$),
we denote it by $(\iota, S^{-1}R)$ (resp.$(\iota, RS^{-1})$ )
and we call the pair, or by abuse of language,
also the ring $S^{-1}R$ ($RS^{-1}$ resp.), the left (right)
{\bf Ore localization} \index{Ore localization}
of $R$ at set $S$, and map $\iota$ the localization map.
An alternative name for $S^{-1}R$ ($RS^{-1}$)
is the left (right) {\bf ring of fractions}~\index{ring of fractions}
(of ring $R$ at set $S$).\vskip .02in

\ppt {\bf Proposition.} {\it
If $\forall (j,Y)$ in ${\cal C}_{l,\catmonoids}^-(R,S)$
the subset $j(S)^{-1} j(R)$ is a submonoid
(i.e.$\exists (j,Y) \in  {\cal C}_{l,\catmonoids}(R,S)$),
then it is so $\forall (j,Y)$ in ${\cal C}_{l,\catmonoids}^-(R,S)$, i.e.
the categories ${\cal C}_{l,\catmonoids}(R,S)$ and
${\cal C}_{l,\catmonoids}^-(R,S)$
coincide. In that case, $S$ is a left denominator set in $R$.
}

{\it Proof.} Let $j(S)^{-1} j(R)$ be a subring for some $(j,Y)$.
Then $j(S)^{-1}j(R)$ is Ore in $j(R)$ by the previous proposition.
Thus for every $s\in S$, $r \in R$ $\exists s' \in S$ $\exists r'
\in R$ such that $j(r)j(s)^{-1} = j(s')^{-1} j(r')$ and therefore
$j(s')j(r)= j(r')j(s)$ which means $j(s'r) = j(r' s)$. That implies
$\exists s^+ \in S$ with $s^+ s'r = s^+r' s)$. Therefore for any
other $(j',Y')$ in ${\cal C}_{l,\catmonoids}(R,S)$ the subset
$j'(S)^{-1} j'(R)$ is a subring. Moreover we have $s^+ s' \in S$
and $s^+ r'$ satisfy $(s^+ s')r = (s^+ r')s$. Since they were
constructed for an arbitrary $s$ and $r$, $S$ is left Ore in $R$.
\newline Left reversibility: Let $r,r'\in R, s\in S$.
Then $rs = r's \Rightarrow j(r)j(s) = j(r')j(s)$,
so by invertibility of $j(s)$ also $j(r) = j(r')$.
But $(j,Y)$ is object in ${\cal C}_{l,\catmonoids}(R,S)$
so $j(r) = j(r') \Rightarrow \exists s'\in S, s'r = s'r'$.

\ppt \ldef{relFrac}
{\bf Lemma.} {\it (i) Let $S$ be a left denominator set.
Define the relation $\sim$ on $S\times R$ by
\[ (s,r) \sim (s',r') \,\,\Leftrightarrow\,\,
(\exists \tilde s\in S \,\,\exists \tilde{r}\in R) \,\,
(\tilde{s}s' = \tilde{r}s\,\mbox{ and }\,\tilde{s}r' =
\tilde{r}r).\] Then $\sim$ is an equivalence relation.

(ii) Let $\Sigma$ be a system of left fractions in a category ${\cal C}$.
For any pair of objects $X,Y$  in ${\cal C}$
let $(\Sigma \times {\cal C})(X,Y)$ be a class
of all diagrams of the form
$\left(X \stackrel{r}{\rightarrow} Z \stackrel{s}{\leftarrow} Y\right)$
in ${\cal C}$. Define a relation $\sim$ on $(\Sigma \times {\cal C})(X,Y)$ by
\[\begin{array}{l}
\left( X \stackrel{r}{\rightarrow} Z \stackrel{s}{\leftarrow} Y\right)
\sim
\left(X \stackrel{r'}{\rightarrow} Z' \stackrel{s'}{\leftarrow} Y\right)
\\ \,\,\,\,\,\,\,\,\,\,\,\,\,\Leftrightarrow
\exists
\left(X \stackrel{\tilde{r}}{\leftarrow} B \stackrel{\tilde{s}}
{\rightarrow} Y\right),\left\lbrace\begin{array}{l}
\tilde{s}\circ s' = \tilde{r}\circ s : B \rightarrow Z \\
\tilde{s}\circ r' = \tilde{r}\circ r : B \rightarrow Z' \end{array}\right. .
\end{array}\]
The latter condition can be depicted by saying that the diagram
$$\xymatrix{
&Z&\\
X\ar[dr]_{r'}\ar[ur]^{r}&B\ar[r]^{\tilde{r}}\ar[l]_{\tilde{s}}
&Y\ar[dl]^{s'}\ar[ul]_{s}\\
&Z'& }$$ commutes. Conclusion: $\sim$ is an equivalence relation.}

Here $\left(X \stackrel{\tilde{r}}{\leftarrow} B \stackrel{\tilde{s}}
{\rightarrow} Y\right)$ is {\it not} a diagram in
$(\Sigma^{\rm op} \times {\cal C})(X,Y)$.

{\it Proof.} \underline{Reflexivity} is trivial.
\newline\underline{Symmetry}:
By Ore $\exists r_1 \in R, s_1\in S$ with $r_1 s = s_1 s'$.
Also $\exists r_2\in R, s_2\in S$ with $r_1 s = s_1 s'$.
Thus \[ r_2 \tilde{r} s' = r_2 \tilde{s} s = s_2 r_1 s = s_2 s_1 s'.\]
In other words $r_2 \tilde{r} - s_2 s_1 \in I_S$. Thus by the
left reversibility, $\exists t \in S$ with $t(r_2 \tilde{r} - s_2 s_1) = 0$.
Therefore $t(r_2 \tilde{r} - s_2 s_1)r' = 0$,
hence $t s_2 s_1 r' = t r_2 \tilde{s}r$. Compare with
$t s_2 s_1 s' = tr_2 \tilde{s}s$ to see that $(s',r')\sim (s,r)$.
\newline\,\,\,\,\underline{Transitivity}:
Assume $(s,r) \sim (s',r')$ and $(s',r') \sim (s'',r'')$. This means
\[ \exists \,\tilde{s}, \titilde{s} \in S\,
\exists \,\tilde{r}, \titilde{r} \in R \,\,\,\,\,\,\,\,
\left\lbrace
\begin{array}{lr}
\tilde{s} s = \tilde{r} s'\,\,\,\,\,
&\,\,\,\,\,\titilde{s} s' = \titilde{r} s''
\\
\tilde{s} r = \tilde{r} r'\,\,\,\,\,
&\,\,\,\,\,\titilde{s} r' = \titilde{r} r''
\end{array}\right.\]
$S$ is  left Ore, hence
$\exists s_* \in S\,\exists r_* \in R$ with
$s_* \tilde{r} = r_* \titilde{s}$. Therefore
\[\begin{array}{l}
(s_* \tilde{s})s = s_* \tilde{r} s'
= r_* \titilde{s} s' = (r_* \titilde{s}) s''\\
(s_* \tilde{s})r = s_* \tilde{r} r'
= r_* \titilde{s} r' = (r_* \titilde{s}) r''
\end{array}\]
Hence $(s,r) \sim (s'',r'')$.

\ppt \ldef{easyequiv}{\bf Simplifying observation.}
Consider a family of arrows $(s,r) \to ((ps)^{-1}, (pr))$
where $ps \in S$. Then $(s,r) \sim ((ps)^{-1}, (pr))$.
If some property ${\cal P}$ of elements of $S\times R$
does not change along such arrows, then ${\cal P}$
is well-defined on classes $s^{-1}r :=[s,r]/\sim \,\in S\times R/\sim$.

{\it Proof.} Clearly every $\sim$-arrow
is a composition of one such arrow
and an inverse of another such arrow.

\ppt {\bf Lemma.} \ldef{fracOgrlak}
{\it If $t\in S$ and $tr = tr'$ (by reversibility
even better if $rt = r't$) then
$(s, r_1 r r_2) \sim (s, r_1 r' r_2)$.
}

{\it Proof.} There are
$t' \in S$, $r'_1\in R$ with  $r'_1 t = t' r_1$.
Then $(s, r_1 r r_2)\sim (t' s, t' r_1 r r_2)\sim (t' s,r'_1 t r r_2)
\sim (t' s, r'_1 t r' r_2)\sim (t' s, t' r_1 r' r_2)\sim (s, r_1 r' r_2)$.

\ppt {\bf Proposition.} \ldef{criteriumUniv} {\it For
$(j,Y)$ in ${\rm Ob}\,{\cal C}_{l,{\frak D}}(R,S)$ the statement
\begin{equation}\label{eq:LU}
( \forall y \in Y \,\,\,\,\exists s \in S \,\,\,\,\,
\exists r \in R ) \,( y = j(s)^{-1} j(r) )   \end{equation}
holds iff $(j,Y)$ is a universal initial object in
${\cal C} = {\cal C}_{{\frak D}}(R,S)$.
}

{\it Proof.} ($\Leftarrow$) Let $(j,Y)\in {\cal C}_l$ be universal
in ${\cal C}$. Suppose $Y_0 \defeq j(S)^{-1}j(R)$ is a proper
subring of $Y$. We'll denote by  $j_0$ the map from $R$ to $Y_0$
agreeing with $j$ elementwise. Then $(j_0,Y_0)$ is an object in
${\cal C}_l$ and the inclusion $i : Y_0 \rightarrow Y$ is a
morphism from $(j_0,Y_0)$ into $(j,Y)$. By universality of $(j,Y)$
there is a morphism $i' : (j,Y) \rightarrow (j_0,Y_0)$. The
composition of morphisms $i\circ i'$ is an automorphism of $(j,Y)$
clearly different from the identity, contradicting the
universality of $(j,Y)$.

($\Rightarrow$) Let $(j,Y)$ satisfy~(\ref{eq:LU}) and let
$(j',Y')$ be any object in ${\cal C}(R,S)$. We want to prove that
there is unique map $i : Y \rightarrow Y'$ which satisfies
$i(j(r)) = j'(r) \mfall r \in R$. Note that $i(j(s)^{-1}j(s)) =
i(j(s)^{-1})j'(s)$ implies $i(j(s)^{-1}) = j'(s)^{-1}$. Thus
$i(j(s)^{-1}j(r)) = j'(s)^{-1}j'(r)$ so that the value of $i$ is
forced for all elements in $Y$ proving the uniqueness.

This formula sets $i$ independently of choice of $s$ and $r$.
Indeed, if $j(s)^{-1}j(r) = j(s')^{-1}j(r')$
then $j(r) = j(s) j(s')^{-1}j(r')$.
As $j(S)$ is left Ore in $j(R)$, we can find $\tilde{s} \in S$ and
$\tilde{r} \in R$ such that $j(\tilde{r}) j(s')= j(\tilde{s}) j(s)$
and therefore $j(s) j(s')^{-1}= j(\tilde{s})^{-1}j(\tilde{r})$.
Thus $j(r) = j(\tilde{s})^{-1} j(\tilde{r}) j(r')$
or $j(\tilde{s})j(r) = j(\tilde{r}) j(r')$ and, finally,
$j(\tilde{s}r) = j(\tilde{r} r')$.
Thus $\exists s^+ \in S,$ $s^+ \tilde{s}r = s^+\tilde{r}r'$.
Starting here and reversing the chain of arguments, but with
$j'$ instead of $j$, we get $j'(s)^{-1}j'(r) = j'(s')^{-1}j'(r')$.

\ppt \ldef{Oreexistencemon}
{\bf Theorem.} {\it If $S$ is a left denominator set in $R$,
then the universal object $(j,Y)$
in ${\cal C}_{l,\catmonoids}(R,S)$ exists.}

{\it Proof.} We will construct a universal object $(j,Y) \equiv
(\iota, S^{-1}R)$. \underline{As a set}, $S^{-1}R \defeq (S \times
R)/\sim$. Let $[s,r]$, and,  by abuse of notation, let $s^{-1}r$
also denote the $\sim$-equivalence class of a pair $(s,r) \in S
\times R$. Notice that $1^{-1}r =1^{-1}r'$ may hold even for some
$r\neq r'$, namely when $\exists s \in S$ and $r,r'$ with $sr =
sr'$. The equivalence relation is forced
by~(\ref{eq:cat_left_denom}).
\newline\underline{Multiplication} \fbox{
$s_1^{-1}r_1\cdot s_2^{-1}r_2 := (\tilde{s}s_1)^{-1} (\tilde{r}r_2)$} where
$\tilde{r} \in R$, $\tilde{s} \in S$
satisfy $\tilde{r}s_2 = \tilde{s}r_1$
(thus $\tilde{s}^{-1}\tilde{r} = r_1 s_2^{-1}$), as in the diagram:
$$
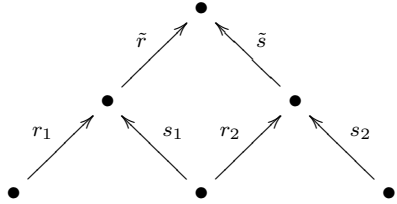
\xymatrix{
&& \bullet &&\\
&\bullet\ar[ur]^{\tilde{r}}&&\bullet\ar[ul]_{\tilde{s}}&\\
\bullet\ar[ur]^{r_1}&&\bullet\ar[ul]_{s_1}\ar[ur]^{r_2}
&&\bullet\ar[ul]_{s_2}}$$
If we choose another pair of multipliers $\titilde{r} \in R,
\titilde{s} \in S$ with $\titilde{r}s_2 = \titilde{s}r_1$ instead,
then by the left Ore condition we can choose $r_* \in R, s_*\in S$
with $s_* \tilde{s} = r_*\titilde{s}$. Then
\[ r_* \titilde{r} s_2 = r_* \titilde{s} r_1
= s_* \tilde{s} r_1 = s_* \tilde{r} s_2 \]
and therefore $r_* \titilde{r} - s_* \tilde{r} \in I_S$.

In other words,
$\exists s^+ \in S$ with $s^+ r_* \titilde{r} = s^+ s_* \tilde{r}$.

Therefore we have
\[\begin{array}{l}
s^+ s_* \tilde{r} r_2 = s^+ r_* \titilde{r} r_2 \\
s^+ s_* \tilde{s} s_1 = s^+ r_* \titilde{s} s_1
\end{array}\]
which proves $(\tilde{s}s_1)^{-1} (\tilde{r} r_2) =
(\titilde{s}s_1)^{-1} (\titilde{r} r_2)$. Thus multiplication is
well defined as a map $\mu_1 : (S\times R) \times (S\times R)\to S^{-1}R$.

We have to show that $\mu_1$ factors to
$\mu :  S^{-1}R\times S^{-1}R\to S^{-1}R$.

By~\luse{easyequiv}, it is sufficient to show that
$a = \mu_1((s_1, r_1),(s_2, r_2))$ equals
$b = \mu_1(((rs)_1,(rr_1)),((ps)_2, (pr_2)))$ whenever $rs\in S$
and $ps\in S$.



$s_2'r_1 = r_1' s_2$ for some $s_2'\in S$ and $r_1' \in R$.
Then $a = (s_2' s_1)^{-1} (r_1' r_2)$.
As $ps_2\in S$, $\exists p' \in R, s_*\in S$ with
$p' (ps_2) = s_* r_1' s_2 = s_* s_2' r_1$.
Furthermore, $s_\sharp r = p_\sharp s_* s_2'$ for some
$s_\sharp \in S$ and $p_\sharp \in R$.
Putting these together, we infer
$s_\sharp (r r_1) = p_\sharp s_* s_2' r_1 = p_\sharp p' ps_2$
and therefore $(rr_1)(ps_2)^{-1} \to s_\sharp^{-1} (p_\sharp p')$,
i.e., be definition, that
$b = (s_\sharp rs_1)^{-1}(p_\sharp p' p r_2)$, hence by above,
$b = (p_\sharp s_* s_2' s_1)^{-1}(p_\sharp p' p r_2)$.
Now use lemma \luse{fracOgrlak} and $(p'p)s_2 = (s_* r_1') s_2$ to conclude
$b = (p_\sharp s_* s_2' s_1)^{-1}(p_\sharp s_* r_1' r_2)=
(s_2' s_1)^{-1}(r_1' r_2) = a$.

Hence $\mu$ is well-defined. The unit element is
clearly $1 = 1^{-1} 1$. We need to show associativity of $\mu$.
The product
$s_1^{-1}r_1 \cdot s_2^{-1}r_2 \cdot s_3^{-1}r_3$
does not depend on the bracketing, essentially because one can complete
the following commutative diagram of elements in $R$
(the composition is the multiplication in $R$:
any pair of straight-line (composed) arrows
with the same target is identified with a pair in $S\times R$):
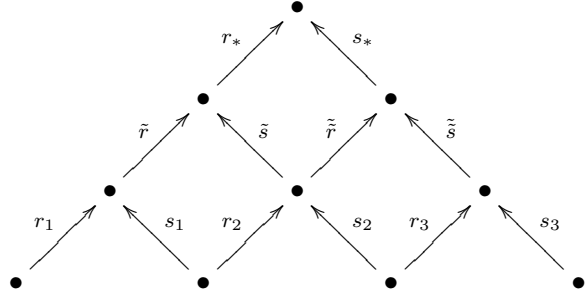
$$\xymatrix{
&&&\bullet&&&\\
&& \bullet\ar[ur]^{r_*} &&\bullet\ar[ul]_{s_*}&&\\
&\bullet\ar[ur]^{\tilde{r}}&&\bullet\ar[ul]_{\tilde{s}}\ar[ur]^{\titilde{r}}
&&\bullet\ar[ul]_{\titilde{s}}&\\
\bullet\ar[ur]^{r_1}&&\bullet\ar[ul]_{s_1}\ar[ur]^{r_2}&&\bullet\ar[ul]_{s_2}
\ar[ur]^{r_3}&&\bullet\ar[ul]_{s_3}
}$$

Finally, the construction gives the universal object because it
clearly satisfies the equivalent condition
in~\luse{criteriumUniv}.

\minsection{Ore localization for rings}

\indent \ppt {\bf Exercise.} The two left Ore conditions
together immediately imply the {\bf combined left Ore condition}:

{\it If $n\in R$ is such that $n s = 0$ for some $s \in S$, then
for every $r\in R$ there are $s'\in R$, $r'\in R$ such that
$r' n = s' r$.}

It is sometimes useful to quote this property in order to
avoid introducing additional variables needed for deriving it.

\ppt {\bf Lemma.} {\it Let $f: R \rightarrow R'$ be a ring morphism and
$S \subset R$ Ore. Then $f(S)$ is an Ore set in R.}

\ppt {\bf Notation.}  In this section we are concerned
only with the category of unital rings.
Thus ${\cal C}(R,S):= {\cal C}_\catRings(R,S)$.

\ppt {\bf Notation.} For any $S \subset R$ let \ldef{IS}\index{$I_S$}
$I_S := \{ n \in R \,|\,\exists s \in S,\, sn = 0\}$.
$I_S$ is clearly a right ideal. If $S$ is a left Ore set,
then $sn = 0$ and the left Ore condition imply that $\forall r \in R$
$\exists s_0 \in S,$ $r_0\in R$ with $r_0 s = s_0 r$, hence
$s_0 rn = r_0 sn = 0$. Thus $I_S$ is then a 2-sided \underline{ideal}.

\ppt {\bf Corollary.} If $S^{-1}R$ exist then
$\forall (j,Y) \in {\cal C}(R,S)$, ${\rm ker}\, j \subset I_S$.
In particular, an Ore localization of a domain is a domain.

\ppt \ldef{Oreexistence}
{\bf Theorem.} {\it If $S$ {\it is a left denominator set} in $R$
then the universal object $(j,Y)$ in ${\cal C}_l(R,S)$ exists.
}

{\it Proof.} In~\luse{Oreexistencemon}, we have constructed a monoid
structure on $Y = S\times R/\sim$. We exhibit an additive
structure on $Y$ such that $j$ is a ring map and
$(j,Y)$ is indeed universal.
\newline\underline{Addition}: Suppose we are given two fractions
with representatives $(s_1,r_1)$ and $(s_2,r_2)$. By the left Ore
condition, $\exists\tilde{s} \in S$, $\exists \tilde{r}\in R$, such that
$\tilde{s} s_1 = \tilde{r} s_2$. The sum is then defined as
\[\fbox{$ s_1^{-1} r_1 + s_2^{-1} r_2 \defeq
 (\tilde{s}s_1)^{-1} (\tilde{s}r_1 + \tilde{r}r_2) $}\]
Suppose we have chosen $(\titilde{s},\titilde{r}') \in S \times R$ with
$\titilde{s} s_1 = \titilde{r} s_2$, instead of $(\tilde{s}, \tilde{r})$.
Then by left Ore we find $s_* \in S$
and $r_* \in R$ such that $s_* \tilde{s} = r_* \titilde{s}$.
Then
\[ r_* \titilde{r} s_2 = r_* \titilde{s} s_1 = s_* \tilde{s} s_1 =
s_* \tilde{r} s_2 \]
hence $(s_* \tilde{r} - r_* \titilde{s}) \in I_S$, i.e.
$\exists s^\sharp \in S$ with
\[ s^\sharp s_* \tilde{r} = s^\sharp  r_* \titilde{s} \]
Then
\[\begin{array}{l}
(s^\sharp s_*) (\tilde{s} r_1 + \tilde{r}r_2) =
        (s^\sharp r_*) (\tilde{s} r_1 +\titilde{s} r_2)\\
(s^\sharp s_*) (\tilde{s} s_1) = (s^\sharp r_*) (\titilde{s} s_2)
\end{array}\]
Conclusion: $(\tilde{s} s, \tilde{s} r_1 + \tilde{r}r_2)
\sim (\titilde{s} s, \titilde{s} r_1 + \titilde{r}r_2)$, as required.

Now let's check that the sum does not depend on the choice of
the representative of the first summand.
Suppose we are given two representatives of
the first fraction $s^{-1}_1 r_1 = s'^{-1}_1 r'_1$.
Then for some $({s}^*,{r}^*) \in S\times R$ we have
\[ s_* s_1 = r_* s'_1\,\,\,\,\mbox{ and }\,\,\,\,s_* r_1 = r_* r'_1\]
Second fraction in $s_2^{-1} r_2$. Choose
\[(\titilde{s},\titilde{r}) \in S\times R\,\,\,\,\,\mbox{ with }\,\,\,\,
 \titilde{s} s'_1 = \titilde{r} s_2. \]
Now choose $(s_\sharp, r_\sharp)\in S\times R$ such that
$s_\sharp r_* = r_\sharp \titilde{s}$. Then
$(r_\sharp \titilde{r}) s_2 = r_\sharp \titilde{s} s'_1 =
s_\sharp r_* s'_1 = (s_\sharp s_*) s_1$
and $(r_\sharp \titilde{s}) r'_1 = s_\sharp r_* r'_1 = (s_\sharp s_*) r_1$.
Therefore
\[ \begin{array}{lcl}
s^{-1}_1 r_1 + s^{-1}_2 r_2 & = & (s_\sharp s_* s_1)^{-1} (s_\sharp s_* r_1
+ r_\sharp \titilde{r} r_2 ) \\
&=& (r_\sharp \titilde{s} s'_1)^{-1} (r_\sharp \titilde{s} r_1'
+ r_\sharp \titilde{r} r_2 ) \\
&=& (\titilde{s} s'_1)^{-1} (\titilde{s} r'_1 + \titilde{r} r_2) \\
&=&  s'^{-1}_1 r'_1 + s^{-1}_2 r_2
\end{array}\]
We should also check that the sum does not depend on the second
summand. This proof would not be symmetric to this one as our
definition of the sum is not. We shall choose an indirect proof.
Denote the set-theoretic quotient map by $p : S \times R
\rightarrow S^{-1}R$. By now we have completed the proof that
addition as a map from
\[ \tilde{+} : S^{-1}R \times (S \times R) \rightarrow  S^{-1}R \]
is well defined. Now we prove that the map
\[ \tilde{+}(p \times {\rm id})\tau : (S \times R) \times
(S \times R)\rightarrow S^{-1}R \] where $\tau$ is the
transposition of factors coincides with $\tilde{+}(p \times {\rm
id})$. Thus we have a well-defined addition as a map defined on
$S^{-1}R \times S^{-1}R$ which is then automatically commutative.
It is sufficient to prove that for any two pairs $(s_1, r_1)$ and
$(s_2,r_2)$ and any
\[ \tilde{s}, \titilde{s} \in S, \tilde{r}, \titilde{r} \in R\,\,
\mbox{ with }\,\,\,\tilde{s} s_1 = \tilde{r} s_2, \,\,\,\,\,\,\,
\titilde{r} s_1 = \titilde{s} s_2, \]
the classes
\[\begin{array}{l}
(\tilde{s} s_1)^{-1} (\tilde{s} r_1 + \tilde{r} r_2 ) \\
(\titilde{s} s_2)^{-1} (\titilde{r} r_1 + \titilde{s} r_2 )
\end{array}\]
coincide in $S^{-1}R$. For that purpose, choose $s_\sharp \in S$
and $r_\sharp \in R$ such that $s_\sharp \tilde{r} = r_\sharp \titilde{s}$.
Then
\[ r_\sharp \tilde{s} s_1 = s_\sharp \tilde{r} s_2
= r_\sharp \titilde{s} s_2. \]
Next $r_\sharp \titilde{r} s_1 = r_\sharp \titilde{s} s_2 =
s_\sharp \tilde{r} s_2 = s_\sharp \tilde{s} s_1$, and therefore
$(r_\sharp \titilde{r} - s_\sharp \tilde{s}) \in I_S$ (\luse{IS}). Thus
$\exists s^+ \in S$ with
\[ s^+ r_\sharp \titilde{r} - s^+ s_\sharp \tilde{s} = 0. \]
In particular,
 $s^+ r_\sharp \titilde{r} r_1 = s^+ s_\sharp \tilde{s} r_1 = 0$.
Thus
\[\begin{array}{lcl}
(\tilde{s} s_1)^{-1} (\tilde{s} r_1 + \tilde{r} r_2 ) & = &
(s^+ s_\sharp \tilde{s} s_1)^{-1} (s^+ s_\sharp\tilde{s} r_1 +
s^+ s_\sharp \tilde{r} r_2 )\\
&=& (s^+ r_\sharp \titilde{s} s_1)^{-1} (s^+ r_\sharp\titilde{r} r_1 +
s^+ r_\sharp \titilde{s} r_2 )\\
&=& (\titilde{s} s_1)^{-1} (\titilde{r} r_1 +  \titilde{s} r_2 )
\end{array}\]
The associativity of addition is left to the reader.

The distributivity law follows by~\luse{distributivenessOre}.

The element $1^{-1} 0$ in $S^{-1}R$ is the zero and thus $S^{-1}R$
is equipped with a natural unital ring structure.

Define $\iota : R \rightarrow S^{-1}R$ by $\iota(r) = [1,r] = 1^{-1}r$.
Check that $\iota$ is a unital ring homomorphism.
Check that $\iota(S)$ consists of units and that
$\iota(S)^{-1}\iota(R) = Y$.
Pair $(\iota,S^{-1}R)$ is a universal object in ${\cal C}_l(R,S)$,
as it clearly satisfies the equivalent condition in~\luse{criteriumUniv}.

\ppt {\it Right} Ore conditions, and right Ore localizations with respect to
$S\subset R$, are by definition the left Ore conditions and
localizations with respect to $S \subset R^{\rm op}$.
The {\it right} ring of fractions is denoted
$RS^{-1} := (S^{-1}R^{\rm op})^{\rm op}$.
It consists of certain equivalence pairs $rs^{-1} :=
[(r,s)]$, where $(r,s)  \in R \times S$.

\minsection{Practical criteria for Ore sets}

This section is to be read only by those who want to test in practice
wheather they have an Ore set at hand.

\ppt {\bf Theorem.} {\it
(i) Let $S$ and $S'$ be multiplicative sets in ring $R$,
where $S$ is also left Ore in $R$. Assume

1. for a map $j : R \rightarrow Y$ of unital rings,
the image $j(S)$ consists of units in $Y$ iff
the image $j(S')$ consists of units in $Y$;

2. $sr = 0$ for some $s \in S$ iff $\exists s' \in S'$ with $s' r = 0$.

Then $S'$ is left Ore as well and
$S^{-1}R$ is canonically isomorphic to $S'^{-1}R$.
}
{\it Proof.} Under the assumptions the categories ${\cal
C}_l(R,S)$ and ${\cal C}_l(R,S')$ are identical, so call them
simply ${\cal C}$. The left Ore condition is equivalent to the
existence of an initial object in ${\cal C}$; and the 2
localizations are just the 2 choices of an initial object, hence
there is a unique isomorphism in ${\cal C}$ between them; its
image under the forgetful functor ${\cal C} \rightarrow R-{\rm
Mod}$ into the category of unital rings, is the canonical
isomorphism as required. \vskip .017in

\ppt The left Ore condition is often checked
inductively on a filtration, or ordered set of generators.
To disseminate this kind of reasoning we will temporarily use
some {\it nonstandard} notation which generalizes
the left Ore condition.
One fixes a (only) multiplicative set $S \subset R$.
For {\it any} subset $A \subset R$,
and any $(s,r)\in S \times A \subset S \times R$,
introduce predicate
\[ {\bf lOre}(s,r\uparrow A) :=
{\bf lOre}_{S,R}(s,r\uparrow A) \equiv
(\exists s' \in S, \,\exists r' \in A, \,s'r = r's),
\]
where $S,R$ in subscripts may be skipped if known from context.
Moreover if $A = R$ then $\uparrow A$ may be skipped from the notation.
For example, ${\bf lOre}(s,r) = {\bf lOre}_{S,R}(s, r \uparrow R)$.

For any subsets $A_0 \subset A$ and $S_0 \subset S$, abbreviate
\[ {\bf lOre}(S_0,A_0\uparrow A) \equiv  \left(\,\forall s \in S_0,
\forall r \in R_0,\,{\bf lOre}_{S,R}(s,r\uparrow A)\,\right),\]
with rules for skipping $\uparrow A$ as before.
For example, ${\bf lOre}(S,R)$ is simply the
left Ore condition for $S \subset R$.

Finally,
\[
{\bf slOre}(S_0,A) \equiv  {\bf lOre}_{S,R}(S_0, A \uparrow A).
\]
For an additive subgroup
$A \subset R$ consider also the relative versions, e.g.
\[
{\bf rel-lOre}_{S,R}(s,r; I) \equiv
(\exists s' \in S, \,\exists r' \in R, \,s'r - r's \in A).
\]
If $A = I$ is an ideal, then this predicate is suitable for study
of some (non-Ore) generalizations (cf.~\cite{Elizarov} for such).

\ppt \ldef{extendOre}{\bf Extending Ore property}
 \index{extending Ore property}
{\it Let $A,B\subset R$ be additive subgroups of $R$,
$A \subset B\subset R$,
and $S\subset R$ multiplicative subset. }
\[\begin{array}{lc}
(i) & ({\bf lOre}(S,A) \mbox{ and }
{\bf rel-lOre}(S,R; A)) \Rightarrow {\bf lOre}(S,R)\\
(ii) & ({\bf lOre}(S,A\uparrow B) \mbox{ and }
{\bf rel-slOre}(S,B; A)) \Rightarrow {\bf lOre}(S,B)\\
(iii) & ({\bf lOre}(S,A\uparrow B) \mbox{ and }
{\bf rel-slOre}(S,B; A) \\
& \,\,\,\,\,\,\,\,\,\,\mbox{ and } SB\subset B)
\Rightarrow {\bf slOre}(S,B)
\end{array}\]
{\it Proof.} (i) is clearly the $B = R$ case of (ii).
Let $b \in B$ and $s \in S$. Then ${\bf rel-slOre}(S,B; A))$ means
that $\exists s'\in S$, $\exists b' \in B$, $\exists a\in A$ such
that $s' b - b' s = a$. Now we compare $s$ and $a$.
There are $b_1 \in B$, $s_1\in S$ such that
$b_1 s = s_1 a$. Thus $s_1 s' b - s_1 b' s = s_1 a = b_1 s$,
and finally, $(s_1 s')b = (s_1 b' + b_1)s$. $S$ is multiplicative
hence (ii), and if $SB \subset R$ then
$s_1 b' + b_1\in B$ hence (iii).

\ppta {\bf Remark.} The above condition is usually checked for
generators only. Also we can iterate the above criterion if we
have a finite or denumerable family of nested subrings, for which
the induction is convenient. One may also need to nest subsets of
$S$, with refined criteria, like ${\bf lOre}(S_1 \uparrow
S_2,A_1\uparrow A_2)$, where the $\uparrow S_2$ means that $s'$
may be chosen in $S_2$.

\ppt \ldef{S1multOre}
{\bf Lemma.} If $S_1$ {\it multiplicatively} generates
$S$, and $A \subset R$ then
\[ {\bf lOre}(S_1, R)
\Leftrightarrow {\bf lOre}(S,R),\]
\[ {\bf slOre}(S_1,A) \Leftrightarrow {\bf slOre}(S,A).\]

{\it Proof.} The first statement is clearly a particular case of
the second. Hence we prove the second statement; the nontrivial
direction is $\Rightarrow$. By assumption, the set $S$ can be
written as a nested union $\cup_{n\geq 0} S_n$ where $S_n$
consists of all those $s \in S$ which can be expressed as a
product $\prod_{k = 1}^{n'} s_k$ with $n' \leq n$ and $s_k \in
S_k$; (hence $S_1$ is as the same as before). The assumption is
${\bf slOre}(S_1,A)$, hence by induction it is enough to prove
that ${\bf slOre}(S_n,A)\Rightarrow{\bf slOre}(S_{n+1},A)$ for all
$n \geq 1$. Take $s = s_1 s_2 \cdots s_n$. Then ${\bf
slOre}(S_n,A)$ means that for any $a\in A$ we have
\[ \begin{array}{l} \exists a' \in A \,\exists s' \in S\mbox{ }
( a's_{2}\cdots s_{n} = s'a ),\\
\exists a'' \in A \,\exists s'' \in S \mbox{ }
( a'' s_{1} = s'' a' )\end{array}\]
and consequently
\[
a'' s_1 s_2 \ldots s_n = s'' a' s_2 \ldots s_n = (s'' s') a,
\]
with desired conclusion by the multiplicative closedness of $S$.

\ppt {\bf Lemma.} \ldef{oreadit}
If $A^+_0, A^+\subset R$ be the
additive closure of $A_0,A$ respectively, then (obviously)
\[ {\bf lOre}(S,A_0 \uparrow A) \Rightarrow
{\bf lOre}(S,A_0^+ \uparrow A^+).\]

\ppt {\bf Lemma.} If $A$ generates $R$ as a ring, then
\[ {\bf lOre}(S,A) \Rightarrow {\bf lOre}(S,R).\]

{\it Proof.} By \luse{oreadit} it is enough to check this
multiplicativity:
\[(\forall i,\,
{\bf lOre}(S,c_i))\Rightarrow {\bf lOre} (S,c_n\cdots c_1),\] when
$c_i \in A$. However, this general statement holds for any choice
of $c_i$ whatsoever. Namely, if we do not require $c_i \in A$, we
see, by induction on $n$, that it is enough to prove this
statement for $n = 2$. For $s \in S$ and $c_1, c_2 \in R$ we can
find $r'_1, r'_2 \in R$ and $s', s'' \in S$, so that $r_1' s = s'
c_1$ and $r_2' s' = s'' c_2$. The result is
\[ r_2' r_1' s = r_2' s' c_1 = s'' c_2 c_1,\]
hence the lemma is proved.

\ppt {\bf Theorem.} {\it
If $A \subset B$ and the subring $\bar{A} \subset R$
is also contained in $B$, then for any $S_1$ multiplicatively
generating a multiplicatively closed set $S$ we have }
\[ {\bf slOre}(S_1,A)\Rightarrow {\bf lOre}(S,\bar{A}\uparrow B),\]
\[ {\bf lOre}(S, A) \Rightarrow {\bf lOre}(S,\bar{A}).\]

{\it Proof.} We know ${\bf slOre}(S_1,A)\Rightarrow {\bf slOre}(S,A)$.
Hence the first assertion follows from the second by $\bar{A}\subset B$.
We proved the second statement in the case $\bar{A} = R$.
If $S\subset \bar{A}$, the statement
clearly does not say anything more than it would
say after replacing $R$ by its subring $\bar{A}$.
The proof of the general case
is exactly the same, as $s \in R$ is never used,
and our calculations and quantifiers
may be taken over a bigger overring.

\ppt {\bf Warning-theorem.} {\it If $A$ generates $R$ as a ring
and $S_1$ generates $S$ multiplicatively,
then it is NOT necessarily true that
\begin{equation}\label{eq:lOrenottrue}
{\bf lOre}(S_1, A)\Rightarrow {\bf lOre}(S,R),
\end{equation}
even if $S_1$ has only one multiplicative generator. } We know
from \luse{S1multOre} that~(\ref{eq:lOrenottrue}) holds if we
replace {\bf lOre} by the stronger condition {\bf slOre}.
Nevertheless, various intermediate conditions, standing between
{\bf lOre} and {\bf slOre}, often utilizing filtrations and
combined arguments, are widely used in practice. However it is
also common to use~(\ref{eq:lOrenottrue}) without proper
justification.

\ppta {\bf Counterexample.} (proving the warning statement above)

Let $R$ be the unital ring generated by 4 generators $z_1, z_2, z_3, D$
modulo the following relations:
\begin{equation}\label{eq:counterOre}
\begin{array}{l}
   D z_1 = z_2 z_3 D,\\
   D^2 z_2 = z_3 z_1 D, \\
   D^3 z_3 = z_1 z_2 D.
\end{array}
\end{equation}
Clearly the~(\ref{eq:counterOre}) are simply the identities
needed to check ${\bf lOre}(S_1, A)$ where $S_1 = \{ D\}$,
and $A = \{ z_1, z_2, z_3\}$. The powers of $D$ on the left-hand side,
which are $1,2,3$ may be replaced by $1,p,q$ respectively, where
$p>0$ and $q >1$ are any integers, and the same proof applies,
but the inequality $q>1$ is indeed essential:
if $p = q = 1$, this is not a counterexample at all.

We claim that for any nonnegative integer $n$,
$D^n z_1 = P D^2$ does not have a solution for $P \in R$, hence
the left Ore condition is {\it not} satisfied for the multiplicative
set $S = \{ 1, D, D^2, \ldots\}$. The proof of the claim will
be by contradiction, but we need first to study a convenient
basis of ring $R$.

A basis of $R$ consists of all ordered monomials where, the
right-hand side of any of the equations~(\ref{eq:counterOre}) does
not appear as a factor. This is obtained using {\sc Shirshov-Bergman}'s diamond
lemma~(\cite{Bergman,BokutKol}) with the reduction system $K$ having 3
reductions corresponding to the relations ~(\ref{eq:counterOre})
with production arrows from right to left. This reduction system
has clearly no ambiguities whatsoever and all reductions send
monomials into monomials in generators $z_1, z_2, z_3, D$. It is
direct to see using this basis that $D$ is not a zero divisor.

Suppose then, that $S$ is left Ore. Then there exist $n$ such that
\begin{equation}\label{eq:Orelevel2}
 D^n z_1 = P D^2, \,\,\,\,\mbox{ for some } P \in R.
\end{equation}
We suppose that $n \geq 3$, and leave the remaining case to the reader.
Equation~(\ref{eq:Orelevel2}) implies $D^{n-1} z_2 z_3 D = PD^2$.
$D$ is not a divisor, hence $D^{n-1} z_2 z_3 = PD$.
Now write $P$ as a sum of linearly independent $K$-reduced monomials $P_i$.
Because $D$ is not a zero divisor, monomials
$P_i D$ are also linearly independent.
Since the reductions in $K$ send monomials
to monomials, and $D^{n-1}z_2 z_3$ is
a $K$-reduced monomial, we see that $D^{n-1}z_2 z_3$ can not
be obtained as a sum of more than one of the linearly independent
monomials $P_i D$, hence $P$ has to be a monomial.
The only way that $PD$ in $K$-reduced form (which is $D^{n-1}z_2 z_3$)
has $z_3$ as the most right-hand side factor is that $P = P' z_1 z_2$
for some $P'$ in $K$-reduced form.
Hence we obtain $D^{n-1} z_2 z_3 = P' z_1 z_2 D = P' D z_3$.
Again using basis one can check directly that $Qz_3 = 0$ implies $Q = 0$,
hence $D^{n-1} z_2 = P' D$.
Now $D^{n-3} D^2 z_2 = D^{n-3} z_3 z_1 D = P' D$
implies $D^{n-3} z_3 z_1 = P'$.
This substituted back in the expression for $P$ and
the equation~(\ref{eq:Orelevel2}) gives
\[ D^n z_1 = P' z_1 z_2 D =
D^{n-3} z_3 z_1 z_1 z_2 D^2 = D^{n-1} z_2 D^2 z_3 D.\] This is a
contradiction as the two sides differ even though they are
$K$-reduced.

\ppt {\bf Proposition.} {\it Let $S$ and $T$ be left Ore sets in
some ring $R$. Then the set of all elements of the form
$st$ where $s \in S$ and $t \in T$ satisfies the left Ore condition in $R$
(but it is not necessarily multiplicatively closed).
}

{Proof.} Suppose the contrary, i.e. there is $st \in ST$ and $r \in R$,
such that we can not find $s'\in S, t'\in T$ and $r' \in B$,
so that $r'st = s't'r$. Set $T$ is Ore so there are
$t' \in T$ and $r_1 \in R$ with $r_1 t = t' r$.
Next we can find $s' \in S$ and $r_2 \in R$ so that
$r_2 s = s' r_1$. Combining, we obtain
\[ s't'r = s'r_1 t = r_2 s t,\]
with contradiction.

\ppt {\bf Proposition.} \ldef{l1tensorS}
{\it Let $S$ be a left Ore set in a $k$-algebra $R$.
The set $1 \otimes S$ of all elements of $R \otimes_k R$
of the form $1 \otimes s$ where $s \in S$ satisfies the left Ore condition.}

{\it Proof.} $1 \otimes S$ is obviously multiplicatively closed.

If the Ore condition is not true, there is an element
$y = \sum_{i = 1}^n a_i \otimes b_i$ and an element $s \in S$ such that
$(1 \otimes S) y \cap (R \otimes R)(1 \otimes s) = \emptyset$.

We use induction by $n$ to find an element in the intersection.
 If $n = 1$ we simply use that $S$ is left Ore to
find $r' \in R$ and $s' \in S$
such that $r' s = s' b_1$ and we have
\[ (1 \otimes s')(a_1 \otimes b_1) = (a_1 \otimes r')(1 \otimes s),\]
which proves the basis of induction.

Suppose we found $s'_j \in S$ and
$z = \sum_{i = 1}^j a_i \otimes r'^j_i$ so that
\[ (1 \otimes s'_j)(\sum_{i = 1}^j a_i \otimes b_i) =
(\sum_{i=1}^j a_i \otimes r'_i)(1 \otimes s).\]
Now we use again the property that $S$ is left Ore to find
$r'^{j+1}_{j+1} \in R$ and $\bar{s}_{j+1} \in S$ such that
\[ r'^{j+1}_{j+1} s = \bar{s}_{j+1} s'_j b_{j+1}.\]
$S$ is a multiplicatively closed set so $s'_{j+1}= \bar{s}_{j+1}s'_j$
is an element of $S$.
Now we have
\[\begin{array}{c}
\begin{array}{ll}
 (1 \otimes s'_{j+1})(\sum_{i = 1}^{j+1} a_i \otimes b_i ) = &
 (1 \otimes \bar{s}_{j+1})(1 \otimes s'_{j})
(\sum_{i=1}^{j} a_i \otimes b_i) \, + \\
& \,\,\,+\,
(1 \otimes \bar{s}_{j+1})(1 \otimes s'_{j})(a_{j+1} \otimes b_{j+1})
\end{array} \\
\begin{array}{l}
\,\,\,\,\,\,\,\,\,\,\,\,\,\,\,\,
= (\sum_{i = 1}^j a_i \otimes  \bar{s}_{j+1} r'^{j}_i)(1 \otimes s)
+ (a_{j+1} \otimes r'^{j+1}_{j+1})(1 \otimes s)\\
\,\,\,\,\,\,\,\,\,\,\,\,\,\,\,\,
= (\sum_{i = 1}^{j+1} a_i \otimes r'^{j+1}_i)(1 \otimes s)
\end{array}\end{array}\]
where we denoted $r'^{j+1}_i = \bar{s}_{j+1}r'^{j}_i$ for
$i < j+1$ and $r'^{j+1}_{j+1}$ has already been defined.

\minsection{Ore localization for modules}

The modern point of view on Ore localization is to express it as a
{\it localization functor} on some category of modules. The
localization map $\iota : R \rightarrow S^{-1}R$ will be replaced
by a localization functor $Q_S^*$ from $R$-Mod to $S^{-1}R$-Mod.

\ppt Let $S$ be a left Ore set in a ring $R$, and $M$ a left
$R$-module. Notice that $S^{-1}R$ may be viewed as an
$S^{-1}R$-$R$-bimodule. The {\bf module of fractions} $S^{-1}M$ of
$M$ with respect to Ore set $S$ is the left $S^{-1}R$-module
\begin{equation}\label{eq:Orelocmod}
Q^*_S(M) = S^{-1} M := S^{-1} R \otimes_R M.
\end{equation}
For each morphism $f : M \rightarrow N$, set
$Q^*_S(f):= 1 \otimes f : S^{-1}R\otimes_R M \rightarrow S^{-1}R\otimes_R N$.
This defines a localization functor
$Q^*_S : R-{\rm mod} \rightarrow S^{-1}R-{\rm mod}$
whose right adjoint is the obvious forgetful functor
$Q_{S*}: S^{-1}R-{\rm mod}\rightarrow R-{\rm mod}$.
If $\iota = \iota_R : R \rightarrow S^{-1}R$ is the localization map,
then define the map of $R$-modules
$\iota_M : M \rightarrow S^{-1}M$ via $\iota_M = \iota_R \otimes_R \id$
i.e. $m \mapsto 1\otimes m$, also called the localization map.
Maps $\iota_M$ make together a natural transformation of functors, namely
the adjunction $\iota : {\rm Id} \rightarrow Q_{S*} Q_S^*$.

\ppta {\bf Remark.} If $S$ is a right Ore set, and $M$ a left $R$-module,
then $M[S^{-1}]:= R[S^{-1}]\otimes_R M$. If $N$ is a right $R$-module
then view $RS^{-1}$ (or $S^{-1}R$) as a $R$-$RS^{-1}$-
(resp. $R$-$S^{-1}R$)- bimodule and define $Q_S^*(N) := N\otimes_R R[S^{-1}]$
(resp. $N \otimes_R S^{-1}R$).
We emphasize that the choice of right vs. left Ore sets is
not correlated with the choice of right or left $R$-module categories,
at least in the principle of the construction.

\ppt {\bf Universal property.} For given $R,M,S$ as above we define
the category ${\cal M} = {\cal M}(R,M,S)$. The objects of ${\cal M}$
are pairs $(N,h)$ where $N$ is a left $S^{-1}R$-module
and $h : M \rightarrow {}_{R}N$ a map of left $R$-modules.
A morphism of pairs $\mu : (N,h) \rightarrow (N',h')$ is
a map $\mu : N \rightarrow N'$ of $S^{-1}R$-modules such that
$h' = \mu \circ h$.

{\bf Proposition.} {\it The pair $(S^{-1}M,\iota_M)$ is 
an initial object in ${\cal M}$.}

{\it Proof.}
For any pair $(N,h) \in {\rm Ob}({\cal M})$ there is a well-defined
morphism
\begin{equation} 
\alpha :(S^{-1} M,\iota_M)\rightarrow  (N,h)
\,\,\mbox{  by  }\,\,
\alpha(s^{-1} \otimes_R m) = s^{-1}h(m).\end{equation}
Let now $\alpha'$ be any morphism from $(S^{-1} M,\iota)$ to $(N,h)$.
By $h = \alpha' \circ i$ we conclude
\[ \alpha'(1 \otimes m) = h(m) = \alpha(1 \otimes m),
\,\, \forall m \in M.\]
The elements of the form $1 \otimes m$ generate $S^{-1} M$
as a module over $S^{-1}R$ and therefore $\alpha = \alpha'$.

\ppt Elements in the tensor product~(\ref{eq:Orelocmod}) are of the form
$\sum_i s^{-1}_i \otimes m_i$ but such can be added up to
a single term of that form, as the fractions can be always put to
the common denominator. Namely, by the left Ore condition
$\forall s,s' \in S \,\exists \tilde{s} \in S \,
\exists \tilde{r} \in R$, $s' s^{-1} = \tilde{s}^{-1} \tilde{r}$,
and therefore
\begin{equation}\label{eq:comdenMod} \begin{array}{lcl}
s^{-1} \otimes m + s'^{-1} \otimes m' &=& s'^{-1} s' s^{-1}\otimes m
+ s'^{-1} \tilde{s}^{-1} \tilde{s}\otimes m' \\ &=&
s'^{-1} \tilde{s}^{-1} \otimes (\tilde{r}m + \tilde{s}m').
\end{array}\end{equation}
Thus every element in $S^{-1}M$ may be written in the form $s^{-1}
\otimes_R m$, hence there is a surjection of sets $\nu : S \times
M \to S^{-1}M$. The set $S\times M$ may be viewed as a set retract
$S \times \{1\} \times M$ of $S \times R\times M$ via the
retraction $(s,r,m) \mapsto (s,rm)$. Clearly $\nu$ extends to
$\nu' : S \times R \times M\to S^{-1}M$. By the universality of
the free \abelian group ${\DDl Z} (S \times R \times M)$ with
basis $S \times R \times M$, $\exists!\tilde{\nu'} : {\DDl Z} (S
\times R \times M)\rightarrow S^{-1}M$ which is additive and
extends $\nu'$. It is clear by checking on the basis elements
$(s,r,m)$ and uniqueness that the composition of the canonical
projections $ {\DDl Z} (S \times R \times M)\rightarrow {\DDl Z}
(S^{-1}R \times M) \rightarrow S^{-1}R \otimes_R M$ equals
$\tilde{\nu'}$.

For $r \in R$ for which $rs \in S$,
$s^{-1} \otimes_R m = (rs)^{-1} r \otimes_R m = (rs)^{-1} \otimes_R rm$
implying that ${\rm ker}\,\nu'\subset {\DDl Z}(S\times M)$
contains all differences $(s,m)- (s',m')$ in ${\DDl Z}(S\times M)$
of pairs in $S\times M$ which are equivalent via
\begin{equation}\label{eq:kernelNu}
\fbox{$
(s,m) \sim (s',m') \,\Leftrightarrow  \,\exists r,r' \in R\,\,\,\,\,
 rs = r's'\in S \,\, \mbox{ and } \,\, rm = r'm'.
$}
\end{equation}
{\bf Lemma.} {\it (i) On $(S\times M)/\sim$
there is a unique binary operation $+$ such that
\begin{equation}\label{eq:smsmsmm} (s,m) + (s,m') \sim (s, m+m').
\end{equation}
(ii) $((S\times M)/\sim,+)$ is an \abelian group.
Hence by the universality of the free \abelian group,
the composition
$S \times R \times M \to S\times M\rightarrow (S\times M)/\sim$
extends to a unique map
$p : {\DDl Z}(S\times R \times M) \to (S\times M)/\sim$
of \abelian groups.

(iii) The map $p$ factors to a map $p' : S^{-1}R \otimes_R M
\rightarrow (S\times M)/\sim$.

(iv) $p'$ is an inverse of $\nu$, hence $p'$ respects addition.
}

{\it Proof.} (i) {\it Uniqueness.} Suppose there are two such
operations, $+_1,+_2$ and two classes $(s_1,m_1)$ and $(s_2,m_2)$
on which $+_1$ and $+_2$ disagree. By the left Ore condition
choose $\tilde{s} \in S, \tilde{r}\in R$ with $\tilde{s}s_1 =
\tilde{r}s_2$. Then $(s_1,m_1) +_i (s_2,m_2) \sim
(\tilde{s}s_1,\tilde{s}m_1) +_i (\tilde{r}s_2,\tilde{r}m_2) \sim
(\tilde{s}s_1, \tilde{s}m_1 + \tilde{r}m_2)$ which agree for $i =
1,2$, giving a contradiction.

{\it Existence.} Take  $(s_1,m_1) + (s_2,m_2) := (\tilde{r}_1 s_1,
\tilde{r}_1 m_1 + \tilde{r}_2 m_2)$ with any choice of
$\tilde{r}_1,\tilde{r}_2$ such that
$\tilde{r}_1 s_1 = \tilde{r}_2 s_2 \in S$.
We verify that the class of the result does not depend on the choices.
If $s_1,m_1$ are replaced by $rs_1 \in S,rm_1\in M$ we can by the
combined Ore condition choose
$r_*,s_*$ with $s_* r = r_* \tilde{r}_1$, hence
$S \ni s_* r s_1 = r_* \tilde{r}_1 s_1 = r_* \tilde{r}_2 s_2$.
Hence the rule for the sum gives
$(s_* r s_1, s_* r m_1 + r_* \tilde{r}_2 m_2) \sim
(r_* \tilde{r}_1 s_1, r_* \tilde{r}_1 m_1 + r_* \tilde{r}_2 m_2)
\sim (\tilde{r}_1 s_1, \tilde{r}_1 m_1 + \tilde{r}_2 m_2)$.
By symmetry, we have the same independence for choice of $(s_2,m_2)$.
Finally, suppose that instead of choosing $\tilde{r}_1, \tilde{r}_2$
we choose $\tilde{b}_1,\tilde{b}_2$. As $\tilde{r_1} s_1 \in S$.
By the combined Ore condition we may choose
$s_\sharp \in S$, $r_\sharp \in R$, such that
$r_\sharp \tilde{r}_1 = s_\sharp \tilde{b}_1$
with $r_\sharp \tilde{r}_1 s_1 =  s_\sharp (\tilde{b}_1 s_1) \in S$.
Hence
$(\tilde{r}_1 s_1, \tilde{r}_1 m_1 + \tilde{r}_2 m_2)
 \sim
(r_\sharp \tilde{r}_1 s_1, r_\sharp \tilde{r}_1 m_1 +
r_\sharp \tilde{r}_2 m_2)\sim
(s_\sharp \tilde{b}_1 s_1, s_\sharp \tilde{b}_1 m_1 +
r_\sharp \tilde{r}_2 m_2).
$
Now $r_\sharp \tilde{r}_2 s_2 =  r_\sharp \tilde{r}_1 s_1 =
s_\sharp \tilde{b}_1 s_1 = s_\sharp \tilde{b}_2 s_2$, hence
by left reversibility, there is $s_2^*$ such that
$s_2^* r_\sharp \tilde{r}_2 = s_2^*  s_\sharp \tilde{b}_2$.
Thus
$
(s_\sharp \tilde{b}_1 s_1, s_\sharp \tilde{b}_1 m_1 +
 r_\sharp \tilde{r}_2 m_2) \sim
(s_2^* s_\sharp \tilde{b}_1 s_1,
s_2^* s_\sharp \tilde{b}_1 m_1 + s_2^* s_\sharp \tilde{b}_2 m_2)
\sim (\tilde{b}_1 s_1,
\tilde{b}_1 m_1 + \tilde{b}_2 m_2),$ as required.

(ii) In the proof of existence in (i) we have seen that $+$ is
commutative. Notice also that the class of $(s,0)$ (independent on
$S$) is the neutral element. For any pair, and hence for any
triple of elements in $S\times M/\sim$, we can choose
representatives such that all three are of the form $(s,m)$ with
the same $s$. For such triples the associativity directly follows
by applying~(\ref{eq:smsmsmm}).

(iii) As $p$ and the projection
$S \times R \times M \to S\times M\rightarrow (S\times M)/\sim$
are additive it is sufficient to show
that $p$ sends the kernel of the projection
to $0 \in (S\times M)/\sim$. The kernel of the projection
is spanned by the elements of several obvious types, so
we check for generators.

1. $(s,r,m) - (s',r',m)$ where $s^{-1}r = (s')^{-1}r'$.
That means that for some $\tilde{s}\in S$, $\tilde{r} \in R$ we have
$\tilde{s}s= \tilde{r}s'$ and $\tilde{s}r= \tilde{r}r'$.
Compute
$p(s,r,m) - p(s',r',m') = (s,rm) + (s', -r'm') =
(\tilde{s}s,\tilde{s}rm)+(\tilde{r}s', -\tilde{r}r'm')
= (\tilde{s}s, \tilde{s}rm - \tilde{r}r'm') = 0$.

2. Elements $(s,r+r', m) - (s,r,m) - (s,r+r',m)$,
as well as $(s,rr',m)-(s,r,r'm)$ and
$(s,r,m+m') - (s,r,m) - (s,r,m')$
go to $0$ because, by (ii),  in computing $p$ one has to first
act with the second component to the third.

\ppt {\bf Proposition.} \ldef{relonStimesM}{\it $(S\times M)/\sim$
is additively canonically isomorphic to $S^{-1}M$. This
isomorphism equips $(S\times M)/\sim$ with the canonical left
$S^{-1}R$-module structure for which the following formulas can be
taken as defining:
\[ t^{-1}r (s,m) = (s_* t, r_* m) \in (S\times M)/\sim\]
where $s_* \in S$ and $r_* \in R$ such that $s_* r = r_* s$.
} 

{\it Proof.} By the lemma we have the first statement and hence
we can view the class of $(s,m)$ as $s^{-1}\otimes m$ and
the formulas follow.
They are defining because
the formulas agree with the action transferred by an
isomorphism, hence the existence, and by the left Ore condition
for any  $t^{-1}r$ and $(s,m)$ there are
$s_*$, $r_*$ qualifying for the formulas, hence the uniqueness.

\ppta {\bf Corollary.}\ldef{torsionOre}
{\it  Let $S\subset R$ be left Ore
and $M$ a left $R$-module.
Then $m = 0$ for some $s$ iff  $\exists s \in S$ and $s m = 0$.
}

\pptb (long) {\bf Exercise.} Let $S,W$ be two left Ore sets in $R$,
and $M = {}_R M$ a left $R$-module.
The relation $\sim$ on $S \times W \times M$ given by
\[
(s,w,m) \sim (s',w',m') \Leftrightarrow
\left\lbrace
\begin{array}{l}
\exists r,r',u,u',\tilde{r},\tilde{r}'\in R,
\exists \tilde{w},\tilde{w}'\in W,\\
rs = r's' \in S, u\tilde{w} = u'\tilde{w}' \in W,\\
\tilde{w}r = \tilde{r}w,\tilde{w}'r' = \tilde{r}'w',
z\tilde{r}m = z'\tilde{r}'m',
\end{array}
\right.\] is an equivalence relation. Map $(S \times W \times
M/\sim) \to S^{-1}W^{-1}M$ given by $(s,w,m)\mapsto s^{-1}\otimes
w^{-1} \otimes m$ is well-defined and bijective.
\minsection{Monads, comonads and gluing}

\ppt A {\bf monoidal category} is a category ${\cal C}$ equipped
with a 'tensor' product bifunctor $\otimes : {\cal C} \times {\cal
C} \rightarrow {\cal C}$; a distinguished object $1_{\cal C}$, a
family of associativity isomorphisms $c_{ABC} : (A \otimes B)
\otimes C \rightarrow A \otimes (B \otimes C)$, natural in objects
$A,B,C$ in ${\cal C}$; left unit $l_A : 1_{\cal C} \otimes A
\rightarrow A$, and right unit isomorphisms $r_A : A \otimes
1_{\cal C} \rightarrow A$, both indexed by and natural in objects
$A$ in ${\cal C}$; furthermore, require some standard coherence
conditions (pentagon axiom\index{pentagon axiom} for associativity
coherence; left and right unit coherence \index{coherence}
conditions, cf.~\cite{Borceux, MacLane}). A monoidal category
$(\tilde{\cal C},\otimes, 1_{\cal C}, c, r, l)$ is {\bf strict} if
$c_{ABC},l_A,r_A$ are actually all identity morphisms.

\ppt {\bf Monads and comonads}. \ldef{monadscomonads}
Given a diagram of categories ${\cal A,B,C}$,
functors $f_1,f_2,g_1,g_2$ and natural transformations $F,G$ as follows
\begin{equation}\label{eq:horcomp}\diagram{
      {\cal A}\,\,\, \fhUp{\fhd{f_1}{}}{\Uparrow \!\!\!F}{\fhd{}{f_2}}
\,\,\,{\cal B}\,\,\, \fhUp{\fhd{g_1}{}}{\Uparrow \!\!\!G}{\fhd{}{g_2}}
\,\,\,{\cal C},
}\end{equation}
one defines the natural transformation $G\star F : g_2 \circ f_2
\Rightarrow g_1 \circ f_1$ by
\[
 ( G\star F )_A := G_{f_1(A)} \circ g_2(F_A) = g_1(F_A) \circ G_{f_2(A)}
:  g_2(f_2(A)) \rightarrow g_1 (f_1(A)).\]
$(F,G) \mapsto F\star G$ is called the \wind{Godement product}
('horizontal composition'\index{horizontal composition},
cf.~(\ref{eq:horcomp})).
It is associative for triples for which $F\star(G\star H)$ is defined.

Given functors $f,g,h: {\cal A}\rightarrow {\cal B}$ and
natural transformations $\alpha : f \Rightarrow g$,
$\beta : g \Rightarrow h$, define their 'vertical' (or ordinary)
composition $\beta\circ \alpha : f \Rightarrow h$ to be their
composition taken objectwise:
$(\beta\circ \alpha)_A := \beta_A \circ \alpha_A : f(A) \rightarrow h(A)$.

Note the \wind{interchange law}:
$(\alpha \circ \beta) \star (\gamma \circ \delta) =
(\alpha \star \gamma) \circ (\beta\star\delta)$.

{\footnotesize If $T, T'$ are endofunctors in ${\cal A}$
 and $\alpha : T \Rightarrow T'$, $\beta : T' \Rightarrow T$
natural transformations, one may also use the concatenation
notation: $\alpha T : TT \Rightarrow T'T$ is given by $(\alpha
T)_M := \alpha_{TM} : T(TM) \rightarrow T'(TM)$, i.e. $\alpha T
\equiv \alpha \star 1_T$; similarly, $T\alpha$ equals $1_T \star
\alpha$, or, $(T\alpha)_M = T(\alpha_M) : TTM \rightarrow TT'M$.
This notation extends to the sequences with more functors but only
one natural transformation -- it is preferable to specify the product
$\circ$ versus $\star$ among the words if each has some natural
transformation mixed in. Here concatenation is higher binding than
any of the composition signs. Notice that $1_T \star 1_S = 1_{T \circ
S}$. }

Given a strict monoidal category
$\tilde{\cal C} := ({\cal C},\otimes, 1_{\tilde{\cal C}})$
a {\bf monoid} in ${\cal C}$ is a pair $(X,\mu)$ of an object $X$ and
a multiplication morphism $\mu : X \otimes X \to X$
which is associative and there is a 'unit' map
$\eta : 1_{\tilde{\cal C}} \rightarrow X$ such that
$\mu \circ (\eta \otimes \id) = \mu \circ (\id \otimes \eta) \cong
\id$ (here the identification $1_{\tilde{\cal C}} \otimes X \cong X$
is assumed). As this characterizes the unit map uniquely,
one may or may not include $\eta$ in the data,
writing triples $(X,\mu,\eta)$ when convenient.

For a fixed small category ${\cal A}$, the category ${\rm
End}_{\cal A}$ of endofunctors (as objects) and natural
transformations (as morphisms) is a strict monoidal category: the
product of endofunctors is the composition, the product of natural
transformations is the Godement product, and the unit is ${\rm
Id}_{\cal A}$.

A \wind{monad} $(T,\mu,\eta)$ in an arbitrary category ${\cal A}$
is a monoid in ${\rm End}_{\cal A}$,
and a \wind{comonad} $(\perp,\delta,\eta)$
in ${\cal A}$ is a monad in ${\cal A}^{\rm op}$. The natural
transformations $\delta : \perp\to\perp\circ\perp$
and $\epsilon : \perp \to {\rm Id}_{\cal A}$
are also called the coproduct \index{coproduct!of a comonad}
and the counit of the comonad\index{counit of comonad}
respectively.

An {\bf action} of a monad ${\bf T} = (T,\mu,\eta)$
\index{action of monad}
on an object $M$ in ${\cal A}$
is a morphism $\nu : T(M) \rightarrow M$ such that the diagram
\[
\begin{array}{ll}
\begin{array}{ccc}
TT(M) & \stackrel{\mu_M}\to & T(M) \\
T(\nu) \downarrow && \downarrow \nu \\
T(M) & \stackrel{\nu}\to & M
\end{array}
\end{array}
\]
commutes and $\nu \circ \eta_M = {\rm Id}_{M}$.
We say that $(M,\mu)$ is a {\bf module}\index{module over monad}
(older terminology: algebra\index{algebra over monad})
over ${\bf T}$. A map $(M,\mu) \rightarrow (N,\nu)$ is a morphism
$f : M \rightarrow N$ in ${\cal A}$ intertwining the actions in the
sense that $f \circ \nu_M = \nu_N \circ T(f) : T(M) \rightarrow N$.
For a fixed ${\bf T}$, modules and their maps make a category
${\cal A}^{\bf T}\equiv {\bf T}-{\rm Mod}$,
called the \wind{Eilenberg-Moore} category
of ${\bf T}$ (\cite{EilMoore:trip}). The natural forgetful functor
$U^{\bf T} : {\cal A}^{\bf T}\rightarrow {\cal A}$, $(M,\nu) \mapsto M$
is faithful,
reflects isomorphisms and has a left adjoint $F : M \mapsto (M,\mu_M)$.
The unit of adjunction
$\eta : {\rm Id}_{\cal A} \Rightarrow U^{\bf T} F = T$
coincides with the unit of ${\bf T}$,
and the counit
$\epsilon : FU^{\bf T} \rightarrow {\rm Id}_{{\cal A}^{\bf T}}$
is given by $\epsilon_{(M,\nu)} = \nu$.
The essential image of $F$, is a full and replete
subcategory ${\cal A}_{\bf T}\subset {\cal A}^{\bf T}$
and it is called the {\bf Kleisli category} of ${\bf T}$,
while its objects are called {\bf free ${\bf T}$-modules}.

Dually, for a comonad $\comdG = (\comdg,\delta,\epsilon)$ in ${\cal A}$,
a $\comdG$-{\bf comodule} is an object in category
$\comdG-{\rm Comod} :=
(({\cal A}^{\rm op})^\comdG)^{\rm op} = (\comdG-{\rm Mod})^{\rm op}$;
 equivalently it is a
pair $(M,\sigma)$ where $M$ is an object in ${\cal A}$ and
$\sigma : M \rightarrow \comdg(M)$ is a morphism in ${\cal A}$
such that
\[
\begin{array}{ll}
\begin{array}{ccc}
M & \stackrel{\sigma}\rightarrow & \comdg(M) \\
\sigma \downarrow && \downarrow \comdg(\sigma) \\
\comdg(M) & \stackrel{\delta}\rightarrow & \comdg\comdg(M)
\end{array}
\end{array}
\]
and $\epsilon_M \circ \sigma = {\rm Id}_{M}$.
A map $(M,\rho_M) \rightarrow (N,\nu_N)$ is a morphism
$f : M \rightarrow N$ in ${\cal A}$ intertwining the coactions in the
sense that $\sigma_N \circ f = T(f) \circ \sigma_M : M \rightarrow T(N)$.
The forgetful functor $\comdG-{\rm Mod} \rightarrow {\cal A}^{\rm op}$
may be interpreted as a functor
$U_\comdG : \comdG-{\rm Comod}\rightarrow {\cal A}$ which thus
has a {\it right} adjoint $H : M \mapsto (\comdg(M),\delta)$ whose
essential image by definition consists of
{\bf cofree $\comdG$-comodules}\index{cofree $\comdG$-comodule}.
The counit of the adjunction agrees with the counit of the
comonad $\epsilon : \comdg = U_\comdG H \Rightarrow {\rm Id}_{\cal A}$
and the unit $\eta : {\rm Id}_{\comdG-{\rm Comod}}\Rightarrow HU_{\comdG}$
is given by
$\eta_{(M,\sigma)}= \sigma : (M,\sigma) \rightarrow (\comdg(M),\delta)$.

\ppt \ldef{monadfromadjoint}
An archetypal example of a monad is constructed from a pair
\index{monad!associated to adjunction}
of adjoint functors $Q^* \dashv Q_*$ \index{$\dashv$} where
$Q_* : {\cal B}\rightarrow {\cal A}$.\index{adjoint functors}
In other words there are natural transformations
$\eta : {\rm Id}_{\cal A} \Rightarrow Q_* Q^*$ called the unit
\index{unit of adjunction} and
$\epsilon : Q^* Q_* \Rightarrow {\rm Id}_{\cal B}$, called
the counit of the adjunction, \index{counit of adjunction}\index{adjunction}
such that the composites in the two diagrams
\[ Q_* \stackrel{\eta_{Q_*}}{\rightarrow}
Q_* Q^* Q_* \stackrel{Q_*(\epsilon)}{\rightarrow} Q_*,
\,\,\,\,\,\,
Q^* \stackrel{Q^*(\eta)}{\rightarrow} Q^*Q_* Q^*
\stackrel{\epsilon_{Q_*}}{\rightarrow} Q^*,
\]
are the identity transformations.
Then
${\bf T} := (Q_* Q^*, 1_{Q_*} \star \epsilon \star 1_{Q^*}, \eta)$
is a monad in ${\cal A}$. In other words, the multiplication is given  by
\[ \mu_M  =  Q_*(\epsilon_{Q^*(M)})
: Q_* Q^* Q_* Q^* (M) = TT(M) \rightarrow Q_* Q^* (M) = T(M).\]
The {\bf comparison functor} \index{comparison functor}
$K^{\bf T}: {\cal B} \rightarrow {\cal A}^{\bf T}$ is defined by
\[ M \mapsto (Q_*(M), Q_*(\epsilon_M)),\,\,\,\,\,\,\,F \mapsto Q_*(f).\]
It is full and $Q^*$ factorizes as
${\cal B} \stackrel{K^{\bf T}}{\rightarrow} {\cal A}^{\bf T}
\stackrel{U^{\bf T}}{\rightarrow} {\cal A}$.
More than one adjunction (varying ${\cal B}$)
may generate the same monad in ${\cal A}$ in this vein.

Dually, ${\bf G} := (Q^*Q_*, Q^*\eta Q_*, \epsilon)$
is a comonad, i.e. a monad in
${\cal B}^{\rm op}$.
The comparison functor
$K^{\bf G} : {\cal A}^{\rm op}\rightarrow ({\cal B}^{\rm op})^{\bf T}$
is usually identified with a 'comparison functor'
$K_{\bf G} : {\cal A} \rightarrow  (({\cal B}^{\rm op})^{\bf T})^{\rm op}
\equiv \comdG-{\rm Comod}$ which is hence given by
$N \mapsto (Q^*(N),\eta_{Q^*(N)})$. $K_{\bf G}$ is full
and $Q_*$ factorizes as
${\cal A} \stackrel{K_\comdG}{\rightarrow} \comdG-{\rm Comod}
\stackrel{U_\comdG}{\rightarrow} {\cal B}$.

\ppt A map of monoids $f : (A,\mu,\eta)\rightarrow (A',\mu',\eta')$
in a monoidal category $({\cal A},\otimes, 1_{\cal A}, a, l, r)$
is a morphism $f : A \rightarrow A'$ in ${\cal A}$,
commuting with multiplication: $\mu\circ(f\otimes f) = f \circ \mu$;
and with the unit map: $\eta' \circ \mu = \eta \otimes \eta$, where
on the left the application of one of the isomorphisms
$l_{1_{\cal A}}, r_{1_{\cal A}} :
1_{\cal A} \otimes 1_{\cal A} \rightarrow 1_{\cal A}$
is assumed.
In particular, the morphism \index{morphism of monads}
$\phi : (T, \mu, \eta) \rightarrow (T', \mu', \eta')$
of monads in ${\cal A}$ is a
natural transformation $\phi : T \rightarrow T'$ such that
$\mu \circ (\phi \star \phi) = \phi \circ \mu : TT \Rightarrow T'$
and $\eta' \circ \mu = \eta \star \eta : TT \Rightarrow {\rm Id}_{\cal A}$.
If $M$ is an object in ${\cal A}$ and $\nu$ a $T'$-action on $M$,
then $\nu' \circ \phi_M : TM \Rightarrow M$ is a $T$-action on $M$.
More precisely, a natural transformation $\phi : T \Rightarrow T'$
and rules ${\cal A}^\phi (M, \nu) = (M, \nu' \circ \phi_M)$
and ${\cal A}^\phi(f) = f$ define a functor
\[
 {\cal A}^\phi : {\cal A}^{{\bf T'}} \rightarrow {\cal A}^{{\bf T}}
\]
iff $\phi$ is a morphism of monads and every functor
${\cal A}^{\bf T'} \rightarrow {\cal A}^{\bf T}$ inducing
the identity on ${\cal A}$ is of that form.

\ppt Let ${\bf \Delta}$ be the 'simplicial' category:
its objects are nonnegative integers viewed as
finite ordered sets ${\bf n} := \{ 0 < 1 <  \ldots < n\}$ and
its morphisms are nondecreasing monotone functions.
Given a category ${\cal A}$, denote by ${\rm Sim}{\cal A}$
the category of {\bf simplicial objects} in ${\cal A}$, i.e. functors
$F : {\bf \Delta}^{\rm op}\rightarrow {\cal A}$.
Represent $F$ in ${\rm Sim}{\cal A}$ as a sequence
$F_n := F({\bf n})$ of objects, together
with the face maps $\partial_i^n : F_n \rightarrow F_{n-1}$
and the degeneracy maps $\sigma_i^n : F_n \rightarrow F_{n+1}$
for $i \in {\bf n}$ satisfying
the familiar simplicial identities~(\cite{Weibel, MacLane}).
Notation $F_\bullet$ for this data is standard.

Given a comonad $\comdG$ in ${\cal A}$ one defines
the sequence $\comdG_\bullet$ of endofunctors
${\DDl Z}_{\geq 0} \ni n \mapsto \comdG_n := \comdg^{n+1} :=
\comdg \circ \comdg \circ \ldots \circ \comdg$,
together with natural transformations
$\partial_i^n : \comdg^i \epsilon \comdg^{n-i} :
\comdg^{n+1} \rightarrow \comdg^n$
and $\sigma_i^n : \comdg^i \delta \comdg^{n-i} :
\comdg^{n+1} \rightarrow \comdg^{n+2}$,
satisfying the simplicial identities.
Use
$\epsilon \comdg \circ \delta = \comdg \epsilon \circ \delta
= {\rm Id}_{\cal A}$
in the proof.
Hence any comonad $\comdG$ canonically induces a
simplicial endofunctor\index{simplicial endofunctor},
i.e. a functor $\comdG_\bullet : \Delta^{\rm op}
\rightarrow {\rm End}{\cal A}$, or equivalently, a functor
$\comdG_\bullet : {\cal A} \rightarrow {\rm Sim}{\cal A}$.
The counit $\epsilon$ of the comonad $\comdG$ satisfies
$\epsilon \circ \partial_0^1 = \epsilon \circ \partial_1^1$, hence
$\epsilon : \comdG_\bullet \rightarrow {\rm Id}_{\cal A}$ is in fact
an augmented simplicial endofunctor.

This fact is widely used
in homological algebra~(\cite{Barr:composite, MacLane, Weibel}),
and now also in the cohomological study of
noncommutative spaces~(\cite{Ros:NcSch}).

\ppt {\bf Barr-Beck lemma.}~(\cite{BarrWells,MLMoerd:Sheaves})
\index{Barr-Beck lemma}
Let $Q^* \dashv Q_*$ be an adjoint pair
${\bf T}$ its associated monad,
and ${\bf G}$ its associated comonad
(as in~\luse{monadfromadjoint}).
Recall the notions of preserving
and reflecting (co)limits from~\luse{presrefl}.

{\it If $Q_*$ preserves and reflects coequalizers of all parallel pairs
in ${\cal A}$ (for which coequalizers exists) and if any parallel
pair mapped by $Q_*$ into a pair having a coequalizer in ${\cal B}$
has a coequalizer in ${\cal A}$ then the comparison
functor $K : {\cal B} \rightarrow {\cal A}^{\bf T}$
is an {\it equivalence} of categories.

If $Q^*$ preserves and reflects equalizers of all parallel pairs
in ${\cal B}$ (for which equalizers exists) and if any parallel
pair mapped by $Q^*$ into a pair having an equalizer in ${\cal A}$
has an equalizer in ${\cal B}$ then the comparison functor
$K' : {\cal A} \rightarrow {\bf G}-{\rm Comod}$ is an equivalence
of categories.
}

Left (right) exact functors by definition preserve finite limits
(colimits) and faithful functors clearly reflect both. In particular,
this holds for (co)equalizers of parallel pairs.
In \abelian categories such (co)equalizers always exist. Hence

{\bf Corollary.} {\it Consider an adjoint pair
$Q^*\dashv Q_*$ of additive functors between \abelian categories.
If $Q_*$ is faithful and exact, then the comparison functor for
the associated comonad is an equivalence. If $Q^*$ is
faithful and exact, then the comparison functor for the associated
monad is an equivalence.
}

Given a functor  $U : \tilde{\cal M}\to {\cal M}$
one may ask when there is a monad $\bf T$ in $\cal M$ and an
equivalence $H : \tilde{\cal M}\to{\cal M}^{\bf T}$ such that $U^T H = U$.
The conditions are given by the {\it Beck monadicity
(=tripleability) theorem(s)~(\cite{BarrWells,Beck,MacLane,MLMoerd:Sheaves})}.
If we already know that $U$ has left adjoint,
this may be rephrased by asking if
the comparison functor for the {\it associated} monad
is an equivalence. Barr-Beck lemma gives only sufficient conditions
for this case, it is easier to use, and widely applicable.

\ppt (\cite{Ros:NcSch}) \ldef{comonadfamily}
{\bf A comonad associated to a family
of continuous functors.}
Let $\{ Q^*_\lambda : {\cal A} \rightarrow
{\cal B}_\lambda\}_{\lambda \in \Lambda}$
be a small family of continuous
(= having a right adjoint) functors.
The categories ${\cal B}_\lambda$ are {\it not}
necessarily constructed from ${\cal A}$ by a localization.

One may consider the category
${\cal B}_\Lambda : = \prod_{\lambda \in \Lambda}{\cal B}_\lambda$
whose objects are families
$\prod_{\lambda \in \Lambda} M^\lambda$ of objects $M^\lambda$ in
${\cal B}_\lambda$ and morphisms are families
$\prod_{\lambda \in \Lambda}f_\lambda :
\prod_{\lambda \in \Lambda} M^\lambda \rightarrow
\prod_{\lambda \in \Lambda} N^\lambda$ where
$f_\lambda : M^\lambda \rightarrow N^\lambda$
is a morphism in ${\cal B}_\lambda$, with componentwise
composition. This makes sense
as the family of objects is literally a function from
$\Lambda$ to the disjoint union $\coprod_\lambda {\rm Ob}\,{\cal B}_\lambda$
which is in the same Grothendieck universe.

The family of adjoint pairs
$Q^*_\lambda \dashv Q_{\lambda *}$ defines an inverse image functor
${\bf Q}^* = \prod Q^*_\lambda : {\cal A} \rightarrow {\cal A}_\Lambda$
by ${\bf Q}^*(M) :=  \prod_{\lambda \in \Lambda}Q^*_\lambda(M)$ on objects
and ${\bf Q}^*(f) := \prod_\lambda Q^*_\lambda(f)$ on morphisms.
However, a direct image functor may not exist.
We may naturally try
${\bf Q}_* : \prod'_\lambda M^\lambda \mapsto \prod'_\lambda Q_*(M^\lambda)$
where $\prod'$ is now the symbol for the {\it Cartesian} product
in ${\cal A}$ which may not always exist.
For {\it finite} families, with ${\cal A}$ \abelianM,
these trivially exist.
Let ${\cal A}^\Lambda = \prod_{\lambda\in\Lambda} {\cal A}$
be the power category.
Assume a fixed choice of the Cartesian product for all $\Lambda$-tuples
in ${\cal A}$.
Then $\{M^\lambda\}_\lambda \mapsto
\prod'_{\lambda\in \Lambda} M^\lambda$
extends to a functor ${\cal A}^\Lambda\to {\cal A}$,
and the universality of products implies that the projections
$p'_{\nu \{M^\lambda\}_\lambda} : \prod_\lambda M^\lambda \rightarrow M^\nu$
form a natural transformation of functors
$p'_\nu : \prod'_{\lambda} {\rm Id}_{\cal A}\Rightarrow {\rm p}_{\nu}$
where $p_\nu : {\cal A}^\Lambda \rightarrow {\cal A}$
is the $\nu$-th formal projection $\prod_\lambda M^\lambda\to M^\nu$.
The unique liftings
$\bm{\eta}_M : M \rightarrow {\bf Q}_* {\bf Q}^* (M)$
of morphisms $\eta_{\nu M} : M \rightarrow Q_{\nu *} Q^*_\nu (M)$
in the sense that
$(\forall \nu)\, \eta_{\nu M} = p'_{\nu M} \circ \eta_M$ hence
form a natural transformation
$\bm{\eta} : {\rm Id}_{\cal A} \Rightarrow {\bf Q}_* {\bf Q}^*$.

Define $\bm{\epsilon} \equiv \prod_\lambda \epsilon_\lambda
  : \prod_\lambda Q^*_\lambda Q_{\lambda *} \Rightarrow
\prod_\lambda {\rm Id}_{{\cal B}_\lambda} = {\rm Id}_{{\cal B}_\Lambda}$
componentwise.
This way we obtain an adjunction ${\bf Q}^* \dashv {\bf Q}_*$.
If $Q_{\lambda *}$ is faithful and exact for every
$\lambda$ then ${\bf Q}_*$ is as well.

Consider the comonad $\comdG$
in ${\cal B}_\Lambda$ associated to ${\bf Q}^* \dashv {\bf Q}_*$.
We are interested in situation when
the comparison functor $K_\comdG$ is an equivalence of categories.
That type of a situation arises in practice in two different ways:
\begin{itemize}
\item 1) All categories ${\cal A}, {\cal B}_\lambda$
and flat localization functors $Q^*_\lambda, Q_{\lambda *}$ are given
at start and the construction is such that
we know the faithfulness of ${\bf Q}_*$.
\item 2) Only categories ${\cal B}_\lambda$ are given (not ${\cal A}$)
but equipped with {\bf gluing morphisms} i.e. the family $\Phi$ of
flat functors (not necessarily localizations)
$\phi^*_{\lambda,\lambda'} : {\cal B}_\lambda \rightarrow {\cal
B}_{\lambda'}$ for each pair $\lambda,\lambda'$, where $\Phi$
satisfies some cocycle condition.
\end{itemize}

\ppta In 1), to ensure the faithfulness of ${\bf Q}_*$
we require that the family  $\{Q^*_\lambda\}_{\lambda \in \Lambda}$
is a flat {\bf cover} of ${\cal A}$.
That means that this is a small {\it flat} family
of functors with domain ${\cal A}$ which is
{\bf conservative}\index{conservative family}, i.e. a morphism
$f \in {\cal A}$ is invertible iff $Q^*_\lambda(f)$ is invertible
for each $\lambda \in \Lambda$.
A flat map \index{flat morphism}
whose direct image functor is conservative is called
\wind{almost affine}. In particular, this is true for
adjoint triples $f^* \dashv f_* \dashv f^!$
coming from a map $f : R \rightarrow S$ of rings.
Adjoint triples where the direct image functor is
conservative are called {\bf affine morphisms}.

\pptb \ldef{conservative} In 2), we {\it a posteriori} construct
${\cal A}$ to be ${\cal B}_\Lambda$ as before but
equipped with functors $Q^*_{\lambda}\dashv Q_{\lambda *}$
constructed from $\phi$-functors.
The cocycle condition for gluing morphism is equivalent to
the associativity of the associated comonad~(\cite{Borceux}).
The remaining requirements are made to ensure
that the comparison functor is an equivalence and
the other original data may be reconstructed as well.
The Eilenberg-Moore category of the associated
monad may be constructed directly from gluing morphisms,
and it appears to be just a reformulation of the descent category.
In a generalization,
the category ${\cal B}_{\lambda\lambda'}$ which is the essential image
of $\phi_{\lambda\lambda'}$ in ${\cal B}_{\lambda'}$ may be
replaced by any 'external' category ${\cal B}_{\lambda\lambda'}$,
but then, instead of $\phi^*_{\lambda\lambda'}$ one
requires not only flat functors
$\phi^*_{\lambda\lambda'} = : \phi^{\lambda *}_{\lambda\lambda'} :
{\cal B}_\lambda \rightarrow {\cal B}_{\lambda\lambda'}$, but
also flat functors
$\tilde{\phi}^*_{\lambda\lambda'} := \phi^{\lambda' *}_{\lambda\lambda'}:
{\cal B}_{\lambda'} \rightarrow {\cal B}_{\lambda\lambda'}$.
This generalization is essentially
more general only if we allow
the direct image functors 
(of the second type, i.e. $\phi^{\lambda'}_{\lambda\lambda'*}$),
to be not necessarily fully faithful
(hence ${\cal B}_{\lambda\lambda'}$ may not be viewed as
a full subcategory of ${\cal B}_{\lambda'}$).
Another generalization of this descent situation,
which can be phrased as having a pseudofunctor \index{pseudofunctor}
from a finite poset $\Lambda$
(viewed as a 2-category with only
identity 2-cells) to the 2-category of categories has been
studied by {\sc V.~Lunts}~(\cite{Lunts:def}).
This analogue of the descent category is called
a {\it configuration category}\index{configuration category}.
\vskip .012in

\pptc The usual formalism of descent is via fibred categories,
cf.~\cite{Vistoli:descent}. For the correspondence between the two
formalisms see, e.g.~\cite{Ros:NcSS}.
\vskip .022in

\ppt {\bf Globalization lemma.} (version for Gabriel filters:
in \cite{Rosen:88} p. 103)\ldef{glob}
{\it Suppose  $\{Q_\lambda^* : R-{\rm Mod}
\rightarrow {\cal M}_\lambda \}_{\lambda \in \Lambda}$
is a finite cover of $R-{\rm Mod}$ by flat localization
functors (e.g. a conservative family of Ore localizations
$\{S_\lambda^{-1}R \}_{\lambda \in \Lambda}$).
Denote $Q_\lambda := Q_{\lambda *}Q_\lambda^*$ where $Q_{\lambda *}$,
is the right adjoint to $Q_\lambda^*$.
Then for every left $R$-module $M$ the sequence
$$
0 \to M  \stackrel{\iota_{\Lambda,M}}\longrightarrow
\prod_{\lambda \in \Lambda} Q_\lambda M
\stackrel{\iota_{\Lambda\Lambda,M}}\longrightarrow
\prod_{(\mu,\nu) \in \Lambda\times \Lambda}
Q_{\mu}Q_\nu M
$$
is exact, where $\iota_{\Lambda M} : m \mapsto \prod \iota_{\lambda,M}(m)$,
and
}
\[ \iota_{\Lambda\Lambda M} := \prod_\lambda m_\lambda \mapsto
\prod_{(\mu,\nu)} (\iota^\mu_{\mu,\nu,M} (m_\mu) -
\iota^\nu_{\mu,\nu,M} (m_\nu)). \]
Here the order matters: pairs with $\mu = \nu$
may be (trivially) skipped, but,
unlike in the commutative case, {\bf we can not
confine to the pairs of indices with $\mu < \nu$ only}. {\it Nota bene!}

{\it Proof.} A direct corollary of Barr-Beck lemma. For
proofs in terms of Gabriel filters and torsion 
see~\cite{Rosen:book}, pp.~23--25,
and~\cite{JaraVerschorenVidal,Versch:96,vOyst:assalg}.\vskip .02in

\ppt \ldef{idempmonad}
A monad ${\bf T} = (T,\mu,\eta)$ in ${\cal A}$ is 
{\bf idempotent} if the multiplication
$\mu : TT \Rightarrow T$ is an equivalence of endofunctors.
As $\mu_M$ is the left inverse of $\eta_{TM}$, and of $T(\eta_M)$,
then $\mu_M$ is invertible iff any of them is, hence both,
and then $\eta_{TM} = T(\eta_M) = \mu^{-1}_M$.

If $\nu : TM \rightarrow M$ is a ${\bf T}$-action,
then by naturality
$\eta_M \circ \nu = T(\nu) \circ \eta_{TM} = T(\nu) \circ T(\eta_M)
= T(\nu \circ \eta_M) = {\rm Id}_{TM}$, hence $\eta_M$ is 2-sided
inverse of $\nu$ in ${\cal A}$,
hence {\it every ${\bf T}$-action is an isomorphism}.
Conversely, If every ${\bf T}$-action is an isomorphism,
$\mu_M$ is in particular, and ${\bf T}$ is idempotent.
Moreover, if every action $\nu : TM \rightarrow M$ is an isomorphism,
then its right inverse $\eta_M$ must be the 2-sided inverse, hence
there may be {\it at most} one action on a given object $M$ in ${\cal A}$.
By naturality of $\eta$, its inverse $\nu$ which is well-defined on the
full subcategory of ${\cal A}$ generated by objects in the image of
$U^{\bf T} : {\cal A}^{\bf T}\rightarrow {\cal A}$, is also natural, i.e.
it intertwines the actions, hence it is in fact a
morphism in ${\cal A}^{\bf T}$,
hence the forgetful functor $U^{\bf T} :{\cal A}^{\bf T}\rightarrow {\cal A}$
is not only faithful but also full. Its image is strictly full
(full and closed under isomorphisms) 
as the existence
of ${\bf T}$-actions on $M$ depends only on the isomorphism class of
$M$ in ${\cal A}$. To summarize, 
${\cal A}^{\bf T}$ includes via $U ^{\bf T}$ as a
strictly full subcategory into ${\cal A}$
and the inclusion has a left adjoint $F^{\bf T}$.

In general, a {\bf (co)reflective subcategory}
${\cal B} \hookrightarrow {\cal A}$ is a strictly full subcategory,
such that the inclusion $U : {\cal B}
\hookrightarrow {\cal A}$ has a (right) left adjoint, say $F$. As
$F^{\bf T}\dashv U^{\bf T}$, we have just proved that ${\cal
A}^{\bf T}$ is canonically isomorphic to a reflective subcategory
of ${\cal A}$ via inclusion $U^{\bf T}$ if ${\bf T}$ is
idempotent. On the other hand, it may be shown that for any
reflective subcategory $U : {\cal B} \hookrightarrow {\cal A}$ the
corresponding monad $(UF, U(\epsilon_F), \eta)$ is idempotent and
the comparison functor ${\cal K} : {\cal B} \cong {\cal A}^{\bf
T}$ is an isomorphism. Similarly, coreflective subcategories are
in a natural correspondence with idempotent comonads.

\minsection{Distributive laws and compatibility}

\ppt A {\bf distributive law} \index{distributive law}
{\bf from a monad  ${\bf T} = (T, \mu^T, \eta^T)$ to an endofunctor
$P$} is a natural transformation $l : TP \Rightarrow PT$ such that
\begin{equation}\label{eq:distr0}
l \circ (\eta^T)_P = P(\eta^T),\,\,\,\,\,\,\,\,\,
l \circ (\mu^T)_P = P(\mu^T) \circ l_T \circ T(l).
\end{equation}
Then $P$ {\bf lifts}
to a unique endofunctor $P^{\bf T}$ in ${\cal A}^{\bf T}$,
in the sense that $U^{\bf T} P^{\bf T} = P U^{\bf T}$.
Indeed, the endofunctor $P^{\bf T}$
is given by $(M,\nu) \mapsto (PM,P(\nu)\circ l_M)$.

\ppta A {\bf distributive law}\index{distributive law}
{\bf from a monad  ${\bf T} = (T, \mu^T, \eta^T)$ to a monad
${\bf P} = (P, \mu^P, \eta^P)$ in ${\cal A}$} (\cite{Beck:distr})
(or ``of ${\bf T}$ over ${\bf P}$'')
is a distributive law from ${\bf T}$ to the endofunctor $P$,
compatible with $\mu^P,\eta^P$ in the sense that
\[ l \circ T(\eta^P) = (\eta^P)_T,\,\,\,\,\,\,\,\,\,
l \circ T(\mu^P) = (\mu^P)_T \circ P(l) \circ l_P.\]
For clarity, we show the commutative diagram for one of the relations.
$$\xymatrix{
TPP \ar[r]^{l_P} \ar[d]_{T(\mu^P)}& PTP \ar[r]^{P(l)}
& PPT\ar[d]^{P(\mu^T)}\\
TP \ar[rr]^l && PT
}$$
Then ${\bf P}$ lifts to a unique monad
${\bf P}^{\bf T} = (P^{\bf T},\tilde\mu,\tilde\eta)$
in ${\cal A}^{\bf T}$, such that $P^{\bf T}$ lifts $P$,
and for all $N\in {\cal A}^T$ we have
$U^{\bf T}(\tilde\eta_N) = (\eta^P)_{U^{\bf T}N}$
and $U^{\bf T}(\tilde\mu_N) = (\mu^P)_{U^{\bf T}N}$.
Indeed, such a lifting is
defined by the formulas $P^{\bf T}(M,\nu):= (PM,P(\nu)\circ l_M)$,
$\mu^P_{(M,\nu)} = \mu_M$, $\eta^P_{(M,\nu)} = \eta_M$.
On the other hand, if ${\bf P}^{\bf T}$ is a lifting of ${\bf P}$
then a distributive law $l = \{l_M\}$ is defined, namely
$l_M$ is the composition
$$TPM \stackrel{TP(\eta_M)}\longrightarrow
TPTM \stackrel{U^{\bf T}(\epsilon_{PF^{\bf T}M})}\longrightarrow PTM.$$
where $F^{\bf T} : {\cal A}\to {\cal A}^{\bf T}$ is
the free $\bf T$-module functor from~\luse{monadscomonads}.

Distributive laws from ${\bf T}$ to ${\bf P}$
are in a canonical bijective correspondence with
those monads in ${\cal A}$ whose underlying functor is $PT$,
whose unit is $\eta_P\star\eta_T$, such that
$P\eta^T : P\to PT$, and $\eta^P_{T}:T\to PT$
are triple maps and which satisfy the
{\bf middle unitary law}\index{middle unitary law}
$\mu \circ (P(\eta^T)\star\eta^P_T) = \id : PT\to PT$
(cf.~\cite{Beck:distr}).
In this correspondence, $\mu^{PT}_M : PTPTM \to PTM$ is obtained
by $\mu^{PT}_M = (\mu^P \star \mu^T) \circ P(l_{TM})$, and conversely,
$l_M$ by composition $TPM \stackrel{TP(\eta_M)}\to TPTM
\stackrel{\eta_{TPTM}}\to PTPTM \stackrel{\mu^{PT}_{M}}\to PTM$.

\pptb {\bf Distributive law from a comonad
${\bf G} = (G,\delta^G, \epsilon^G)$ to a comonad
${\bf F} = (F,\delta^F, \epsilon^F)$} is a natural
transformation $l : F\circ G\Rightarrow G\circ F$ such that
\[
\begin{array}{cc}G(\epsilon^F)\circ l = (\epsilon^F)_{G}, &
(\delta^G)_F \circ l = G(l)\circ l_G \circ F(\delta^G),\\
(\epsilon^G)_{F} \circ l = F(\epsilon^G), &
G(\delta^F) \circ l = l_F \circ F(l) \circ (\delta^F)_G.\end{array}
\]

\pptc {\bf Mixed distributive law from a monad ${\bf T}$
to a comonad ${\bf G}$}:
a natural transformation $l : TG \Rightarrow G T$ such that
\[
\begin{array}{cc} \epsilon^T \circ l = T(\epsilon), &
l\circ \mu_G  = G(\mu) \circ l_T \circ T(l),\\
l \circ \eta_G = G(\eta), &
l\circ T(\delta) = \delta_T \circ G(l) \circ l_G.
\end{array}
\]
Such an $l$ corresponds to a lifting of the comonad ${\bf G}$
to a comonad ${\bf G}^{\bf T}$ in ${\cal A}^{\bf T}$,
where $G^{\bf T}(M,\nu) = (MG, G(\nu)\circ l_M)$.

\pptd {\bf Mixed distributive law
from a comonad ${\bf G}$ to a monad  ${\bf T}$}:
a natural transformation $l : GT \Rightarrow TG$ such that
\[
\begin{array}{cc}
T(\epsilon)\circ l = \epsilon_T,&
l\circ G(\mu)  = \mu_G \circ T(l) \circ l_T\\
l \circ \eta_G = G(\eta), &
T(\delta)\circ l = l_G \circ G(l)\circ \delta_T.
\end{array}
\]
Such distributive laws are in a correspondence with liftings of
a monad $\bf T$ to a monad ${\bf T}^{\bf G}$ in ${\bf G}-{\rm Comod}$.

\ppt Examples are
abundant~(\cite{Borceux,Skoda:distr};~\cite{Street:monads},II).
In a common scenario, the objects in a category of interest
are in fact objects in a simpler 'base' category,
together with multiple extra structures, satisfying
``compatibility'' conditions between the structures,
which correspond to a fixed choice of distributive laws.
Consider, left $R$-modules and right $S$-modules, where the
base category is the category ${\bf Ab}$ of \abelian groups. The
rebracketing map $(R\otimes ?) \otimes S \mapsto R\otimes (? \otimes
S)$ gives rise to a distributive law from $R \otimes ?$ and $? \otimes S$
in ${\bf Ab}$. Thus, it induces a monad $\bf V$ with underlying functor
$V = (S\otimes ?)\circ (T\otimes ?)$.
$\bf V$-modules over $U$ are precisely $R-S$-bimodules.
Similarly, for a group (or Hopf algebra) $G$,
one can describe $G$-equivariant
versions of many standard categories of sheaves or modules
with extra structure, by considering one (co)monad
for the underlying structure and another expressing the $G$-action.

\ppt A monad $(T,\mu,\eta)$ in arbitrary
2-category (even bicategory) $\cal C$ has been studied
~(\cite{Street:monads}): $T : X \to X$ is now a 1-cell,
where $X$ is a fixed $0$-cell,
and $\mu, \eta$ 2-cells should satisfy analogous axioms as
in the ordinary case, which corresponds to ${\cal C} = {\rm Cat}$.
On the other hand, if $\cal C$ is a bicategory with a single object $X$,
it may be identified with a monoidal 1-category.
The distributive laws in that case supply a notion of
compatibility of monoids and comonoids in
an arbitrary monoidal category.
The distributive laws between monoids
and comonoids in ${\rm Vec}_\genfd$
are called {\it entwining structures}~(\cite{BrzMajid:ent}).

\ppt Let ${\bf T} := (T, \mu, \eta)$ be a monad, and
$Q^* : {\cal A}\to {\cal B}$ a localization functor. The monad
${\bf T}$ is {\bf compatible} with the localization if its underlying
endofunctor $T$ is compatible with the localization, i.e. there is
a functor $T_{\cal B} : {\cal B}\to {\cal B}$
with $Q^* T = T_{\cal B} Q^*$, cf.~\luse{compat}.
In that case, $T_{\cal B}$ is
the underlying endofunctor of a unique monad
${\bf T}_{\cal B} := (T_{\cal B}, \mu^{\cal B}, \eta^{\cal B})$
in ${\cal B}$ such that  $(\mu^{\cal B})_{Q^*N} = Q^*(\mu_N)$
for every $N$ in ${\rm Ob}\,{\cal A}$.

{\it Proof.} Let $f : N\to N'$ be a morphism in ${\cal A}$,
and $g : Q^* N' \to M$ an isomorphism in ${\cal B}$.
Consider the diagram
$$\xymatrix{ T_{\cal B}T_{\cal B} Q^* N \ar[r]^=
\ar[d]^{T_{\cal B}T_{\cal B} Q^* f} &
Q^* TT N \ar[r]^{Q^*\mu_N}\ar[d]^{Q^* TT f} &
Q^* TN \ar[r]^= \ar[d]^{Q^* T f} &
T_{\cal B} Q^* N \ar[d]^{T_{\cal B} Q^* f}\\
T_{\cal B}T_{\cal B} Q^* N' \ar[r]^= \ar[d]^{T_{\cal B}T_{\cal B}g}
& Q^* TT N' \ar[r]^{Q^*\mu_{N'}} &
Q^* TN' \ar[r]^{=} & T_{\cal B} Q^* N' \ar[d]^{T_{\cal B}g}\\
T_{\cal B}T_{\cal B} M \ar[rrr]^{(\mu^{\cal B})_M} &  &
 & T_{\cal B} M.
}$$
The upper part of the diagram clearly commutes. In particular,
if $Q^* f$ is identity,  then $Q^*\mu_N = Q^* \mu_{N'}$. The vertical
arrows in the bottom part are isomorphisms, so there is a map
$(\mu^{\cal B})_M$ filling the bottom line. One has to
show that this map does not depend on the choices and that
such maps form a natural transformation. The localization functor
is a composition of a quotient functor onto the quotient category and
an equivalence. We may assume that $Q^*$ is the functor onto
the quotient category. Then, by the construction of the
quotient category, every morphism $g$ is of the zig-zag
form as a composition of the maps of the form $Q^*f$ and
formal inverses of such maps, and if $g$ is an isomorphism,
both kinds of ingredients are separately invertible in ${\cal B}$.
To show that
$(\mu^B)_M =
(T_{\cal B} g)\circ Q^*(\mu_{N'})\circ (T_{\cal B}T_{\cal B}g)^{-1}$
for every isomorphism $g : Q^*N'\to M$ is the consistent choice,
we use the upper part of the diagram repeatedly
(induction by the length of zig-zag) for the zig-zag
isomorphism $h = (g_1)^{-1}g_2: Q^*N_2\to Q^*N_1$
where $g_i : Q^* N_i \to M$.
One obtains
$T_{\cal B} h \circ Q^*\mu_{N_2}\circ (T_{\cal B}T_{\cal B}h)^{-1}
= Q^*\mu_{N_1}$, hence
$(T_{\cal B} g_2)\circ Q^*(\mu_{N_2})\circ (T_{\cal B}T_{\cal B}g_2)^{-1}
= (T_{\cal B} g_1)\circ Q^*(\mu_{N_1})\circ (T_{\cal B}T_{\cal B}g_1)^{-1}$.

Upper part of the diagram also shows the naturality for $\mu^{\cal B}$
with respect to each arrow of the form $Q^*f$ and with respect to
formal inverses of such. For any morphism $h : M\to M'$ in ${\cal B}$,
using its zig-zag representation, we extend this to
the naturality diagram $(\mu^{\cal B})_{M'} \circ T_{\cal B}T_{\cal B}h
=  T_{\cal B}T_{\cal B}h \circ (\mu^{\cal B})_M$.
Uniqueness of $\mu^{\cal B}$ is clear by
the requirement $(\mu^{\cal B})_{Q^* N} = Q^*(\mu_N)$ and the naturality.
The unit morphism $\eta : 1_{\cal B}\to T_{\cal B}$ satisfies
 $g \circ \eta^{\cal B}_M = T_{\cal B}(g) \circ Q^*(\eta_{N})$
for every isomorphism $g : Q^*N \to M$ in ${\cal B}$
such that $N \in {\rm Ob}\,{\cal A}$.
In particular, $(\eta^{\cal B})_{Q^*N} = Q^*(\eta_N)$.
The very axioms of a monad may be checked in a similar
vein.    

\pptb If $\comdG = (\comdg,\delta,\epsilon)$ is a comonad and
the endofunctor $\comdg$ is compatible with each localization in family
$\{Q^*_\lambda : {\cal A}\to {\cal B}_\lambda\}_{\lambda \in \Lambda}$,
then there is a unique family of comonads
$\{\comdG_\lambda = (\comdg_\lambda,\delta_\lambda,\epsilon_\lambda)
\}_{\lambda \in \Lambda}$
such that $Q_\lambda^* G = G_\lambda Q^*_\lambda$ for each $\lambda$.
We have then $(\delta_{\lambda})_{Q^*_\lambda M} = Q^*_\lambda (\delta_M)$
and $(\epsilon_\lambda)_{Q^*_\lambda M} = Q^*_\lambda (\epsilon_M)$
for every $M \in {\rm Ob}\,{\cal A}$.

\ppt If $Q^*: {\cal A}\to {\cal B}$ is a localization
with right adjoint $Q_*$, and $T'$ is an endofunctor in ${\cal B}$,
then $Q_* T' Q^*$ is compatible with $Q^*$.
Indeed $\epsilon$ is
an isomorphism by~\luse{absloc}, hence
$\epsilon_{T' Q^*} : Q^* Q_* T' Q^*\Rightarrow T' Q^*$
is an isomorphism, and the assertion follows by~\luse{compat}.\vskip .015in

\ppt {\bf Example from Hopf algebra theory}.
Let $B$ be a $\genfd$-bialgebra and $(E,\rho)$
a right $B$-comodule together with
a multiplication $\mu : E \otimes_\genfd E \to E$ making it
a $B$-comodule algebra, i.e. an algebra
in the category of right $B$-comodules.
The coaction $\rho : E \to E\otimes H$ is
{\bf compatible} with a fixed Ore
localization $\iota : E \to S^{-1}E$ if there is a
coaction $\rho_S : S^{-1}E \to S^{-1}E\otimes H$
which is an algebra map and such that
$\rho_S \circ \iota = (\iota \otimes \id_B)
\circ \rho$.  $B$ induces a natural comonad $T = T(B,E)$
in $E-{\rm Mod}$, such that $(E-{\rm Mod})^T$ is a
category of so-called $(E,B)$-Hopf modules. The compatibility
above ensures that the localization lifts to a localization
of $(E-{\rm Mod})^T$~(\cite{Skoda:ban}, 8.5).
Hence, $T$ is compatible with the localization in the usual sense,
with numerous applications of this type of situation
~(\cite{Skoda:coh-states,Skoda:hg,Skoda:locco,Skoda:ban}).\vskip .015in

\ppt The compatibility of certain localizations of noncommutative
spaces with differential functors is central in the
treatment~\cite{LuntsRosMP,LuntsRos1,LRloc1,LRloc2} of ${\cal
D}$-modules on noncommutative spaces.\vskip .01in

\ppt Distributive laws become a much simpler issue when both
monads in question are idempotent in the sense
of~\luse{idempmonad}. In the localization literature this is
roughly the  situation treated under the name of {\bf mutual
compatibility of localizations}.

Let $S,W$ be two Ore sets. Then the set $SW$ of products
$\{ sw | s \in S, w \in W\}$ is not necessarily
multiplicatively closed.

Suppose $SW$ is multiplicatively closed. This means that $\forall
s \in S$, $\forall w\in W$ if the product $ws$ is in $SW$ then
$\exists w'\in W$, $\exists s'\in S$ such that $ws = s' w'$.
Suppose now $M := {}_R M\in R-{\rm Mod}$. Each element in
$S^{-1}R\otimes_R W^{-1}R \otimes_R M$ is of the form
$s^{-1}\otimes w^{-1}\otimes m$ with $s \in S$, $w \in W$ and $m
\in M$. By a symmetric argument, $ws \in SW$,
and thus $\exists w'\in W$, $\exists
s'\in S$ such that $ws = s' w'$. Choosing $s', w'$ by this rule
we obtain an assignment $s^{-1}\otimes w^{-1}\otimes m \mapsto
(w')^{-1}\otimes (s')^{-1}\otimes m$. We claim that this
assignment is well-defined and a map of left $R$-modules
$S^{-1}R\otimes_R W^{-1}R \otimes_R M \to W^{-1}R\otimes_R S^{-1}R
\otimes_R M$. As $M$ runs through $R-{\rm Mod}$, such maps form a
natural transformation $Q_S Q_W \to Q_W Q_S$ of functors, which is
clearly an isomorphism.

In fact, this natural transformation is a distributive law.
Although the compatibility of $Q_S$ and $Q_W$ is symmetric, the
converse does not hold: compatibility does not mean that $SW$ is
multiplicatively closed. Indeed, let $R$ be a $\DDl C$-algebra
with two generators $a,b$ and relation $ab = qba$ where $q \neq
\pm 1, 0$. Then the set multiplicatively generated by $A$ and set
multiplicatively generated by $B$ are 2-sided Ore sets, and the
corresponding localization functors are compatible; however $AB$
is not multiplicatively closed.

If $S, W$ are left Ore in $R$, that does not mean that
$\iota_W(S)$ is left Ore in $W^{-1}R$. Namely, the left Ore
condition for (the image of) $S$ in $W^{-1}R$ includes the
following: $\forall s\in S, \forall t\in T, \exists s'\in S,
\exists (w')^{-1}r' \in W^{-1} R, (w')^{-1}r' s = s' w^{-1}$. If
$R$ is a domain, this means that $r' sw = w's'$. This is almost
the same condition as that $SW$ is multiplicatively closed
(above), except that one can choose extra $r'$. In the same away
as in the former case, we derive the compatibility of $Q_S$ and
$Q_W$. If we change left Ore sets to right Ore sets, or $S$ being
Ore in $W^{-1}R$ to $W$ being Ore in $S^{-1}R$ we get similar
``Ore'' equations $swr = w's'$, $wsr = s'w'$ etc. From the
abstract point of view (say torsion theories) these
compatibilities are indistinguishable.

The compatibility implies that the localization at
the smallest multiplicative set
generated by $S$ and $W$ is isomorphic to the consecutive
localization by $S$ and then by $W$.
This simplifies the formalism of localization
(cf. semiseparated schemes and \v{C}ech resolutions
of~\cite{Ros:NcSch}, cf.~\luse{glob},
\cite{JaraVerschorenVidal} etc.).

\minsection{Commutative localization} \label{sec:commloc}

Here we describe specifics in the commutative case, and further
motivation from commutative algebraic geometry, and its
abstractions.

\ppt Suppose $R$ is a unital associative ring, $Z(R)$ its
center, and $S \subset Z(R)$ a multiplicative subset. Obviously,
$S$ is automatically a left and right Ore subset in $R$, with
simpler proofs for the construction and usage of the Ore
localization. We say that $S^{-1}R = S \times R/\sim$ is the
\wind{commutative localization} of $R$ at $S$. {\it The
equivalence relation $\sim$ (\luse{relFrac}) simplifies to}
\begin{equation}\label{eq:relComFrac}
\fbox{$s^{-1} r \sim s'^{-1} r' \,\Leftrightarrow\,
\exists \tilde{s} \in S, \,\, \tilde{s} (s r' - s' r) = 0.$}\end{equation}
Proof. By the definition, $\exists \tilde{s}\in S$,
$\exists \tilde{r} \in R$,
such that $\tilde{s} r = \tilde{r} r'$ and $\tilde{s} s = \tilde{r} s'$.
Therefore,
\[ \tilde{s} sr' = \tilde{r} s' r' = \tilde{r} r' s' = \tilde{s} r s' =
\tilde{s} s' r.\,\,\,\,\,\,\,\,\Box \]

Unlike sometimes (mis)stated in the literature
(e.g.~\cite{vOyst:assalg},p.14), the commutative
formula~(\ref{eq:relComFrac}) (and variants of it) is
inappropriate even for mildly noncommutative rings and even
2-sided Ore sets which are not in center. E.g. take the unital
${\DDl C}$-algebra generated by two elements $b$ and $d$ with $bd
= qdb$, where ${\DDl C} \ni q \neq 1$. That algebra has no zero
divisors. Let $S$ be the 2-sided Ore set multiplicatively
generated by $b$ and $d$. Formula $b^{-1} = (db)^{-1} d$, and the
criterion above would imply that $db=bd$ with contradiction.

For general $R$ and $S$, formula~(\ref{eq:relComFrac})
is actually not even an equivalence relation on $R\times S$.

\ppt \ldef{genschemes} {\bf General requirements on scheme-like
theories.} One wants to mimic several major points from the
classical case. We first decide which geometric objects constitute
the category ${\cal C}$ of affine schemes; then find a suitable
larger geometric category ${\frak Esp}$ of spaces, in the sense
that it is equipped with a fully faithful functor ${\cal C}
\hookrightarrow {\frak Esp}$, where the objects in the image will
be called here {\it geometric affine schemes};\index{geometric
affine schemes} finally there is a gluing procedure which assigns
to a collection $\{{\frak C}_h\}_{h \in H}$ of geometric affine
schemes and some additional 'gluing' data $Z$, a space which may
be symbolically denoted by $\left(\coprod_h {\frak C}_h \right)/ Z
\in {\frak Esp}$, together with canonical morphisms ${\frak C}_h
\rightarrow \left(\coprod_h {\frak C}_h \right)/ Z$ in ${\frak
Esp}$. For a fixed (type of) gluing procedure \index{gluing
procedure} ${\cal G}$, and including all the isomorphic objects,
one constructs this way a subcategory of {\bf locally affine
spaces} of type $({\cal C}, {\cal G})$ in ${\frak Esp}$.

Additional requirements and typical choices are in place.

\ppta {\bf Abstract affine schemes.}\index{abstract affine scheme}
Most often one deals with some monoidal category
$\tilde{\cal A} = ({\cal A},\otimes, 1_{\tilde{A}})$.
Then ${\cal C} = {\rm Aff}(\tilde{A})
:= ({\rm Alg}(\tilde{A}))^{\rm opp}$ \index{${\rm Aff}(\tilde{A})$}
is the {\bf category of affine schemes in ${\cal A}$} i.e.
the {\it opposite to the category of algebras (monoids) in $\tilde{A}$}.
The basic example is the monoidal category of $R$-bimodules, where
$R$ is a $\genfd$-algebra over a commutative ring $\genfd$. The affine
schemes in this category are given by $\genfd$-algebra maps
$f : R \rightarrow R'$ (making such $R'$ an $R$-bimodule).

If ${\cal A}$ is {\it symmetric} \index{monoidal
category!symmetric} via a symmetry $\tau$ where $\tau_{AB} : A
\otimes B \rightarrow B \otimes A$, then one may consider only
$\tau$-commutative algebras $(A,\mu,\eta)$, i.e. for which $\mu
\circ \tau = \mu$. $\tau$-{\bf commutative affine schemes} in
${\cal A}$ are the objects of the opposite of the category of
$\tau$-commutative algebras
 ${\cal C} = {\rm cAff}(\tilde{A},\tau)
:= ({\rm cAlg}(\tilde{A},\tau))^{\rm opp}$ \index{${\rm cAff}$}
(one often skips $\tau$ in the notation). Examples are
(super)commutative affine schemes in a $\otimes$-category of
$\genfd$-modules and also the opposite to the category
$\underline{\rm cdga}_\genfd$ of commutative differential graded
$\genfd$-algebras, which is important in recent 'derived algebraic
geometry' \index{derived algebraic geometry} program
~(\cite{Toen:higher,ToenVez:mixt}).

\pptb {\bf Gluing for ringed topological spaces} (and a version
for local (l.) rings). A (locally) ringed space is a pair $(X,F)$
consisting of a topological space $X$, and a ('structure') sheaf
$F$ of (l.) rings on $X$. A morphism of (l.) ringed spaces is a
continuous map $f : X \to X'$ with a {\it comorphism} i.e. a map
of sheaves of (l.) rings $F'\to f^*X'$ over $X'$. We obtain a
category ${\frak rSp}$ (${\frak lSp}$). \index{${\frak
rSp}$}\index{${\frak lSp}$} Given a full subcategory ${\cal A}$ of
${\frak rSp}$, considered as a category of (heuristic term here)
``local models''  we may consider all (l.) ringed space $X$ for
which there is a cover (in usual sense) $X = \cup_\alpha X_\alpha$ of
underlying topological spaces and for each $\alpha$ an isomorphism
$X^0_\alpha \cong X_\alpha$ in ${\frak rSp}$ (${\frak lSp}$) where
$X^0_\alpha$ is in ${\cal A}$. More abstractly, but equivalently,
consider all families of morphisms $\{i_\alpha : X_\alpha \to
X\}_\alpha$, which are 'covers by embeddings': topologically
covers of $X$ by families of monomorphisms (continuous, open and
injective), and sheaf-wise isomorphisms on stalks. The
intersections $X_\alpha \cap X_\beta$ represent fibred products
$X_\alpha \times_X X_\beta$ in ${\frak rSp}$. Spaces glued from
objects in ${\cal A}$ are nothing but the colimits of the diagrams
of the type $\xymatrix{\coprod_{\alpha\beta} Y_{\alpha\beta}
\ar[r]<-2pt> \ar[r]<2pt> &\coprod Y_\alpha}$ where each of the
morphisms $Y_{\alpha\beta}\to Y_\alpha$ and $Y_{\alpha\beta}\to
Y_\beta$ are embeddings, and $Y_\beta \in {\rm Ob}\,{\cal A}$.
There is a natural condition on ${\cal A}$: each $X$ in ${\cal A}$
as a topological space has a basis of topology made out of (spaces
of) some family of objects in ${\cal A}$ (or isomorphic to them);
and this family may be chosen so that the restrictions of the
structure sheaf agree. In the theory of schemes, affine schemes
are such models: affine subschemes make a basis of topology, but
not every open subset is affine; nor their intersections. Still
the intersections and colimits exist from the start, in our
ambient category of ringed spaces which is big enough.

\pptc {\bf Grothendieck (pre)topologies}~(G.(p)t.)
\index{Grothendieck (pre)topology} and gluing in commutative
algebraic geometry~(\cite{MLMoerd:Sheaves,Vistoli:descent}).
Schemes are glued from affine schemes in the Zariski topology
(which may be considered both as an ordinary topology and a
G.(p)t.); useful generalizations (e.g. algebraic spaces) in flat
and \'{e}tale G.(p)t. etc.
A \wind{sieve} is an assignment of a collection $J(R)$ of
morphisms in ${\cal C}$~(\cite{MLMoerd:Sheaves}) with target $R$,
such that if the target of $f$ is $J(R)$ and $g : R\to R'$ then
$g\circ f \in J(R')$. A G.t. is a collection of sieves $J(R)$ with
target $R$ for each $R$ in ${\rm Ob}\,{\cal C}$, such that if $f
\in J(R)$ and $g : R\to R'$ then $g\circ f \in J(R')$), satisfying
some axioms~(\cite{MLMoerd:Sheaves}). A G.pt. in a category ${\cal
C}$ with fibred products is a class ${\cal T}$ of families $\{
U_\alpha \to U\}_\alpha$ of morphisms (with one target per
family), such that $\{ {\rm Id} : U\to U\} \in {\cal T}$; if $\{
f_\alpha : U_\alpha \to U\}_\alpha \in {\cal T}$ and $\forall
\alpha$ $\{ g_{\alpha\beta} : U_{\alpha\beta} \to
U_\alpha\}_\beta\in {\cal T}$ then $\{f_\alpha \circ
g_{\alpha\beta}\}_{\alpha,\beta} \in {\cal T}$; and finally if $\{
f_\alpha : U_\alpha \to U\}_\alpha \in {\cal T}$ and $g : V \to U$
is a morphism then $\{g^*(f_\alpha) : V\times_U U_\alpha \to V\}
\in {\cal T}$. Elements of ${\cal T}$ are called {\bf covers} (in
${\cal T}$), and the pair $({\cal C},{\cal T})$ a site.
\index{site} To any (ordinary) topological space $X$ one
associates a ``small'' site ${\frak Ouv}_X$: objects are open
subsets in $X$; morphisms are inclusions; $\{U_\alpha
\hookrightarrow U\}_\alpha$ is a cover if $\cup_\alpha U_\alpha =
U$.

A presheaf $F$ of sets on $X$
is a functor $F : ({\frak Ouv}_X)^{\rm op}\to{\rm Sets}$;
a presheaf $F$ on any site $({\cal C},{\cal T})$
with values in a category ${\cal D}$ with products is a functor
$F : {\cal C}^{\rm op}\to {\cal D}$. Given a cover
$\{ U_\alpha \to U\}_\alpha \in {\cal T}$, there are two
obvious embeddings $\xymatrix{\prod U_\alpha \times_U U_\beta
\ar[r]<-2pt> \ar[r]<2pt> &\prod_\alpha U_\alpha }$.
A presheaf is a sheaf on
$({\cal C},{\cal T})$ if $\{ U_\alpha \to U\}_\alpha \in {\cal T}$
if for every such cover the induced diagram
$\xymatrix{ F(U) \ar[r] & \prod_\alpha F(U_\alpha)
\ar[r]<-2pt> 
\ar[r]<2pt> 
&  \prod_{\alpha\beta}
F(U_\alpha \times_U U_\beta)}$
is an equalizer diagram. 

For gluing, one again needs some bigger ambient category (or,
instead, some universal construction). Our local models are now
(commutative) affine schemes ${\rm Aff} := {\rm Aff}({\rm Ab})$
with a G.pt. ${\cal T}$. The Yoneda embedding $X \mapsto X\hat{}
:= {\rm Aff}(?,X)$ is a fully faithful functor from ${\rm Aff}$
into the category ${\rm PFas}\,({\rm Aff})$ of presheaves of sets
on ${\rm Aff}$. One typically deals with {\it subcanonical G.t.}
which means that the presheaves in the Yoneda image (representable
functors) are sheaves. As in the case of ordinary topologies, to
construct the global locally ${\cal T}$-affine spaces, one needs
colimits of certain diagrams of the form
$\xymatrix{\prod_{\alpha\beta} V_{\alpha\beta} \ar[r]<-2pt>
\ar[r]<2pt> &\prod_\alpha U_\alpha }$, where, in the addition, the
colimit cone $\prod U_\alpha\to U$ corresponds to a cover in
${\cal T}$. As for example, the nonseparated schemes in Zariski
topology, some locally ${\cal T}$-affine spaces may not be
produced this way with $V_{\alpha\beta}$ being in ${\rm Aff}$.
Similar problems for biflat covers by localizations in
noncommutative geometry are known (lack of compatibility of
localizations; nonsemiseparated covers). Furthermore, one needs to
extend the notion of ${\cal T}$-covers to the target category of
sheaves. We hope that the reader sees at this point the meaning of
this abstract machinery. We won't proceed with the full
construction. Namely, in the commutative case, it is usually
replaced by equivalent constructions. For example, to construct
the algebraic spaces, one usually does not glue affine schemes
``over intersections'' but rather starts with an equivalence
relation in the category of all separated schemes.
In noncommutative case, G.t. are elaborated in~\cite{KontsRos:MP, Ros:NcSS}
and partly in~\cite{Ros:USpSch}. There is also an approach to G.t. and
quasicoherent sheaves for noncommutative affine schemes by {\sc
Orlov}~(\cite{orlov:qcsh}), utilizing ringed sites and sheaves of
groupoids, but implicit application to the construction of noncommutative
locally affine spaces is not given.

\ppt Though we assume that the reader has been exposed to the
commutative scheme theory, we here sketch the basic construction
of an affine scheme, as a part of a widely-known easy
generalization to pairs of the form
(noncommutative ring, central subring), cf.~\cite{semiquantum}.

A  left ideal ${\frak p} \subset R$
is {\bf prime} \index{prime ideal}
if for any two left ideals $I,J$,
if $IJ \subset {\frak p}$ then
$I\subset {\frak p}$ or $J\subset {\frak p}$.
A left ideal ${\frak p} \subset R$ is
{\bf completely prime} \index{completely prime ideal}
if $fg \in {\frak p}$ then $f \in {\frak p}$ or $g \in {\frak p}$;
equivalently $R/{\frak p}$ is a domain; or $R\backslash {\frak p}$
is multiplicative. Each completely prime ideal is
prime: otherwise one could find $f \in I\backslash{\frak p}$
and $g \in J \backslash {\frak p}$, such that $fg \in IJ\subset {\frak p}$
with contradiction.
If $R$ is commutative the converse holds, as one can see by
specializing the definitions to the principal left ideals
$I = Rf$, $J = Rg$, $IJ = RfRg = Rfg$.

Consider the category ${\cal R}$,
whose objects are pairs $(R,C)$, of a unital ring $R$
and a central subring $C \subset Z(R)$; and
morphisms $(R,C) \to (R',C')$ are maps of rings $\phi: R\to R'$
such that $\phi(C) \subset C'$.

${\rm Spec}\,C$ is the set of all prime ideals $\frak p$ of $C$.
For any ideal $I \subset C$, define $V(I)\subset {\rm Spec}\,C$
as the set of all ${\frak p} \subset C$, such that $I \subset {\frak p}$.
Sets of the form $V(I)$ depend only
(contravariantly with respect to inclusions)
on the radical $\sqrt{I}$ (the intersection
of all prime ideals containing $I$) and satisfy the axioms of antitopology.
Complements of such sets hence form a topology on ${\rm Spec}\,C$,
called \wind{Zariski topology}.
{\bf Principal open sets} \index{principal open set}
are the sets of the form $U_f = V((f))$,
where $(f)$ is the (principal) ideal generated by $f \in C$.
They make a basis of Zariski topology,
i.e., any Zariski open set is a union of sets of that form.
Open sets and inclusions form category ${\frak Ouv}_C$.
Principal open sets and inclusions
form its full subcategory ${\frak OuvP}_C$.

Define ${\cal O}'_C(U_f):= C[f^{-1}]$ and ${\cal O}'_{R,C}(U_f):= R[f^{-1}]$.
Every inclusion $U_f \hookrightarrow U_g$ induces the unital ring maps
$\phi_{f,g,i}: {\cal O}'_i(U_g) \to {\cal O}'_i(U_f)$, $i \in \{C, (R,C)\}$.
Hence we have contravariant functors
${\cal O}_i' : {\frak OuvP}_C\to \catRings$.
Natural inclusions ${\bf in}_f : C[f^{-1}]\rightarrow R[f^{-1}]$ form
a natural transformation, i.e.,
${\cal O}'_C$ is a subfunctor of ${\cal O}'_{C,R}$.
Functors ${\cal O}'_i$ extend naturally to functors
${\cal O}_i : {\frak Ouv}_C \to \catRings$ which are sheaves, and this
requirement fixes sheaves ${\cal O}_i$ uniquely up to isomorphism of sheaves.
Namely, represent any open set $U$ as the union $\cup U_f$ of
a family ${\cal U}$ of (some or all) $U_f \subset U$.
Define a diagram $\Delta({\cal U})$ in $\catRings$ as follows.
Vertices are of the form ${\cal O}_i(U_f)$ and of the form
${\cal O}_i(U_{fg}) = {\cal O}_i(U_f \cap U_g)$
where $U_f, U_g \in {\cal U}$, and the arrows
are $\phi_{f,fg,i}$. Then ${\cal O}_i(U)$
has to be isomorphic to the inverse limit of the diagram
$\Delta({\cal U})$ and may be consistently set so.
Moreover the natural transformation ${\bf in}$ extends making
${\cal O}_C$ a subfunctor of ${\cal O}_{C,R}$. In fact,
${\cal O}_{C,R}$ is an algebra in the monoidal category of (2-sided)
${\cal O}_C$-modules. All the stalks $({\cal O}_C)_{\frak p}$
of ${\cal O}_C$ are local rings, namely the localizations
$C_{\frak p}:= C[\{f^{-1}\}_{f \notin {\frak p}}]$
'at prime ideal ${\frak p}$'.

Let  $\psi : R \rightarrow R'$ be a map of unital rings.
The inverse image $\psi^{-1}({\frak p})$
of a completely prime ideal is completely prime.
Let $\psi : (R,C) \to (R',C')$ be a morphism in ${\cal R}$,
and let map $\phi : C \to C'$ agrees with $\psi$ on $C$.
One has map $\phi^* : {\rm Spec}\,C' \to {\rm Spec}\,C$
given by $\phi^* : {\frak p} \mapsto \phi^{-1}({\frak p})$.
If $U \subset {\rm Spec}\,C'$ is open, then
$\phi^{-1}(U) \subset {\rm Spec}\,C$
is open as well, because $\phi^*(V(I)) = V(\phi^{-1}(I))$ for each ideal
$I \subset C'$. Hence, $\phi^*$ is continuous.
If $g \notin \phi^{-1}({\frak p})$ then $\phi(g) \notin {\frak p}$. Thus
all elements $g \in C$, newly inverted in $C_{\phi^{-1}({\frak p})}$
(and $R_{\phi^{-1}({\frak p})}$)
are also invertible in $C'_{\frak p}$ (and $R'_{\frak p}$).
Hence, by the universality of localization, one has a unique map
$\psi_{\frak p} : (R_{\phi^{-1}({\frak p})}, C_{\phi^{-1}({\frak p})})
\to (R'_{\frak p}, C'_{\frak p})$
such that
$\psi_{\frak p} \circ \iota_{\phi^{-1}({\frak p})}
= \iota_{\frak p} \circ\psi$.
Define $\phi_U :  {\cal O}_{R,C}(\phi^* U) \rightarrow {\cal O}_{R',C'}(U)$
by $(\phi_U(r))_{\frak p} = \psi_{\frak p} (r_{\phi^{-1}({\frak p})})$.
One has to check that $\phi_U(r)$ is indeed in ${\cal O}_{R',C'}(U)$, i.e.,
that ${\frak p}\mapsto  \psi_{\frak p} (r_{\phi^{-1}({\frak p})})$ is
indeed a section. For this consider $U_f \subset U$ affine, i.e.,
$U_f = \{ {\frak p} \,|\, f \notin {\frak p} \}$ and
$\phi^* U_f := \{\phi^{-1}({\frak p}),{ {\frak p}\not\ni f}\}$.
An argument as above gives map
$\psi_f : {\cal O}_{R,C}(\phi^* U_f) \rightarrow {\cal O}_{R',C'}(U_f)$
satisfying $\psi_{f} \circ \iota_{\phi^{-1}(U_f)}
= \iota_{U_f} \circ\psi$. It is easy to check then that $\psi_f$ induces
$\psi_{\frak p}$ in stalk over ${\frak p}\in U_f$.
As a result, we obtain a map
$\psi^\sharp : \phi_* {\cal O}_{R,C} \rightarrow {\cal O}_{R',C'}$
of sheaves over ${\rm Spec}\,C'$.

\ppta Let ${\frak lSp}_2$ be a category whose object are
locally ringed spaces $(X,{\cal O})$ in ${\frak lSp}$
together with a sheaf ${\cal O}^{\rm nc}$ of noncommutative
algebras in the category of ${\cal O}$-modules.
A morphism $\psi : (X,{\cal O}_X,{\cal O}^{\rm nc}_X) \to
(Y,{\cal O}_Y,{\cal O}^{\rm nc}_Y)$
in ${\frak lSp}_2$ is a morphism
$\psi^{\rm c} : (X,{\cal O}_X) \to (Y,{\cal O}_Y)$
in ${\frak lSp}$, together with a map of sheaves
of ${\cal O}_X$-modules
$\psi^\sharp : {\cal O}^{\rm nc}_Y\to \psi^* {\cal O}^{\rm nc}_X$.

We have above constructed a functor
${\rm Spec_2} : {\cal R} \rightarrow {\frak lSp}_2$.

\pptb In fact, ${\cal O}^{\rm nc}$ is in a smaller category of
{\it quasicoherent sheaves of ${\cal O}$-modules} (shortly: quasicoherent
modules). A presheaf $F$ of ${\cal O}$-modules on a ringed
space $(X,{\cal O})$ is quasicoherent~(EGA 0.5.1.1) if $\forall x\in X$
$\exists U^{\rm open} \ni x$ with an exact sequence
${\cal O}^I \to {\cal O}^J\to F \to 0$ where ${\cal O}^I$,
${\cal O}^J$ are free ${\cal O}$-modules (of possibly infinite rank).
If the ringed space is locally ${\cal T}$-affine for some
Grothendieck topology ${\cal T}$ on the category of
commutative affine schemes, then one may take $U$ affine,
and an equivalent definition of quasicoherence is that
for any pair of affine open subsets $W \subset V \subset X$,
\begin{equation}\label{eq:quasicoherence}
F(W) = {\cal O}(W) \otimes_{{\cal O}(V)} F(V).\end{equation} The
noncommutative structure sheaf ${\cal O}^{\rm nc}$ of ${\rm
Spec}_2(R,C)$ is a quasicoherent ${\cal O}$-module. For presheaves
of  ${\cal O}^{\rm nc}$ left modules one may use the same
formula~(\ref{eq:quasicoherence}). For bimodules one distinguishes
left and right
quasicoherence~\cite{Konts:defqalgvar,Lunts:def,Yek:def} (in the
right hand version the tensoring with ${\cal O}(W)$ in
formula~(\ref{eq:quasicoherence}) is from the RHS instead).
Formula~(\ref{eq:quasicoherence}) is nothing else but the formula
for a localization functor $Q^* : F(V)\mapsto F(W)$ from ${\cal
O}(V)$-modules to  ${\cal O}(W)$-modules. In this case, both $Q^*$
and its right adjoint $Q_*$ are exact functors. Ring theorists
call such localizations of full module categories {\bf perfect
localizations}~(\cite{JaraVerschorenVidal}). Equivalently, the
canonical forgetful functor from the localized category to the
category of modules over the localized ring is an equivalence of
categories.

{\sc F.~van Oystaeyen}~(\cite{vOyst:assalg})\ldef{hereditary}
\index{hereditary torsion theory}
defines quasicoherent presheaves on the lattice
of hereditary torsion theories (= localizations where $Q^*$ is exact
and the torsion subcategory is coreflective)
over a ${\DDl Z}_{\geq 0}$-graded rings, by using the
appropriate localization instead of tensoring. For affine
case see~\cite{Rosen:book}, I 6.0.3 (b) and I 6.2.
One has theorems on gluing of such modules using Barr-Beck lemma.

{\sc D.~Orlov} (\cite{orlov:qcsh}) defines quasicoherent (pre)sheaves
on $X$ where $X$ is a presheaf of sets (in particular, by Yoneda,
on any object) on a given ringed site $({\cal C},{\cal T})$.
\vskip .02in

\ppt \ldef{delgab}
Now we will quote two theorems. A theorem of {\sc Deligne} shows that
in the context of commutative affine schemes, a formula which can
be recognized as a  localization at a Gabriel filter (see below),
describes the behaviour of the category of quasicoherent sheaves
under passage to a not necessarily affine subset. Hence the
``noncommutative localization'' is already seen there! The
theorem may be proved directly, and we suggest to the reader
to at least convince oneself that the RHS formula is indeed an
$A$-module.
``Deligne theorem''~(\cite{Hart:ag})
in fact, cf.~(\cite{Hart:res}, Appendix)\footnote{
I thank Prof. {\sc Hartshorne} for an email on the history of the result.},
has been originally inferred from the
Gabriel theorem below --
the general statement that quasicompact open subsets of
(nonaffine) schemes correspond always to exact localizations
of \abelian categories; combined with the Gabriel's formulas on how
such localizations look like~(\cite{Gab:catab}).

\ppta {\bf Deligne's theorem.} \index{Deligne's theorem}
Let $X = {\bf Spec}\,A$ be an affine Noetherian scheme, i.e.
$A$ is a commutative Noetherian ring. Let $U$ be a Zariski open set
(not necessarily affine!), and $I$ an ideal such that
$V(I) = X \backslash U$.
Let $M$ be an $A$-module, and $\tilde{M}$ the corresponding
quasicoherent ${\cal O}_X$-module. Denote by $\tilde{M}|_U$ its
restriction to $U$.
Then
\[
  \Gamma(\tilde{M}|_U) \cong
\lim_{n \rightarrow \infty} {\rm Hom}_A (I^n, M).
\]
This isomorphism of $A$-modules is natural in $M$.

\pptb ({\sc P.~Gabriel} \cite{Gab:catab}, VI.3)
\index{localizing to open subscheme}
Let $(X, {\cal O}_X)$ be
a scheme, and $U$ an open subset of $X$, such that the canonical
injection $j : U \hookrightarrow X$ is quasicompact.
Let $Q^*_U : {\frak Qcoh}_X \rightarrow {\frak Qcoh}_U$
be a functor which associates to every
quasicoherent ${\cal O}_X$-module $M$ the restriction $M |_U$
of $M$ to subscheme $(U, {\cal O}_X |_U)$. Functor $Q^*_U$
canonically decomposes as ${\frak Qcoh}_X \rightarrow {\frak Qcoh}_X/
{{\rm Ker}\,Q^*_U} \stackrel{\cong}{\rightarrow} {\frak Qcoh}_U$
into the canonical projection onto
the quotient category by ${\rm Ker}\,Q^*_U$
and an isomorphism.

\minsection{Ring maps vs. module functors}

\ppt \ldef{functorsMod} As usual, we write ${}_R M$ when
we want to emphasize that $M$ is (understood as)
a left $R$-module; $R-{\rm Mod}$ is the category of
left $R$-modules. Let $f : R \rightarrow S$
be any map of (not necessarily unital) rings. We have the
following functors induces by map $f$:
\begin{itemize}
\item {\bf extension of scalars} \index{extension of scalars}
$f^* : R-{\rm Mod} \to S-{\rm Mod}$, $M \mapsto S \otimes_R M$;

\item {\bf restriction of scalars} \index{restriction of scalars}
(forgetful functor)
$f_* :  S-{\rm Mod} \to R-{\rm Mod}$, ${}_S M \mapsto {}_R M$;

\item $f^! : R-{\rm Mod} \to S-{\rm Mod}$,
$M \mapsto {\rm Hom}_R ({}_R S, M)$.

\end{itemize}
Denote $F \dashv G$ when functor $F$ is left adjoint
to functor $G$. Easy fact: $f^* \dashv f_* \dashv f^!$.
Hence $f^*$ is left exact, $f^!$ right exact and $f_*$ exact.

\ppt \ldef{comorph}
It is of uttermost importance to have in mind the geometrical
picture of this situation in the case when $R$ and $S$ are commutative
and unital. Denote by $\catlSp$ \index{$\catlSp$}
the category of
locally ringed spaces.
An object $(X,{\cal O}_X)\in \catlSp$ is
a pair of a topological space $X$
and a sheaf of commutative local rings ${\cal O}_X$ over $X$; and
a morphism is a pair $(f^o, f^\sharp)$
where $f^o : X \to Y$ is a map of topological spaces
and a 'comorphism'
$f^\sharp : f^\bullet {\cal O}_{Y} \rightarrow {\cal O}_{X}$
is a map of sheaves of local rings over $X$, and
$f^\bullet$ a sheaf-theoretic pullback functor.
There is a contravariant functor
${\rm Spec} : \catCommRings \rightarrow \catlSp$ assigning
spectrum to a ring. Map $f : R \rightarrow S$ is replaced by
a map of locally ringed spaces
\[ {\rm Spec}\,f = (f^0,f^\sharp):
{\rm Spec}\, S \rightarrow {\rm Spec}\,R.\]
The comorphism $f^\sharp$ is important:
e.g. if $f^o : X \hookrightarrow Y$ is an inclusion, the difference
between a subvariety $(X, {\cal O}_X)$ and, say, its $n$-th infinitesimal
neighborhood in $(Y,{\cal O}_Y)$,
may be expressed by a proper choice of $f^\sharp$.

\ppt After geometrizing rings, one may proceed to geometrize modules.
The basic fact here is the affine {\sc Serre}'s theorem establishing
a correspondence $M \leftrightarrow \tilde{M}$ between the
$R$-modules and quasicoherent sheaves of ${\cal O}_X$-modules,
for $X = {\rm Spec}\,R$. This correspondence is
an equivalence of categories $R-{\rm Mod} \leftrightarrow \catqcoh X$.
Using this equivalence of categories, functors $f^*$ and $f_*$ may be
rephrased as additive functors
\[ f^* : \catqcoh X \to \catqcoh Y,\,\,\,\,\,\,
f_* : \catqcoh Y \to \catqcoh X;\]
and moreover, these can be defined for any morphism
$f = (f^o, f^\sharp)$ between locally ringed spaces.
In this wider context, functor $f^*$  is called
the {\bf inverse image functor} of $f$, given by
$f^* {\cal F} := {\cal O}_Y \otimes_{{\cal O}_X} f^\bullet {\cal F}$
where $f^\bullet {\cal F}$ is usual, sheaf theoretic, pullback of
sheaf ${\cal F}$ via $f^o$. The restriction of scalars functor
then generalizes to the {\bf direct image functor} for sheaves
which is on presheaf level given by
\[ f_* {\cal F} (U) : = {\cal F} (f^{-1}(U)), \]
and which sends ${\cal O}_X$-modules to ${\cal O}_Y$-modules
via ${\cal O}_Y$-action
given by the composition of $f^\sharp \times \id$ and
the ${\cal O}_X$-action. Functors $f^*$ and $f_*$ are naturally defined
between the bigger categories,
${\cal O}_X-{\rm Mod}$ and ${\cal O}_Y-{\rm Mod}$,
where they simply preserve the quasicoherence.

\ppt Functor $f_*$ (in all settings above) is a right adjoint to $f^*$,
hence it is left exact and the inverse image functor $f^*$ is
right exact. This suggests that a pair of adjoint additive
functors between abelian categories may be viewed as (coming from) a
morphism in geometrical sense. Actually this point of view appears
fruitful. On the other hand, surely the choice of a functor in
its equivalence class is not essential; and the emphasis on the
inverse image vs. direct image functor is the matter of choice as well.

Given two abelian categories ${\cal A}$, ${\cal B}$,
(equivalent to small categories) a {\bf morphism}
$f : {\cal B} \rightarrow {\cal A}$ is an equivalence class
of right exact additive functors from ${\cal A}$ to ${\cal B}$.
An {\bf inverse image functor} $f^*: {\cal A}\rightarrow {\cal B}$
of $f$ is simply a representative of $f$,
which is usually assumed to be made.

An additive functor $f^* : {\cal B} \rightarrow {\cal A}$
between \abelian categories is~(\cite{Ros:NcSch})

$\bullet$ {\bf continuous} \index{continuous functor}
if it has a right adjoint, say $f_*$;

$\bullet$ {\bf flat} if it is continuous and exact;
\index{flat functor}

$\bullet$ {\bf almost affine} \index{almost affine}
if it is continuous and its right adjoint
$f_*$ is faithful and exact;

$\bullet$ {\bf affine} \index{affine functor}
if it is almost affine, and its right adjoint
$f_*$ has its own right adjoint, say $f^{!}$, cf.~\luse{functorsMod}.

Morphism $f$ is continuous (flat, almost affine, affine)
if its inverse image functor $f^*$ is. Some authors say that a functor
$F$ is continuous if it has a left adjoint instead, which means
that they view $F$ as a direct image functor $f_*$ of some
continuous morphism $f$. On the other hand, a continuous
morphism is {\bf coflat} if its direct image is exact,
and {\bf biflat} if it is flat and coflat. Usually one equips
the categories with distinguished objects (``structure sheaves'');
then the morphisms respect the distinguished object.

\ppta {\sc Rosenberg}~(\cite{Ros:NcSch}) introduced an abstract
notion of a \wind{quasicompact relative noncommutative scheme}
$({\cal A},{\cal O})$ over a category ${\cal C}$ as an \abelian
category ${\cal A}$ with a distinguished object ${\cal O}$, finite
biflat affine cover by localizations $Q^*_\lambda : {\cal A}\to
{\cal B}_\Lambda$, with a continuous morphism $g$ from ${\cal A}$
to ${\cal C}$ (think of it as a direct image of
a morphism $X\to {\rm Spec}\,\genfd$) such that
each $g_* \circ Q_{\lambda *} : {\cal B}_\lambda\to {\cal C}$ is affine. This
guarantees application of many usual geometric procedure (for
basic cohomological needs one  does not need $f^!$). Such
'schemes' can be related to some spectra and some Grothendieck
topologies on ${\rm Aff}({\rm Ab})$~(\cite{Ros:USpSch}). Quotient
spaces for comodule algebras over Hopf algebras may be sometimes
constructed as nonaffine noncommutative schemes~\cite{Skoda:hg}.

\ppt If $R$ and $S$ are rings and  ${}_S B_{R}$ a $S-R$-bimodule,
then the functor $f^B : M \mapsto {}_S B_{R} \otimes_R \,M$
is a right exact functor from $R-{\rm Mod}$ to $S-{\rm Mod}$.
If $S = {\DDl Z}$ then $B_{R}$ is
called flat right $R$-module if $f^B : M \mapsto B_{R} \otimes_R \,M$
is exact.

{\bf Proposition.} ({\sc Watts}~\cite{Watts}, {\sc Eilenberg}~\cite{Eil:af})
\index{Eilenberg-Watts' theorem}
{\it (i) Let $R$ be a (not necessarily unital) ring,
and $f^*$ a flat endofunctor in the category of nonunital
left $R$-modules.
Then $f^*$ is equivalent to the functor
\[
  M \mapsto f^*(R_1) \otimes_R M,
\]
where $R_1$ is the corresponding unital ring.
In particular, $f^*(R_1)$ is flat as a left $R$-module.

(ii) Let $R$ be a unital ring, and $f^*$ a flat endofunctor
in the category $R-{\rm Mod}$ of unital left $R$-modules.
Then $f^*$ is equivalent to the functor $M\mapsto f^*(R)\otimes_R M$.

(iii) Let $f^* : R-{\rm Mod} \to S-{\rm Mod}$ be a continuous functor.
Then there is a $S-R$-bimodule ${}_S B_{R}$ such that $f^*$ is
equivalent to the functor $M \mapsto B \otimes_R M$.
}

Notice that when applied to localizations word ``flat'' is here
used in the sense that $Q^*$ is flat (weaker), and in some other
works means that $Q := Q_* Q^*$ is flat.
The latter case, for $R-{\rm Mod}$, is the situation of Watts theorem,
and it is known under the name of {\bf perfect localizations}.
Equivalently~(\cite{Popescu}), the canonical forgetful functor from the
localized category $(R-{\rm Mod})/{\rm Ker Q}$
to the modules over the localized ring $(QR)-{\rm Mod}$
is an equivalence of categories.

\ppt A {\bf bicategory} \index{bicategory}
(= weak 2-category) \index{2-category} ${\cal A}$ consists of

(1) a class ${\rm Ob} \,{\cal A}$ of {\it objects} ('1-cells')\index{1-cell};

(2) for each pair of objects $A,B$ a small category ${\cal A}(A,B)$;
The objects of ${\cal A}(A,B)$ are called {\it arrows}
(or 'morphisms' or '1-cells'),
and the morphisms in ${\cal A}(A,B)$ are called {\it 2-cells}\index{2-cell};

(3) for each triple of objects $A,B,C$ a bifunctor ('composition map')

\[ \Phi_{A,B,C} : {\cal A}(A,B) \times {\cal A}(B,C)
\rightarrow {\cal A}(A,C), \]

(4) for each object $A \in {\cal A}$, an arrow
$1_A \in {\cal A}(A,A)$ ('identity arrow');\index{identity arrow}

together with natural equivalences
\[ \begin{array}{c}
a_{ABCD} : \Phi_{ACD} (\Phi_{ABC} \times {\rm Id})
\Rightarrow \Phi_{ABD} ({\rm Id} \times \Phi_{BCD}), \\
\lambda_{ABC} : {\rm Id} \times \Phi_{ABC}\Rightarrow \Phi_{ABC} ,
\,\,\,\,\,\,\,\,\,\,
\rho_{ABC} : \Phi_{ABC} \times {\rm Id}\Rightarrow \Phi_{ABC},
\end{array}\]
satisfying some natural conditions. If $a_{ABCD}, \rho_{ABC},\lambda_{ABC}$
are all identities,
then the bicategory is called
(strict) \index{2-category!strict} 2-category.

We omit further details in the definition, just to
sketch the most important example to us: the bicategory of rings and
bimodules ${\frak Bim}$. Beforehand, notice that from any monoidal
category $\tilde{\cal C} = ({\cal C}, \otimes, 1_{\cal C}, a, l,
r)$ we can tautologically form a bicategory $\Sigma \tilde{\cal C}$ with
one object $0 := {\cal C}$, and with
$\Sigma\tilde{\cal C} (0,0) := \tilde{\cal C}$ and
with the composition $M \circ_{\Sigma\tilde{\cal C}} N := M\otimes N$;
one further defines $a(M,N,P) :=a_{0,0,0}(M,N,P)$, $l_A := \lambda_{0,0}(A)$
and $r_A := \rho_{0,0}(A)$. The definitions may be reversed to form a monoidal
category out of any bicategory with a distinguished object $0$
(forgetting the other objects). Hence
monoidal categories may be viewed as bicategories with a single object;
the notion of being 'strict' in both senses agrees as well.

The objects of ${\frak Bim}$ are unital rings.
For any $R,S \in {\rm Ob}\,{\frak Bim}$ let
${\frak Bim}(R,S)$ be the category of $R$-$S$-bimodules and
$(R,S)$-bilinear mappings between them. The composition
\[ {\frak Bim}(R,S) \times {\frak Bim}(S,T)\to {\frak Bim}(R,T)\]
is given by the tensor product bifunctor
$(M,N) \mapsto {}_R M_{S} \otimes_S {}_S N_T$, and the rest of the data
is obvious. These data indeed define a bicategory.

The invertible 1-cells of ${\frak Bim}$ are called
{\sc Morita} equivalences. It has been observed in
various applications of noncommutative geometry,
for instance in physics, that Morita invariance
is a common feature of natural constructions.

The Eilenberg-Watts' theorem identifies bimodules with flat
functors. As a pair of adjoint functors, they resemble
geometric morphisms among topoi.
{\sc M.~Van den Bergh} (\cite{vdB:blow})
defines a (generalized) bimodule to be any
pair of adjoint functors between Grothendieck categories.
He also considers right exact functors
as so-to-say weak (version of) bimodules.
Some situations, for example the duality for coherent sheaves
involve functors for which the right and left adjoints~(\cite{Hart:res}).
coincide. They are known as
Frobenius functors~(\cite{caen}).~\index{Frobenius functor}
In the spirit of Van den Bergh's terminology,
{\sc Pappacena} calls {\it Frobenius bimodules}\index{Frobenius module}
those adjoint pairs $F\dashv G$ where $F$ is simultaneously left and
right adjoint of $G$. In abstract homotopy theoretic setting,
the existence of two-sided adjoints is studied with appropriate
(Bousfield-type) localization arguments
(\cite{May:dual,Neeman:GdualityBous}).

One of motivations for this~(\cite{Neeman:GdualityBous})
is to extend the {\sc Grothendieck} duality theory for coherent sheaves
on varieties to ${\cal D}$-modules. This may be viewed
as an example of noncommutative geometry.
Namely, the role of the structure sheaf ${\cal O}$ is played by the
sheaf ${\cal D}$ of regular differential operators which is
a sheaf of noncommutative ${\cal O}$-algebras
(cf.~\cite{YekZh:dualncringScheme, PolRoth:DalgPrep}
for the viewpoint of noncommutative geometry at ${\cal D}$-modules).
In triangulated categories, the Serre-Grothendieck duality is
axiomatized as an existence of so-called {\sc Serre} functor
(\cite{BonKap:Serre,BonOrl:ICM,BonVdB:gen,orlov:derusp}),
with applications at the borderline between
the commutative and noncommutative geometry.

It is a remarkable observation~(\cite{Landsman:bicat}),
that the noncommutative geometry via operator algebras,
could be also organized around similar bicategories.
Roughly speaking, operator algebras
($C^*$-algebras; von Neumann algebras respectively) are
0-cells, appropriate bimodules as 1-cells (Hilbert bimodules;
correspondences), and bimodule morphisms as 2-cells;
while the monoidal products of 1-cells are specific tensor
products which became prominent earlier in noncommutative geometry
a la Connes, and related K-theories (Rieffel
interior tensor product of Hilbert  bimodules; Connes fusion product).
Invertible 1-cells are called {\bf Morita equivalences}
in all these cases. There are also analogues concerning regular
bibundles over Lie groupoids, and also analogues in symplectic and
Poisson geometry. The latter may be viewed as a quasiclassical
limit of noncommutative geometry. For Morita equivalence of
Poisson manifolds and similar notion of symplectic dual pairs
see~\cite{Landsman:bookCQM,SilvaWein}.

\minsection{Ore localization in filtered rings}

After prerequisites on filtrations, we prove some general
lemmas on localizations in filtered rings,
mainly easy generalizations of some statements
quoted without proof in~\cite{LRloc1}
and in the manuscript~\cite{LRloc2}\footnote{I thank {\sc Valery
Lunts} for introducing me into this subject and sharing his notes.}.

We focus on {\it 'negative'} filtrations.
The main application in mind is the noncommutative
deformation\index{deformation} of commutative objects.
Such filtrations arise from expanding the algebra operations in power
series in the deformation parameter $q-1$
(\cite{Artin:ncdef, Konts:defqalgvar, Lunts:def}).
There is a more refined technique in algebraic analysis,
{\it algebraic microlocalization}\index{algebraic microlocalization},
see~\cite{vandenEssen:survey, vOyst:assalg}
and the references in~\cite{vandenEssen:survey}. {\it 'Positive'}
filtrations involve different techniques than ours. In the study
of noncommutative projective algebraic geometry there is a
(negative) filtration related to deformation, but also a
positive grading corresponding to the study of projective schemes.
The latter grading may be refined to ${\DDl Z}^{\times n}$-grading
or grading with respect to the weight lattice $P$, as in the study
of the quantum flag varieties~(\cite{LRloc1,LRloc2}).
If the latter, positive, grading complicates the picture,
one restricts attention to homogeneous Ore sets only.\vskip .02in

\ppt Given a (typically abelian) category ${\cal C}$, e.g. the
category of modules over a ring $\genfd$, a ${\DDl Z}$-{\bf
filtration} on an object $M$ in ${\cal C}$ is a nested (unbounded)
sequence of subobjects $F_* M = \{\cdots \subset F_{n-1} M \subset
F_n M \subset F_{n+1} M \cdots \subset M\}$. A ${\DDl Z}$-filtered
object\index{${\DDl Z}$-filtered object} in ${\cal C}$ an object
with a ${\DDl Z}$-filtration on it. All filtrations in this
article are assumed to be {\bf exhaustive}\index{exhaustive
filtration} i.e. the supremum subobject ${\rm sup}_{n\in {\DDl Z}}
F_n M$ exists and equals $M$ (e.g. for modules $M
= \cup_n F_n M$).

Let $M = \cup_{n\in {\DDl Z}} F_n M$ be a filtered $\genfd$-module.
The {\bf degree} $d(m)$ of an element $m \in M$ is the smallest
integer $n$, if it exists (otherwise $-\infty$), such that $m \in F_n M$
and $m \notin F_{n-1}M$.
Conversely, if $d : M\rightarrow \{-\infty\}\cup {\DDl Z}$
is subadditive $d(a+b) \leq d(a) + d(b)$,
and $d(0) = -\infty$, then $d$
is the degree function of a unique exhaustive filtration on $M$,
indeed the one where $a \in F_n M$ iff $d(a) \leq n$.
An (exhaustive) filtration is {\bf separated}\index{separated filtration}
if $\cap_n F_n M = 0$. Then $d(m)$ is finite for all $m\neq 0$.
{\it This will be our assumption from now on.}\vskip .02in

\ppta A ${\DDl Z}$-{\bf filtered $\genfd$-algebra} is a
$\genfd$-algebra $(E,\mu)$ with a filtration $F_* E$
on $E$ as a $\genfd$-module where
the multiplication $\mu$ restricted to
$F_n E \times F_m E$ takes values within $F_{n+m} E$, for all $n,m$.
This obviously generalizes to algebras in
any monoidal category $\tilde{\cal C} := ({\cal C}, \otimes, 1_{\cal C})$
(notice that the resulting notion is different than
if we consider these algebras as
auxiliary {\it objects} in an abstract category ${\cal C}'$ of algebras,
when \thispta{} applies, rather than as {\it algebras} in a
{\it monoidal} category $\tilde{\cal C}$).
For $\genfd = {\DDl Z}$ we talk about ${\DDl Z}$-{\bf filtered rings}.

Given a filtered $\genfd$-algebra $(F_* E,\mu)$ a
${\DDl Z}$-{\bf filtered} $F_* E$-{\bf module}
is a ${\DDl Z}$-filtered $\genfd$-module $F_* M$ such that $F_n E$
acting on $F_m M$ takes values within $F_{n+m} M$
for all $n$ and $m$. In particular,
$F_* E$ is a ${\DDl}$-filtered module over itself.

Given a filtered $\genfd$-algebra $E = \cup_n F_n E$,
an \wind{associated graded algebra} is the ${\DDl Z}$-graded $\genfd$-module
${\rm gr}\,E = \oplus_n ({\rm gr}\,E)_n := \oplus_n F_{n} E / F_{n-1} E$
with the multiplication defined as follows.
The {\bf symbol map} \index{symbol map}
$\bar{} : E \rightarrow {\rm gr}\,E$,
$e \mapsto \bar{e}$, maps $e$ to the class $\bar{e}$
of $e$ in $F_{d(e)} E / F_{d(e)-1} E$.
An element $c \in {\rm gr}\,E$ is in the image of
the symbol map iff $c$ is homogeneous.
For $c = \bar{e}, c'=\bar{e}'$, the formula
$cc' := \bar{e}\bar{e}' := \overline{ee}'$
does not depend on the choice of $e,e'$.
Therefore it defines a binary operation.
The additive extension of this operation is
the associative multiplication on ${\rm gr}\,E$.

It is always true $d(e + e') \leq {\rm max}\{d(e),d(e')\}$,
with equality if $d(e) \neq d(e')$. On the other hand, if $d(e) = d(e')$
then $d(e+ e')$ may be in general anything less or equal $d(e)$,
as $e$ and $e'$ may cancel in several of the top degrees.
Consequently the symbol map is {\it not} additive.
However...

\ppt \ldef{multsymbol} {\it ...if ${\rm gr}\,E$ is a {\it domain}, then
$d(ab) = d(a)d(b)$, hence the symbol map $E\to {\rm gr}\,E$
{\it is} multiplicative.
}

\ppt For any subset $S\subset E$ not containing $0$,
we can always define a filtration on the set $S \times E$,
by formula $d(s,e) := d(e) - d(s)$.
If $S$ is left Ore, the localized ring
$S^{-1}E$ may be constructed as in~\luse{Oreexistence},
as certain quotient of $S \times E$.
Hence we have a filtration on $S^{-1}E$
as a set with degree function
$d(s^{-1}e) = \inf_{s'^{-1}e' = s^{-1}e} \,d(s',e')$.
Recall that $(s,e)\sim (s',e')$ means $\exists \tilde{s}\in S$,
$\exists \tilde{e} \in E$, $\tilde{s}s = \tilde{e}s'$ and
$\tilde{s}e = \tilde{e}e'$.

If the degree function is multiplicative, e.g. $E$ is a domain, then
\[\begin{array}{lcl}
d(s,e) & = & d(e) - d(s)\\
       & = & d(e) - d(\tilde{e}) - (d(e) - d(\tilde{s})) \\
       & = & d(\tilde{s}r) - d(\tilde{s}s) \\
       & = & d(\tilde{e}e') - d(\tilde{e}s') \\
       & = & d(e') - d(s') \\
       & = & d(s',e'),
\end{array}\]
hence taking the infimum in the expression for $d(s^{-1}e)$ is superfluous,
as all the representatives of $s^{-1}e$ give the same result.
Therefore the degree is well-defined by $d(s^{-1}e) := d(e) - d(s)$.

The symbol image of a set $S \subset E$ is denoted by $\overline{S}$.
If $S$ is left Ore in $E$, and ${\rm gr}\,E$ is a domain,
then $\overline{S}$ is clearly left Ore in ${\rm gr}\,E$.

\ppt {\bf Lemma.}\ldef{locmult}
{\it If the symbol map $E\to {\rm gr}\,E$
is multiplicative, then the induced degree
function on $S^{-1}E$ is multiplicative as well.
}

\ppt {\bf Proposition.} (\cite{LRloc1}, II 3.2)\ldef{theta1}{\it
(i) We have a well-defined map $\theta : \bar{S}^{-1} {\rm gr}\,E
\to {\rm gr}\,S^{-1} E$ given by $(\bar{s})^{-1}\bar{e} \mapsto
\overline{s^{-1} e}$.

(ii) This map is an isomorphism of graded rings.
}

{\it Proof. } (i) Let $s_1, s_2 \in S$ with $\bar{s}_1 = \bar{s}_2$.
By  $s_1^{-1} s_2 = 1 + s_1^{-1} (s_2 - s_1)$
we get $s_1^{-1} e = s_2^{-1} e + s_1^{-1} (s_2 - s_1) s_2^{-1} e$.
Using~\luse{locmult} and  $d(s_2-s_1) < d(s_1)$,
we see that for each $0 \neq e \in E$,
$s_1^{-1} (s_2 - s_1) s_2^{-1} e$ is lower degree than $s_2^{-1} e$.
Thus $\overline{s_1^{-1}e} =
\overline{s_2^{-1}e} \in \overline{S^{-1}E}$.
In the same vein, but easier, we see that $\overline{s^{-1}e}$
does not depend on the choice of $e \in \bar{e}$.

Finally, choose different classes $\bar{t}$ and $\bar{f}$ with
$(\bar{t})^{-1} \bar{f} = (\bar{s})^{-1}\bar{e}$. That is $\exists
\bar{s}_* \in \bar{S}, \exists \bar{e}_* \in {\rm gr}\,E$ with
$\bar{s}_*\bar{t} = \bar{e}_*\bar{s} \in \bar{S}$ and
$\bar{s}_*\bar{f} = \bar{e}_*\bar{e}$. Then $\overline{s_* t} =
\overline{e_* s}$ and $\overline{s_* f} = \overline{e_* e}$ for
any choice of representatives $e_*, s_*$ of $\bar{e}_*,
\bar{s}_*$. Hence there are $r_1, r_2$ of lower degrees than $e_*
s, e_*e$ respectively, such that $s_* t = e_* s + r_1$ and $s_* f = e_* e +
r_2$. Then $t^{-1} f = (s_* t)^{-1} s_* f = (e_* s + r_1)^{-1}
(e_* e + r_2)$ which by the above equals $(e_* s)^{-1} (e_* e) =
s^{-1} e$ up to elements of lower order, {\it provided that} $e_*
s \in S$. As $e_* s + r_1 \in S$, this is always true if $S$ is
{\it saturated}, see below. However, the conclusion follows
without that assumption. Indeed, by the left Ore condition, choose
$s^\sharp, e^\sharp$ with $S \ni s^\sharp (e_* s + r_1) = e^\sharp
s$. Then $e^\sharp e = e^\sharp s s^{-1} e = s^\sharp e_* e +
s^\sharp r_1 s^{-1} e$, where, by the multiplicativity, $s^\sharp r_1
s^{-1} e$ is of lower order. Consequently, $t^{-1} f = [s^\sharp
(e_* s + r_1)]^{-1} s^\sharp (e_* e + r_2) = (e^\sharp s)^{-1}
e^\sharp e + \mbox{lower order} = s^{-1} e  + \mbox{lower order}$,
as required.

(ii) Since both the degree of $(\bar{s})^{-1}\bar{e}$ and
of $\overline{s^{-1} e}$ are $d(e)-d(s)$, this map
respects the grading. The obvious candidate $\overline{s^{-1} e}\mapsto
(\bar{s})^{-1}\bar{e}$ for the inverse is well-defined by more
straightforward reasons than the map $\theta$. Namely, if
$\overline{t^{-1}f} = \overline{s^{-1}e}$ then $\exists h$ of
lower order with $t^{-1} f = s^{-1}e + h = s^{-1}(e+sh)$. As
$\overline{e} = \overline{e+sh}$ it is enough to check the case $h
= 0$. For some $s_* \in S$, $r_* \in r$ we have $s_* t = r_* s \in
S$ and $s_* f = r_* e$. Then $\bar{s}_* \bar{t} = \bar{r}_*
\bar{s} \in \bar{S}$ and $\bar{s}_* \bar{f} = \bar{r}_* \bar{e}$,
hence $\bar{t}^{-1} \bar{f} = \bar{s}^{-1} \bar{e}$, as required.

\ppt Let  $N = \cup_{k \in {\DDl Z}} F_k N$, be a right and
$M =  \cup_{k \in {\DDl Z}} F_k M$ a left filtered $E$-module,
then $N\otimes_\genfd M$ is filtered with respect to
the unique degree function additively extending
formulas $d(n\otimes_k m) = d(n) + d(m)$.
The canonical quotient map $p_E : N \otimes_\genfd M\rightarrow N\otimes_E M$
induces the filtration $F_k (N\otimes_E M):= p_E(F_k (N \otimes_\genfd M))$.
If $N$ is a filtered $E'-E$-bimodule, one obtains
a filtration of $N\otimes_E M$ as a left $E'$-module.
In particular, {\it given $E = \cup_k F_k E$, where ${\rm gr}\,E$ is a domain,
and given a filtered left $E$-module $M = \cup_k F_k M$,
any Ore localization $S^{-1}M = S^{-1}E \otimes_E M$ is a filtered
left $S^{-1}E$-module with the degree function

$d(s^{-1}m) = d(s^{-1}\otimes_E m) = d(s^{-1}1_E) + d(m) = d(m)-d(s)$}.

\ppt {\bf Lemma.} \ldef{multact}
{\it If the symbol map $E\to {\rm gr}\,E$
is multiplicative, and $M$ a filtered left $E$-module,
then the degree functions are compatible with action
in the sense that $d_M(e.m) = d_E(e) d_M(m)$.
Furthermore, for any left Ore set $S\subset E$,}
$$\begin{array}{lcl}
d_{S^{-1}M}(s^{-1}e.t^{-1}m) &=& d_{S^{-1}E}(s^{-1}e)d_{S^{-1}M}(t^{-1}m) \\
&=& d_E(e) + d_M(m) - d_E(s) - d_E(t).\end{array}$$

\ppt {\bf Proposition.} \ldef{theta2}{\it (i) For a filtered ring $E$,
for which ${\rm gr}\,E$ is a domain, and
any filtered $E$-module $M$, we have a well-defined map
$\theta_M : \bar{S}^{-1} {\rm gr}\,M
\to {\rm gr}\,S^{-1} M$ given by $(\bar{s})^{-1}\bar{m} \mapsto
\overline{s^{-1} m}$.

(ii) $\theta_M$ is an isomorphism of graded
${\rm gr}\,S^{-1}E = \overline{S}^{-1}{\rm gr}\,E$-modules.
}

The proof is by the same techniques as~\luse{theta1}. The
compatibility with the action~\luse{multact} replaces the
multiplicativity, and the formula~(\ref{eq:kernelNu}) for the
equivalence relation $\sim$ on $S\times M$ (with $\left(S\times
M/\sim \right)\,\cong S^{-1}M$) replaces the equivalence relation
$\sim$ from~\luse{relFrac} on $S\times R$ in that proof.

\ppt \ldef{filteredOre}{\bf Ore conditions recursively.} {\it
(i) Let $S$ be a multiplicative set in a ring $E$ with
an exhaustive filtration
\[ F_* E = \{ \ldots \subset F_{-r}E \subset F_{-r+1}E \subset \ldots
\subset F_{-1}E \subset F_0 E \subset F_1 E \subset
\cdots \subset  E\}.\]
Let $S$ satisfy the {\bf bounded below filtered-relative left Ore condition}
\index{filtered-relative Ore condition} in $F_* E$:

\indent{
$\exists r$, $\infty > r \geq -n$,
$\forall s \in S$, $\forall k$, $-r \leq k \leq n$,
$\forall e \in F_k E$, $\exists s' \in S$, $\exists e'\in E$
such that $s'e - e's \in F_{k-1} E$ if $k > -r$,
and $s'e-e's = 0$ if $k = -r$.
}

Then $S$ satisfies the left Ore condition for $S$ in $E$.

(ii) Assume that $S$ is {\bf bounded filtered left reversible}
in $F_* E$:

$\exists r < \infty$, $\forall e_k \in F_k E$, if $\exists s \in S$ with
$e_k s \in F_{k-1} E$ then $\exists s' \in S$ such that
$s' e_k \in F_{k-1} E$ if $k > -r$, and
$s' e_k = 0$ if $k = -r$.

Then $S$ is left reversible in $E$.
} 

{\it Proof.}  (i)
Let $s \in S$ and $e = e_n \in F_n E$. By induction, we can complete
sequences $e_{n}, \ldots, e_{-r}$, $e'_{n},\ldots, e'_{-r}$
(here $e_k, e'_k \in F_k E$) and
$s'_n,\ldots, s'_{-r} \in S$,
with $e'_k s =  s'_k e_k - e_{k-1}$ for all $k$ with $e_{-r-1} := 0$.
By descending induction on $k$,
\[ (e'_k + s'_{k} e'_{k+1} +
 \ldots + s'_{k} \cdots s'_{n-1} e'_n )\, s\,
=\, s'_{k} \cdots s'_n e_n - e_{k-1},  \]
for each $k > -r$, and finally,
\[ (e'_{-r} + s'_{-r} e'_{-r+1} + \ldots +
s'_{-r} s'_{-r + 1} \cdots s'_{n-1} e'_n )\, s
\,= \, s'_{-r} \cdots s'_n e_n. \]

(ii) Suppose $e \in F_k E = E$ and $es = 0$ for some
$s \in S$. It is sufficient to inductively choose a descending sequence
and $s'_{k+1} = 1, s'_k, s'_{k-1}, \ldots, s'_{-r} \in S$, with requirements
$s'_j e \in F_{j-1}E$ for all $j> -r$ and $s'_{-r} e = 0$.
Suppose we have chosen $s_k,\ldots, s_{j+1}$. Then
$(s'_{j+1} e) s = s'_{j+1} (es) = 0$ with $s'_{j+1} e \in F_j E$,
hence by the assumption there exist some $\tilde{s} \in S$ such
that $\tilde{s} s'_{j+1} e \in F_{j-1}E$. Set therefore
$s'_j := \tilde{s} s'_{j+1} \in S$.

\ppta
Let $F_* E$ be an exhaustive filtration of $E$ with $F_{-r} E = 0$
for some finite $r$, and $S \subset E$ be
a multiplicative set.
If its image $\bar{S}$ under the symbol map
satisfies the left Ore condition in ${\rm gr}\,E$,
then the conditions in~\thispt{} hold.
Hence $S$ satisfies the left Ore condition in $E$ as well.

\pptb \ldef{EnDefined}
 Let $t \in E$ be a regular element
\index{regular element}($tE = Et$) in a ring $E$.
Then for each $n > 0$ the ideal $t^n E$ is 2-sided, hence
$E_n := E/(t^n E)$ is a quotient ring in which the element $t$ is
nilpotent of order less or equal $n$.
Rule $F_{-k} E_n = (t^k E)/(t^n E)\subset E/(t^n E) \equiv E_n$
defines a bounded 'negative' filtration
\[ F_* E_n = \{ 0 = F_{-n} E_n \subset \ldots \subset F_{-1}E_n
\subset F_0 E_n = E_n\}\]
in which (the image of)  $t$ is of degree $-1$.
If ${\rm gr}\,E$ is a domain then {\it both} the symbol map
$E\rightarrow {\rm gr}\,E$, and its truncation
$E_n \rightarrow {\rm gr}\,E_n$ are multiplicative.

\ppt \ldef{deponSbarOnly}{\bf Theorem.} {\it
Let $S$ be a left Ore set in some ring $E_n = F_0 E_n$
with a bounded negative filtration $F_\bullet E_n$.
Suppose $S'\subset E_n$ is a multiplicative set such that
$s' \in S'\cap F_{j} E$ iff $\exists b \in F_{j-1} E$ such that
$s' = s - b$. Then $S'$ is left Ore as well
and $S^{-1} E_n = (S')^{-1} E_n$ as graded rings.
}

{\it Proof.} Since the left Ore localization is a universal object
in the category ${\cal C}_l(E_n,S)$ (cf. Chap.4) it is enough to see that a map
of rings $j : E_n \rightarrow Y$ is in it iff it is in
${\cal C}_l(E_n, S')$. If $j(s)$ is invertible in $Y$,
let $c = j(s)^{-1}j(b)$. Mapping $j$ induces a (non-separated in
general) filtration on $Y$ such that $j$ is a map of filtered
rings, by taking the degree to be the infimum of expressions $d(e)-d(t)$
for elements which can be represented in the form $j(t)^{-1} j(e)$
and $-\infty$ otherwise.
With our numerical constraints on the degree, for nonvanishing $e\in E_n$
this difference can not be less than $-n$.
As $d(c) < 0$ we obtain $d(c^{n}) < n-1$,
hence $c^n = 0$. Thus we can invert
$j(s)^{-1}(j(s - b)) = 1 - c$ to obtain the geometrical progression
$ \sum_{j = 0}^{n-1} c^j$.
Then $\sum_{j = 0}^{n-1} c^j j(s)^{-1} j(s-b) = 1$
hence $j(s-b)$ is invertible in $Y$.

It remains to check that $se = 0$
for some $s \in S$ iff $\exists s'' \in S'$ with $s'' e = 0$.
We proceed by induction on the degree $j$ of $e$ starting at $-n$
where $s'e - se \in F_{-n-1} = 0$ for $s' = s-b$ with
the degree of $b$ smaller than of $s'$ hence negative.
For any $j$,
$s' e = (s-b)e = -be$ has the degree at most $j-1$.
On the other hand, by the left Ore condition,
we can find $s_* \in S$, and $b_* \in E$ with
$s_* be = b_* s e = 0$, hence $s_* (s' e) = 0$.
Set $e' := s'e$. Since $s_* e' = 0$ with $d(e')< d(e)$,
by the inductive assumption there exists
$s'_* \in S'$ with $s'_* s' e = s'_* e' = 0$.
Set $s'' := s_*' s'$.

The induced grading on the two localized rings is the same after the
identification,
because the symbol maps evaluate to the same element on
$s$ and $s' = s+b$ (or, alternatively, after the identification,
the gradings on the localization are induced by the {\it same} ring map).

{\bf Definition.} {\it A multiplicative subset $S \subset E_n$ is

$\bullet$ {\bf admissible} \index{admissible subset}
if $\forall s \in S$, $0 \neq \bar{s} \in {\rm gr}\,E_n$;

$\bullet$ {\bf saturated} \index{saturated multiplicative subset}
if $S = \{ s \in E_n \,|\, \bar{s} \in \bar{S} \}$.
}

\ppt {\bf Corollary.} {\it Let $E, t, E_n, F_\bullet E_n$
be as in~\luse{EnDefined}.
Suppose ${\rm gr}\,E$ is a commutative domain.
Let $S$ be a {\em multiplicative} subset in $E_n$. Then

a) $S$ is left and right Ore.

b) $S^{-1} E_n \neq 0$ iff $S$ is admissible.

c) $S^{-1} E_n$ depends only on $\bar{S}\subset \bar{E}$.

d) $S^{-1} E_n$ is filtered by powers of $t$ and
$(S^{-1}E_n)/\langle t \rangle \cong \bar{S}^{-1} E_n$.

e) Any two saturated Ore sets $S$, $T$ are compatible, i.e.
$S^{-1} T^{-1} E_n \cong T^{-1} S^{-1} E_n$
and $ST = \{ st \,|\,s\in S, t\in T\}$ is also saturated.

f) Let $S$ be admissible.
Then ${\rm gr}(S^{-1}E_n) \cong \bar{S}^{-1} {\rm gr} E_n$.
In particular, $\overline{S^{-1}E_n} \cong \overline{S}^{-1}\overline{E}_n$.
}

{\it Sketch of the proof.} a) follows as a simple case of~\luse{filteredOre};
b) is trivial; c) follows by~\luse{deponSbarOnly}; d) is evident;
f) follows from \luse{theta2} after truncating both
sides from $E$ to the quotient filtered ring $E_n$ (it is not
a special case of \luse{theta2}, though, as $E_n$ is {\it not}
a domain); e) Because $T$ is saturated, $(T^{-1}E)_n \cong T^{-1}E_n$.
In $T^{-1}E_n$, set $S$ is still multiplicative,
hence by a) applied to $(T^{-1}E)_n$
it is left Ore. This is equivalent to compatibility
(cf. Sec.~10).

\minsection{Differential Ore condition}

An extensive literature is dedicated to differential structures
of various kind associated to objects of noncommutative geometry:
derivations and rings of regular differential operators on NC rings,
$1^{\rm st}$ and higher order differential calculi,
with and without (bi)covariance conditions,
NC connections and de Rham complexes etc.

\ppt \ldef{extendingDiffCalc} {\it
Let $\partial : R \rightarrow R$
be an $R$-valued derivation on $R$
and $S$ a left Ore set in $R$.
Then the formula
\begin{equation}\label{eq:derOnLoc}
 \bar\partial(s^{-1}r) = s^{-1}\partial(r) - s^{-1}\partial(s)s^{-1}r,
\,\,\,\,\,s\in S, r \in R,
\end{equation}
defines a derivation $\bar\partial : S^{-1}R\rightarrow S^{-1}R$.

The same conclusion if we
started with $\partial : R \rightarrow S^{-1}R$ instead.
}

{\it Proof.} 1. {\sc $\bar\partial$ is well defined.}

Suppose $s^{-1}r = t^{-1}r'$ for some $r,r' \in R$, $s,t \in S$. Then
\[ \exists \tilde{s}\in S,
\,\exists\tilde{r}\in R,\,\,\,\,\, \tilde{s}t = \tilde{r}s,
\,\,\,\,\,\,\,\,\tilde{s}r' = \tilde{r}r. \]
\[ s^{-1} = t^{-1} \tilde{s}^{-1} \tilde{r} \]
\[\begin{array}{lcl} t^{-1}\partial(r')
& =& t^{-1}\tilde{s}^{-1}\tilde{s}\,\partial(r')
\\&=& t^{-1} \tilde{s}^{-1}
[ \partial(\tilde{s}r') - \partial(\tilde{s})r' ] \\
&=& t^{-1} \tilde{s}^{-1} [ \partial(\tilde{r}r) - \partial(\tilde{s})r' ]
\\ t^{-1}\partial(t) & = & t^{-1}\tilde{s}^{-1}\tilde{s}\,\partial(t)
\\ &=& t^{-1}\tilde{s}^{-1}[ \partial(\tilde{s}t) - \partial(\tilde{s})t ]
\\ &=&  t^{-1}\tilde{s}^{-1}[ \partial(\tilde{r}s) - \partial(\tilde{s})t ]
\end{array}\]

\noindent
\[\begin{array}{lcl}
\bar\partial(t^{-1}r') &=& t^{-1}\partial(r') - t^{-1}\partial(t)t^{-1}r'\\
&=&  t^{-1} \tilde{s}^{-1} [ \partial(\tilde{r}r) - \partial(\tilde{s})r' ]
- t^{-1}\tilde{s}^{-1}[ \partial(\tilde{r}s) -
\partial(\tilde{s})t ]t^{-1}r'\\
&=&  t^{-1} \tilde{s}^{-1}\partial(\tilde{r}r) -
t^{-1}\tilde{s}^{-1}\partial(\tilde{r}s)t^{-1}r' \\
&=& t^{-1} \tilde{s}^{-1}\partial(\tilde{r}r) -
t^{-1}\tilde{s}^{-1}\partial(\tilde{r}s)s^{-1}r \\
&=& t^{-1} \tilde{s}^{-1}\partial(\tilde{r})r
+ t^{-1} \tilde{s}^{-1}\tilde{r}\partial(r) \\
& & \,\,\,\,\,\,\,\,\,\,\,\,\,\,\,
- t^{-1}\tilde{s}^{-1}\partial(\tilde{r})s s^{-1}r
- t^{-1}\tilde{s}^{-1}\tilde{r}\partial(s) s^{-1}r
\\ &=& s^{-1}\partial(r) - s^{-1} \partial(s)s^{-1}r
\\ &=&  \bar\partial(s^{-1}r)
\end{array}\]
 2. {\sc $\bar\partial$ is a derivation.} We have to prove that
for all $s,t \in S$ and $r,r' \in R$
\begin{equation}\label{eq:derprodfrac1}
\bar\partial(s^{-1}rt^{-1}r') =
\bar\partial(s^{-1}r)\,t^{-1}r' + s^{-1}r\bar\partial(t^{-1}r').
\end{equation}
The argument of $\bar\partial$ on the left hand side
has to be first changed into a left fraction form
before we can apply the definition of $\bar\partial$.
By the left Ore condition, we can find $r_*\in R$, $s_*\in S$
such that $r_* t = s_* r$ i.e. $rt^{-1} = s_*^{-1}r_*$.

We first prove identity~(\ref{eq:derprodfrac1}) in the case
$s= r' = 1$ i.e.
\begin{equation}\label{eq:derprodfrac2}
\bar\partial(rt^{-1}) =
\partial(r)\,t^{-1} + r\bar\partial(t^{-1}).
\end{equation}
The left-hand side of~(\ref{eq:derprodfrac2}) is
\[ \begin{array}{lcl}
\bar\partial(rt^{-1}) &=& \bar\partial(s_*^{-1}r_*)
\\ &=& s_*^{-1}\partial(r_*) - s_*^{-1} \partial(s_*) s_*^{-1} r_*
\\ &=&  s_*^{-1}\partial(r_*) + \bar\partial(s_*^{-1}) r_*.
\end{array}\]
The right-hand side of~(\ref{eq:derprodfrac2}) is
\[\begin{array}{lcl}
\partial(r) t^{-1} - rt^{-1}\partial(t)t^{-1}
&=& \partial(r)t^{-1} - s_*^{-1} r_* \partial(t) t^{-1}
\\&=& \partial(r) t^{-1} - s_*^{-1} \partial(r_* t) t^{-1}
- s_*^{-1} \partial (r_*) t t^{-1}
\\&=& \partial(r) t^{-1} - \bar\partial(s_*^{-1} r_* t) t^{-1}
+ \bar\partial(s_*^{-1})r_* - s_*^{-1} \partial (r_*)
\\&=& \partial(r) t^{-1} - \partial(r) t^{-1} -
\bar\partial(s_*^{-1})r_* - s_*^{-1} \partial (r_*)
\\ &=& \bar\partial(s_*^{-1})r_* - s_*^{-1} \partial (r_*).
\end{array}\]
Hence~(\ref{eq:derprodfrac2}) follows.
Using~(\ref{eq:derprodfrac2}),
we prove~(\ref{eq:derprodfrac1}) directly:
\[\begin{array}{l}
\begin{array}{lcl}
\bar\partial(s^{-1}r t^{-1} r') &=&  \bar\partial((s_* s)^{-1} r_* r')
\\&=& (s_*s)^{-1}\partial(r_* r')
- (s_* s)^{-1} \partial(s_* s) (s_* s)^{-1} r_* r'
\\ &=& s^{-1} s_*^{-1} \partial(r_*) r' + s^{-1} s_*^{-1} r_* \partial(r')
\\ && \,\,\,\,\,\,\,\,\,\,\,\,\,\,
- s^{-1} s_*^{-1} \partial(s_*) s_*^{-1} r_* r'
- s^{-1} \partial(s) s^{-1} s_*^{-1} r_* r'\end{array}
\\
\begin{array}{lcl}
&=&
s^{-1} s_*^{-1}\partial(r_*)r' + s^{-1} t^{-1} r \partial(r') +
s^{-1} \bar\partial(s_*^{-1})r_* r'
+ \bar\partial(s^{-1})s_*^{-1}r_* r'
\\ &=&
s^{-1} \bar\partial(s_*^{-1} r_*) r' - s^{-1} \bar\partial(s_*^{-1})r_* r
+ s^{-1} rt^{-1} \partial(r') + \bar\partial(s^{-1})rt^{-1}r'
\\ &=& s^{-1}\bar\partial(rt^{-1}) - s^{-1} \bar\partial(s_*^{-1})r_* r
+ s^{-1} rt^{-1} \partial(r') + \bar\partial(s^{-1})rt^{-1}r'
\\ &\stackrel{(\ref{eq:derprodfrac2})}{=}&
 s^{-1} \partial(r) t^{-1} r' + s^{-1} r \bar\partial(t^{-1})r'
+ s^{-1}r t^{-1} \partial(r') + \bar\partial(s^{-1}) r t^{-1} r'
\\ &=& \bar\partial(s^{-1} r) t^{-1} r' +
s^{-1} r \bar\partial(t^{-1}r').
\end{array}\end{array}\]
Standard textbooks have incomplete proofs
of \luse{extendingDiffCalc}, e.g.~\cite{Dixmier:env, Rowen}.

\vskip .02in

\ppt {\bf Definition.} {\it
A \wind{Poisson bracket} on a unital associative $\genfd$-algebra is an
antisymmetric bilinear operation $\{,\} : A \otimes A\rightarrow A$
satisfying the Jacobi identity
$\{f,\{g,h\}\}+\{h,\{f,g\}\}+\{g,\{h,f\}\}= 0$ for all $f,g,h \in A$
and such that for each $f$, $\genfd$-linear map
$X_f : g \mapsto \{f,g\}$ is a $\genfd$-derivation of $A$.
A \wind{Poisson algebra} is a {\em commutative} algebra with a Poisson
bracket.
}

{\bf Proposition.} {\it
Let $A$ be a $\genfd$-algebra with a Poisson bracket $\{,\}$,
and $S\subset A\backslash \{0\}$ a {\em central} multiplicative set.
Then

(i) $S^{-1}A$ posses a bilinear bracket
$\{,\} = \{,\}^S$ such that the localization map
$\iota_S : A \rightarrow S^{-1}A$ intertwines the brackets:
$\{,\}^S\circ(\iota_S\otimes_\genfd \iota_S)=\iota_S\circ\{,\}$.

(ii) If either $\{s,t\} \in {\rm Ker}\,\iota_S$
for all $s,t \in S$, or if $A$ is commutative,
then there is a unique such bracket $\{,\}^S$ which is,
in addition, skew-symmetric.

(iii) If $A$ is commutative then this unique  $\{,\}^S$
is a Poisson bracket.
}

{\it Proof.} (i) Each $X_b$ by~\luse{extendingDiffCalc} induces a
unique derivation $X^S_b =\bar\partial$ on
$S^{-1}A$ by~(\ref{eq:derOnLoc}) for $\partial = X_b$.
Map $b \mapsto X^S_b$ is $\genfd$-linear
by uniqueness as $X^S_b + X^S_{c}$ is a derivation extending
$X_{b + c}$ as well.
For each $s^{-1}a \in S^{-1}A$ define $\genfd$-linear map
$Y_{s^{-1}a} : A \to S^{-1}A$ by
\[ Y_{s^{-1}a} : b \mapsto -X^S_b (s^{-1}a)
= - s^{-1} \{b, a\} + s^{-1}\{b,s\}s^{-1} a.
\]
{\em Because $s$ is central}, $Y_{s^{-1}a}$ is a
$\genfd$-linear derivation. Namely,
\[\begin{array}{lcl}
Y_{s^{-1}a}(b c)
&=& - s^{-1} \{bc, a\}  + s^{-1}\{bc,s\}s^{-1} \\
&=& - s^{-1} \{b, a\} c - s^{-1} b\{c, a\} + \\ &&
  \,\,\, + \, s^{-1}b\{c,s\}s^{-1}a + s^{-1}\{b,s\}c s^{-1}a,
\end{array}\]
and, on the other hand,
\[\begin{array}{lcl}
Y_{s^{-1}a}(b) c + b Y_{s^{-1}a}(c)
&=& - s^{-1} \{b, a\} c + s^{-1}\{b,s\}s^{-1} a c -\\
&& \,\,\,-\, b s^{-1} \{c, a\}  + b s^{-1}\{c,s\}s^{-1}a.
\end{array}\]
Hence $Y_{s^{-1}a}$ extends to
a derivation $Y^S_{s^{-1}a}$ on $S^{-1}A\rightarrow S^{-1}A$ by formula
~(\ref{eq:derOnLoc}) as well. Define
$\{s^{-1}a, t^{-1}b\}:=  Y^S_{s^{-1}a}(t^{-1}b)$.

(ii) To show the skew-symmetry, we calculate,
\[\begin{array}{lcl}
Y^S_{s^{-1}a}(t^{-1}b)
&=& t^{-1}Y_{s^{-1}a}(b) - t^{-1}Y_{s^{-1}a} t^{-1}b \\
&=&  -t^{-1}X^S_b(s^{-1}a) + t^{-1}X^S_t(s^{-1}a) t^{-1}b\\
&=& -t^{-1} s^{-1} X_b (a) + t^{-1}s^{-1} X_b(s) s^{-1}a \\
&&\,\,\,+ t^{-1}s^{-1} X_t (a) t^{-1}b - t^{-1}s^{-1}X_t(s) s^{-1}a t^{-1}b.
\end{array}\]
\[\begin{array}{lcl}
Y^S_{t^{-1}b}(s^{-1}a)
&=& s^{-1}Y_{t^{-1}b}(a) - s^{-1}Y_{t^{-1}b}(s)s^{-1}a\\
&=&  -s^{-1}X^S_a(t^{-1}b) + s^{-1}X^S_s(t^{-1}b) s^{-1}a\\
&=& -s^{-1} t^{-1} X_a (b) + s^{-1}t^{-1} X_a(t) t^{-1}b \\
&&\,\,\,+ s^{-1}t^{-1} X_s (b) s^{-1}a - s^{-1}t^{-1}X_s(t) t^{-1}b s^{-1}a.
\end{array}\]
Using $X_b(a) = - X_a(b)$ etc. and centrality of $s, t$ we see that
the first 3 terms in $Y^S_{s^{-1}a}(t^{-1}b)$ match with
negative sign the
first 3 terms (in order 1,3,2) in expression for
$Y^S_{t^{-1}b}(s^{-1}a)$. If $a$ and $b$ mutually commute,
the 4th term agrees the same way, and if they don't but
$X_s(t) = 0$ in the localization $S^{-1}A$,
then they are simply $0$, implying skew-symmetry
$\{s^{-1}a,t^{-1}b\} + \{t^{-1}b,s^{-1}a\} = 0$.

Uniqueness: $Z_{s^{-1}a}(t^{-1}b) := \{s^{-1}a, t^{-1}b\}$
defines a derivation $Z_{s^{-1}a}$ on $S^{-1}A$, which restricts
to a derivation $Z_{s^{-1}a}| : A \rightarrow S^{-1}A$. On the other hand,
$s^{-1}a \mapsto Z_{s^{-1}a}(b)$ is $-X^S_b$ by its definition.
Hence the value of $Z_{s^{-1}a}|$ is determined at every $b$, and by
\luse{extendingDiffCalc} this fixes $Z_{s^{-1}a}$.

(iii) We'll prove that if the Jacobi rule holds for given $(a,b,c)$
and $(s,b,c)$ in $S^{-1}A^{\times 3}$, then it follows for
$(s^{-1}a,b,c)$ provided $s$ is invertible. By symmetry of the
Jacobi rule and by renaming $s^{-1}a \mapsto a$ we infer that it
follows for $(s^{-1}a,t^{-1}b,c)$, as well, and finally for the
general case by one more application of this reasoning. Thus we
only need to show that Jacobi$(a,b,c)$ implies
Jacobi$(s^{-1}a,b,c)$. For commutative $S^{-1}A$ this is a
straightforward calculation, using the Jacobi identity, lemma above
and skew-symmetry. We name the summands:
\[ \begin{array}{lcl}
\{s^{-1}a, \{b,c\}\} &=& s^{-1}\{a,\{b,c\}\} - s^{-2}\{s,\{b,c\}\}a
=: (A1) + (A2)\\
\{b , \{c, s^{-1}a\}\} &=& s^{-1} \{b,\{c,a\}\}
- s^{-2}\{b,s\}\{c,a\} - s^{-2}\{b,\{c,s\}\}a
\\ &&\,\,\,-\, s^{-2}\{c,s\}\{b,s\}a -s^{-1}\{c,s\}\{b,a\}
\\
&=:& (B1) + (B2) + (B3) +(B4) +(B5)
\\
\{c , \{s^{-1}a,b\}\} &=&
s^{-1}\{ c, \{a,b\}\} - s^{-2} \{c,s\}\{a,b\} -
s^{-2} \{c,\{s,b\}\}a \\
&& \,\,\,+\, s^{-3}\{s,b\}\{c,s\}a - s^{-2}\{s,b\}\{c,a\}
\\ &=:& (C1) + (C2) +(C3) +(C4) + (C5).
\end{array}\]
Then $(A1) + (B1) + (C1) =0$ and $(A2) + (B3) + (C3) =0$
 by Jacobi for $(a,b,c)$, and $(b,c,s)$ respectively.
By skew-symmetry $(B2) + (C5) = 0$,
$(B5) + (C2) = 0$ and $(B4) + (C4) = 0$ which finishes the proof.

\vskip .017in
{\footnotesize
This fact for $A$ (super)commutative
is used for example in the theory of integrable systems,
sometimes in connection to 'quantization' which is a
rich source of examples in noncommutative geometry.
}

\vskip .02in
\ppt Let $(R,\cdot,+)$ be a ring ($\genfd$-algebra), not necessarily unital.
A {\it first order differential calculus} (FODC)
\index{first order differential calculus}\index{FODC}
is a $R-R$-bimodule $\Omega^1(R)$ together with an additive ($\genfd$-linear)
map $d: R \rightarrow \Omega^1(R)$ satisfying
Leibnitz identity \index{Leibnitz identity}
\[ d (ab) = d(a)b + ad(b), \,\,\,\,\,\,a,b \in R, \]
and such that $\Omega^1(R)$ is generated by differentials $dr$, $r
\in R$ as a {\it left module}. Define a category $\catfodc$:
objects are pairs of a ring $R$ and a FODC $(\Omega^1(R),d)$ on
$R$. A morphisms is a pair $(f,e) : (R,\Omega^1(R),d) \rightarrow
(R',\Omega^1(R'), d')$ of a ring map $f : R \rightarrow R'$ and a
map $e : \Omega^1(R) \rightarrow \Omega^1(R')$ of $R-R$-bimodules
such that $e \circ d = d'\circ f$. Fixing $R$ and allowing only
morphisms of the form $({\rm Id}_R, d)$ we obtain a (non-full)
subcategory $\catfodc_R$ of $\catfodc$. If $R$ is unital, then
$({\rm Ker}(R \otimes_\genfd R \stackrel{\cdot}\rightarrow R,d)$
where $da = 1 \otimes a - a \otimes 1$, and the $R$-bimodule
structure is ${}_R R\otimes_\genfd R_{R}$, is an initial object of
that category.

Two objects $c_R = (R,\Omega^1(R),d)$, $c_{R'} = (R',\Omega^1(R'), d')$
in $\catfodc$
are {\bf compatible along} $f : R \rightarrow R'$
\index{compatible calculi}
if there is an $e$ such that $(f,e)\in \catfodc(c_R,c_{R'})$.

{\bf Differential calculi restrict}: Given $c_{R'}\in \catfodc_{R'}$
and $f$ as above,
define $f^1_\Omega \Omega^1(R')$ to be the smallest additive subgroup
of  $\Omega^1(R')$ containing all the elements of the form
$f(a)\partial'(f(b))$, $a, b \in R$. It appears to be an $R-R$-bimodule.
Define $f^\sharp(c_{R'}) := (R, f^1_\Omega \Omega^1(R'), d' \circ f)$.
Then $f^\sharp(c_{R'}) \in \catfodc_R$ because
$ \partial(b).c = \partial'(f(b))f(c) =
\partial'(f(bc)) - f(b)\partial'(f(c)) = \partial(bc) -
b.\partial(c) \in f^\sharp \Omega^1(R').$
where $\partial = \partial' \circ f : R \rightarrow f^\sharp
\Omega^1(R')$ is the restricted differential.
Note the decomposition
of $(f,e): c_R \rightarrow c_{R'}$
into $(f,e) : c_R \rightarrow f^\sharp c_{R'}$ and
$({\rm id}_{R'}, {\rm incl})
: f^\sharp c_{R'} \rightarrow c_{R'}$, where
${\rm incl}:f^1_\Omega \Omega^1(R')\rightarrow \Omega^1(R')$ is
the inclusion of $R'$-bimodules.

Unlike restricting, there is {\it no general recipe for extending}
\index{extending calculi} the calculus
along ring maps $f : R \rightarrow R'$,
except for the special case when $R' = S^{-1}R$ and
$\Omega^1{R} = {}_R R_{R}$, treated in~\luse{extendingDiffCalc}.
That case is of central importance in study of the regular
differential operators and D-modules over noncommutative spaces
~(\cite{LuntsRosMP, LuntsRos1, LRloc1}).
We'll just mention a slight generalization.

\ppt {\bf Theorem.} {\it Let $S \subset R$ be a left Ore set in a ring $R$,
and suppose $\{ x \in \Omega^1(R)\,|\,\exists t \in S, xt = 0\} = 0$.

The following are then equivalent:

(i) The $S^{-1}R$-$R$-bimodule structure on
$S^{-1}\Omega^1(R) \equiv S^{-1}R \otimes_R\Omega^1(R)$
extends to an (actually unique) $S^{-1}R$-bimodule structure
which may carry a differential $d_S : S^{-1}R \rightarrow S^{-1} \Omega^1(R)$
such that the pair of localization maps
$(\iota_S, \iota_{S,\Omega^1(R)})$ is a morphism in $\catfodc$
(i.e. 'the calculi are compatible along the localization').

(ii) The {\bf differential Ore condition} \index{differential Ore condition}
{\it is satisfied}:
\[\fbox{$ \forall t \in S, \,\forall r \in R, \,
\exists s\in S, \,\exists \omega \in \Omega^1(R),\,\,\,\,\,
s\,dr = \omega t.
$}\]
Proof.} (i) $\Rightarrow$ (ii).
If $S^{-1}\Omega^1(R) \equiv S^{-1}R\otimes_R \Omega^1(R)$ is a
$S^{-1}R$-bimodule then
$(dr) t^{-1} \in S^{-1} \Omega^1(R)$ for $t \in S$,
$r \in R$. All the elements in
$S^{-1} \Omega^1(R)$ are of the form
$s^{-1} \omega$ where $s\in S$ and $\Omega \in \Omega^1(R)$.
Hence $\exists s\in S, \,\exists \omega \in \Omega^1(R)$
such that $s\,dr = \omega t$ in the localization.
By~\luse{torsionOre} this means
$s\,dr = \omega t + \omega'$ in $\Omega^1(R)$,
where $s'\omega' = 0$ for some $s' \in S$.
Pre-multiplying by $s'$ we obtain $(s's) dr = (s' \omega) t$,
with required form.

(ii) $\Rightarrow$ (i).
The right $S^{-1}R$ action if it exists is clearly forced by
\begin{equation}\label{eq:extrightaction}
 s_1^{-1} a d(r) t^{-1} b = s_1^{-1} a s^{-1} \omega b \end{equation}
for $s,\omega$ chosen as above.
On the other hand, if~(\ref{eq:extrightaction}) holds,
this right action does extend the right $R$-action.
One has to prove that~(\ref{eq:extrightaction}) can be taken as a definition
of right $S^{-1}R$-action (compatible with the left action),
i.e. it does not depend on choices. If we choose $s',\omega'$ such that
$s'd(r) = \omega' t$ then $s^{-1}\omega t = (s')^{-1}\omega' t$.
As $t$ does not annihilate from the right,
$s^{-1}\omega = (s')^{-1} \omega'$. Other cases are left to the reader.
Hence $S^{-1}\Omega^1(R)$ is a bimodule; its elements are of the form
$s^{-1}adb$.

To prove that it is sufficient, define $d_S$ from $d$ by
the generalization of
formula~(\ref{eq:derOnLoc}) by
$\bar\partial$ and $\partial$
replaced by $d_S$ and $d$
and proceed with the rest of the proof
as in~\luse{extendingDiffCalc} --
all the calculations there make sense.

\minsection{Gabriel filter ${\cal L}_S$ for any $S \subset R$}

\ppt A \wind{lattice} is a poset $(W,\succ)$ such that for any two elements
$z_1,z_2$ the least upper bound $z_1 \vee z_2$
and the greatest lower bound $z_1 \wedge z_2$ exist.
In other words, the binary operations
{\it meet}\index{meet} $\wedge$ and {\it join}\index{join} $\vee$
are everywhere defined.
A poset is {\bf bounded} \index{bounded poset}
if it contains a maximum and a minimum element,
which we denote $1$ and $0$ respectively.
A ('proper') \wind{filter} in a bounded lattice $(W,\succ)$
is a subset $\filtf \subset W$ such that
$1\in \filtf$, $0 \notin \filtf$,
$(z_1,z_2 \in \filtf\Rightarrow z_1 \wedge z_2 \in \filtf)$
and $(z \in \filtf, z' \succ z \Rightarrow z' \in \filtf)$.

E.g. in any bounded lattice $(W,\succ)$, given $m \in W$,
the set ${}^m W$ of all $n \succ m$ is a filter.

\ppt {\it Notation.} Given a left ideal $J \in I_l R$
and a subset $w \subset R$ define
\[\fbox{$(J : w) := \{ z \in R \,|\, zw \subset J \}$ } \]
Then $(J:w)$ is a left ideal in $R$. If $w =: K$ is also a left ideal,
then $(J : K)$ is 2-sided ideal. In particular, if $w = K = R$,
then $(J:R)$ is the maximal 2-sided ideal contained in $J$.
For $r \in R$ we write $(J:r)$ for $(J,\{r\})$.

Given subsets $v,w \subset R$, set $((J:v):w)$
contains precisely all $t_1$ such that
$t_1 w \subset (J:v)$, i.e. $t_1 w v \subset J$.
Hence $((J:v):w) = (J: wv)$.

\ppt {\bf Preorders on left ideals}.
Let $I_l R$ be the set of all left ideals in a ring $R$.
It is naturally a preorder category
with respect to the inclusion preorder\index{inclusion preorder}.
This category is a lattice.
For the localization and spectral questions another partial order $\succ$ on
$I_l R$ is sometimes better: $K \succ J$
(category notation: $J \to K$) iff either $J \subset K$,
or there exist a finite subset $w \subset R$ such that $(J:w) \subset K$.
Any filter in $( I_l R, \succ )$ is called a {\bf uniform filter}.
~\index{filter!uniform}\ldef{defUF}

\ppt Let $R$ be a unital ring and
$S \subset R$ a multiplicative set. Consider
\begin{equation}\label{eq:Gabriel-filter-for-S}
 \filtf_S := \{ J\mbox{ left ideal in } R\,|\, \forall r,\,
(J : r) \cap S \neq \emptyset   \}\subset I_l R.
\end{equation}
We make the following observations:
\begin{itemize}
\item As $(R:r) = R$, $R \in \filtf_S$.
\item Suppose $J,K \in \filtf_S$. Given $r \in R$, $\exists s,t$,
such that $s \in (J : r) \cap S$ and $t \in (K : s r) \cap S$.
Hence $t s r \in J \cap K$. Set $S$ is multiplicative, hence
$ts \in S$ and $ts \in (J : r) \cap (K:r) \cap S = (J \cap K : r)\cap S$.
Thus $J \cap K \in \filtf_S$.
\item $(J : r) \cap S \neq \emptyset$ then, a fortiori,
$(K:r)\cap S \neq \emptyset$ for $K \supset J$.
\item If $J \in \filtf_S$ then $\forall r$ $\,(J : r) \cap S \neq \emptyset$.
In particular, this holds with $r$ replaced by $rr'$.
Using $((J:r):r') = (J : r'r)$
we see that $(J:r) \in \filtf_S$ for all $r \in R$.
\item If $\forall r' \in R$ $(J : r') \cap S \neq \emptyset$ and
$((J':j) : r) \cap S \neq \emptyset$ for all $j \in J$, $r \in R$,
then $\exists s \in S$
such that $srj \in J'$ and $\exists s' \in S$ such that $s'r' \in J$.
In particular for $r = 1$ and $j = s'r'$ we have $ss'r' \in J$.
Now $ss' \in S$ and $r'$ is arbitrary so $J' \in \filtf_S$.
\end{itemize}
These properties can be restated as the axioms for a {\bf Gabriel filter}
\index{Gabriel filter} \index{radical filter}
$\filtf\subset I_l R$ (synonyms ``radical set'', ``radical filter'',
``idempotent topologizing filter''):
\begin{itemize}
\item (F1) $R \in \filtf$ and $\emptyset \notin \filtf$.
\item (F2) If $J,K \in \filtf$, then $J \cap K \in \filtf$.
\item (F3) If $J \in \filtf$ and $J\subset K$ then $K \in \filtf$.
\item (UF) $J \in \filtf \Leftrightarrow \,
(\forall r \in R,\,(J:r) \in \filtf)$.
\item (GF) If $J \in \filtf$ and $\forall j \in J$
the left ideal $(J':j)\in \filtf$, then $J' \in \filtf_S$.
\end{itemize}
Axioms (F1-3) simply say that a set $\filtf$ of ideals
in $R$ is a filter in $(I_l R, \subset)$.
Together with (UF) they exhaust the axioms
for a uniform filter 
(cf.~\luse{defUF}).
\def\glc{Q} 
Axioms (GF) and (UF) imply (F2):
If $j \in J$, $(I\cap J:j)=(I:j)\cap(J:j)=(I:j)\in\filtf$ by (UF).
Since $\forall j\in J$ $(I\cap J:j)\in \filtf$, 
(GF) implies $I\cap J\in\filtf$. \newline
(GF) $\&$ (F1) imply (F3): $(\forall j \in J\subset K)$
$(K:j)=R\in\filtf$, hence $K\in\filtf$.

There are examples of Gabriel filters $\filtf$, even for commutative $R$,
which are {\it not} of the form $\filtf_S$ for a multiplicative $S \subset R$.
Moreover, for rings without unity (F1-3,UF,GF) still make sense, whence
a good notion of a multiplicative set and filters $\filtf_S$ fails to exist.

Notice that if a multiplicative set $S$ satisfies the left Ore
condition, then $\filtf_S = \filtf'_S := \{ J \mbox{ is left ideal }
\,|\, J \cap S \neq \emptyset \}$. Namely, $(J:1) \cap S = J \cap S$
for {\it any} S, hence $\filtf_S \subset \filtf'_S$;
and the {\it left} Ore condition implies that given
an element $s \in J \cap S$ and $r \in R$ we can find $s' \in S$,
and $r' \in R$ with $s'r = r's \in r'J \subset J$,
hence $s' \in (J:r) \cap S$; hence $\filtf'_S \subset \filtf_S$.

\ppt {\bf Exercise.} Check that the intersection of any family of
Gabriel filters is a Gabriel filter.

{\bf Remark}: this is {\it not} always true
for the union: (GF) often fails.

\ppt For given $R$-module $M$ and a filter $\filtf$ in $(I_l R, \subset)$,
the inclusions $J\hookrightarrow J'$ induce 
maps ${\rm Hom}_R (J',M)\to{\rm Hom}_R (J,M)$ for any $M$,
hence we obtain an inductive system of \abelian{} groups.
The inclusion also induce the projections
$R/J\to R/J'$ and hence, by composition, 
the maps ${\rm Hom}_R (R/J',M)\to{\rm Hom}_R (R/J,M)$.
This gives another inductive system of \abelian{} groups.
If a filter $\filtf$ is uniform, 
we consider the same systems and limits of groups 
(without new morphisms), and use (UF)
as ingenious device to define the $R$-module structure on them.

\ppt {\bf Proposition.} {\it Let $\filtf$ be a uniform filter 
and $M$ a left $R$-module.

(i) The inductive limit of \abelian{} groups taken over
downwards directed family of ideals
\[
  H_\filtf(M) :=
{\rm lim}_{J \in \filtf} {\rm Hom}_R (J,M)
\]
has a canonical structure of an $R$-module. $H_\filtf$ extends to
an endofunctor.

(ii) The abelian subgroup
\[ \sigma_\filtf(M) := \{ m \in M \,|\, \exists J \in \filtf, \,J m = 0\}
\subset M
\]
is a $R$-submodule of $M$.

(iii) If $f:M\to N$ is a map of $R$-modules, 
${\rm Im}\,f|_{\sigma_\filtf(M)} \subset \sigma_\filtf(N)$, hence
the formula $f\mapsto\sigma_\filtf(f):= f|_{\sigma_\filtf(M)}$ extends
$\sigma_\filtf$ to a subfunctor of identity. 

(iv) The inductive limit of abelian group 
$\sigma'_\filtf(M) := {\rm lim}_{J \in \filtf} {\rm Hom}_R (R/J,M)$
has a structure of a left $R$-module.

(v) If $1 \in R$ then the
endofunctors $\sigma_\filtf$ and $\sigma'_\filtf$ (on the
categories of modules $M$ with $1_R m = m$, where $m \in M$)
are equivalent.
}

{\it Proof.} (i) Given $f \in  H_\filtf(M)$, represent it as
$f_J$ in ${\rm Hom}_R (J,M)$ for some $J\in\filtf$.
By (UF), $\forall r\in R$, $(J:r)\in \filtf$.
The rule $x \mapsto f_J(xr)$ defines a map $(rf)_{(J:r)}$ in
${\rm Hom}_R ((J:r),M)$ which we would like to
represent class $rf$.
Suppose we have chosen another representative $f_I$, then
there is $K\in\filtf$, $K \subset I\cap J$, such that 
$f_I|_K = f_J|_K =: h$. Then $(K:r)\subset (I\cap J:r) = (I:r)\cap(J:r)$
and the map $x\mapsto h(xr) : K\to M$ agrees with $(rf)_{(J:r)}|_K$
and $(rf)_{(I:r)}|_K$ hence the class $rf$ is well defined.

This is a left action:
$((rr')f)_{(J : rr')}(x) = f_{J}(xrr') = (r'f)_{(J:r')}(xr) 
= (r(r'f))_{((J:r'):r)}(x) = (r(r'f))_{(J : rr')}(x)$. We used
$((J:r'):r) = (J : rr')$.

(ii) Suppose $m \in \sigma_\filtf(M)$, i.e. $Jm = 0$ for some $J \in \filtf$.
For arbitrary $r \in R$ the ideal $(J:r) \in \filtf$ by (UF).
Let $k \in (J:r)$. Then $kr \in J$, hence $krm = 0$. This is true for any
such $k$, hence $(J:r)rm = 0$ and $rm \in \sigma_\filtf(M)$.
As $r$ was arbitrary, $R\sigma_\filtf(M) \subset \sigma_\filtf(M)$.

(iii) If $m \in \sigma_\filtf(M)$ then $0 = f(0) = f(Jm) = Jf(m)$
for some $J$ in $\filtf$. Hence $f(m)\in \sigma_\filtf(N)$.

(iv) Let $r \in R$ and
$f \in {\rm lim}\,{\rm Hom}_R (R/J,M)$. 
Take a representative $f_J \in {\rm Hom}_R (R/J,M)$.
Let $(rf)_{(J:r)} \in {\rm Hom}_R (R/(J:r),M)$ be given by
$(rf)_{(J:r)}(r' + (J:r)) = f_J(r'r + J)$. 
This formula does not depend on $r'$ because 
changing $r'$ by an element $\delta r' \in (J:r)$ changes
$r'r$ by an element $(\delta r') r$ in $(J:r)r \subset J$.
Suppose $f_I \sim f_{J}$. In this situation, with projections
as connecting morphisms, this means that $f_I(x + I) = f_J(x+J)$
for all $x\in R$, and in particular for $x = r'r$, hence
$(rf)_{(J:r)} \sim (rf)_{(I:r)}$ and $rf$ is well defined.


Finally, $f \mapsto rf$ is a left $R$-action.
Indeed, for all $r,r',t \in R$,
\[\begin{array}{lcl} ((rr')f)_{(J:rr')}(t+ (J:rr')) 
&=& (f_J)(trr' + J)
\\&=& (r'f)_{(J:r')}( tr + (J:r'))
\\&=& (r(r'f))_{((J:r'):r)}(t + ((J:r'):r))
\\ &=& (r(r'f))_{(J:rr')}(t + (J:rr'))\end{array}\]
If $R$ and $M$ are unital, then $1_R f = f$ as well.

(v) To make the statement precise,
we should first extend $\sigma_\filtf$ to a functor by
defining it on morphisms as well ($\sigma'_\filtf$ is obviously
a functor as the formula on object is explicitly written in terms of a
composition of functors applied on $M$).
As $\sigma_\filtf(M) \subset M$, it is sufficient
to show that $f(\sigma_\filtf(M)) = \sigma_\filtf(f(M))$ and then
define $\sigma_\filtf(f) := \sigma_\filtf \circ f$.
Element $f(m) \in f(\sigma_\filtf(M))$
iff $Jm = 0$ for some $J \in \filtf$. This is satisfied
iff $f(Jm)= Jf(m) = 0$, i.e.
$f(m) \in  \sigma_\filtf(f(M))$.

The equivalence $\nu : \sigma'_\filtf \Rightarrow \sigma_\filtf$
is given by $\nu_M([f_J]) := f_J(1_R + J)\in \sigma_\filtf(M)$
(because $Jf_J(1_R + J) = f_J(0 + J) = 0$),
with inverse $m \mapsto [\theta^{(m)}_J]$ where
$\theta_J^{(m)} : r + J \mapsto rm$ and any $Jm = 0$
(such $J$ exists and the formulas for $\theta^{(m)}_J$
for different $J$ agree, 
hence {\it a fortiori} define a limit class).
Starting with $[f_J]$ with $m := f_J(1_R + J)$ and
$\theta_J^{(m)} : r+ J \mapsto rf_J(1_R + J) = f_J(r+J)$, hence
$\theta_J^{(m)} = f_J$. Other way around, start with $m\in \sigma_\filtf(M)$,
then $\nu_M(\theta_J^{(m)}) = \theta_J^{(m)}(1_R + J) = 1_R m = m$.
Hence we see that each $\nu_M$ is an isomorphism of modules.
The reader may check that $\nu$, $\nu^{-1}$ are natural transformations.

\vskip .026in
\ppt If ${\cal A}$ is any \abelian{} category, then a
subfunctor $\sigma$ of the identity
(i.e. $\sigma(M)\subset M$ and $\sigma(f)|_{\sigma(M)} = f|_{\sigma(M)}$,
cf.~\luse{subobject}) with the property
$\sigma(M/\sigma(M)) = 0$ is called a \wind{preradical} in ${\cal A}$.
A preradical $\sigma$ in $R-{\rm Mod}$ is left exact iff $J\subset K$ implies 
$\sigma(J)=\sigma(K)\cap J$. A \wind{radical} is a left exact preradical. 

\ppt \ldef{idemrad} {\bf Proposition.} {\it If $\filtf$ is
Gabriel filter, $\sigma_\filtf$ is an \wind{idempotent radical} in
the category of left $R$-modules, i.e. it is a radical
and $\sigma_\filtf \sigma_\filtf = \sigma_\filtf$.
}


\ppt To any Gabriel filter $\filtf$, one associates a localization
endofunctor $Q_{\filtf}$ on the category of left modules by
the formula
\begin{equation}\label{eq:GfQ} Q_\filtf(M) := H_\filtf (M/\sigma_\filtf(M)) =
{\rm lim}_{J \in \filtf} {\rm Hom}_R (J,M/\sigma_\filtf(M)).\end{equation}

Left multiplication by an element $r \in R$ defines a class
$[r] \in Q_\filtf(R)$. There is a unique ring
structure on $Q_\filtf(R)$, such that the
correspondence $i_\filtf : r \mapsto [r]$
becomes a ring homomorphism $i_\filtf : R \rightarrow Q_\filtf(R)$.

Notice that~(\ref{eq:GfQ}) generalizes the RHS of Deligne's formula,
\luse{delgab}{\bf a}.

\ppt Not only every Gabriel filter defines an idempotent radical,
but also~(\cite{JaraVerschorenVidal}):

{\bf Proposition}. {\it Every radical defines a Gabriel filter by the rule
\[ \filtf_\sigma := \{ J \mbox{ left ideal in } R\,|\, \sigma(R/J) = R/J \}.\]
More generally, if $M$ be a left $R$-module and $\sigma$ a radical, define
\[ \filtf_{M,\sigma} := \{ L \mbox{ left }R\mbox{-submodule in } M\,|\,
\sigma(M/L) = M/L \}.\]
Then $\filtf_M := \filtf_{M,\sigma}$ satisfies the following properties
\begin{itemize}
\item (GT1) $M \in \filtf_{M}$.
\item (GT2) If $L,K \in \filtf_{M}$, then $L \cap K \in \filtf$.
\item (GT3) If $L \in \filtf_M$, $K\subset M$ a left submodule
and $L\subset K$, then $K \in \filtf_{M}$.
\item (GT4) If $J \in \filtf_M$ and $K \in \filtf_J$
the left submodule $K \in \filtf_M$.
\end{itemize}
}

\ppta When we restrict to the {\it idempotent radicals},
then the rule $\sigma \mapsto \filtf_\sigma$
gives a {\it bijection} between the idempotent radicals
and Gabriel filters.

\minsection{Localization in \abelian categories}

The language of Gabriel filters is not suited for some other
categories where additive localization functors are useful.
Subcategories closed with
respect to useful operations (e.g. extensions of objects)
are often used as the localization data,
particularly in \abelian and triangulated categories.
We confine ourselves just to a summary of basic notions
in \abelian setting and comment on the connection to the
language of Gabriel filters, as a number of references is
available~(\cite{Borceux,Gab:catab,GZ,JaraVerschorenVidal,
Popescu,Ros:localalg}).

\ppt Let ${\cal A}$ be an additive category.
Let ${\cal P}$ be a full subcategory of ${\cal A}$.
Define the left and right orthogonal to ${\cal P}$
\index{left orthogonal subcategory}
\index{right orthogonal subcategory}
to be the full subcategories
${}^\perp{\cal P}$ and ${\cal P}^\perp$ consisting
of all objects $A \in {\cal A}$ such that ${\cal A}(P,A) = 0$
(resp. ${\cal A}(A,P) = 0$) for all $P \in {\cal P}$.
Zero object is the only object in ${\cal P} \cap {}^\perp{\cal P}$.
It is clear that taking (left or right) orthogonal
reverses inclusions and that ${\cal P} \subset {}^\perp({\cal P}^\perp)$
and ${\cal P}\subset ({}^\perp{\cal P})^\perp$. We leave as an
exercise that ${\cal P}^\perp = ({}^\perp({\cal P}^\perp))^\perp$
and ${}^\perp{\cal P} = {}^\perp(({}^\perp{\cal P})^\perp)$.

\ppt A {\bf thick subcategory} \index{thick subcategory}
of an \abelian category ${\cal A}$
is a strictly full subcategory ${\cal T}$ of ${\cal A}$ which is closed under
extensions, subobjects and quotients.
In other words, an object $M'$ in a short sequence
$0\rightarrow M \rightarrow M' \rightarrow M''\rightarrow 0$
in ${\cal A}$ belongs to ${\cal T}$ iff $M$ and $M''$ do.

Given a pair $({\cal A}, {\cal T})$ where ${\cal A}$ is
\abelian and ${\cal T}\subset {\cal A}$ is thick,
consider the class
\[ \Sigma({\cal T}) := \{ f  \,|\,
{\rm Ker}\,f \in {\rm Ob}\,{\cal T},\mbox{ and }
{\rm Coker}\,f \in {\rm Ob}\,{\cal T}
\}.\]
The quotient category ${\cal A}/{\cal T}$ is defined as follows.
${\rm Ob}\,{\cal A}/{\cal T} = {\rm Ob}\,{\cal A}$
and ${\rm Mor}\,{\cal A}/{\cal T} := {\rm Mor}\,{\cal A}
\coprod \Sigma^{-1}({\cal T})$, where $\Sigma^{-1}({\cal T})$
is the class of formal inverses of morphisms $f \in \Sigma$;
impose the obvious relations.
${\cal A}/{\cal T}$ is additive in a unique way making the
quotient functor additive. In fact~(\cite{Gab:catab,GZ}),
it is abelian.

{\bf Proposition.} ({\sc Grothendieck}~\cite{Tohoku}) {\it
Let ${\cal T}$ be a thick subcategory in ${\cal A}$ and
$\Sigma({\cal T})$  as above.
Then $\Sigma$ is a left and right calculus of fractions in ${\cal A}$
and ${\cal A}[\Sigma({\cal T})]^{-1}$ is naturally isomorphic to
${\cal A}/{\cal T}$.
}

A thick subcategory ${\cal T}$
is a {\bf localizing subcategory} \index{localizing subcategory}
if the morphisms which are {\it invertible} in the quotient
category ${\cal A}/{\cal T}$ are precisely the images of
the morphisms in $\Sigma({\cal T})$.

Every exact localization functor $T^*: {\cal A}\to {\cal B}$
(i.e. an exact functor with fully faithful right adjoint $T_*$)
of an \abelian category ${\cal A}$ is the localization at the localizing
subcategory $\Sigma$ consisting of those morphisms $f$
such that $T^*f$ is either kernel or a cokernel morphism
of an invertible morphism in ${\cal B}$.

If $T^* : {\cal A}\to {\cal B}$ is any exact localization functor,
then set ${\cal T} := {\rm Ker}\,T^*$ to be the full subcategory
of ${\cal A}$ generated by all objects $X$ such that $T^*(X) = 0$.
Then $T^*$ factors uniquely as
${\cal A}\stackrel{{\cal Q}^*}\to {\cal A}/{\cal T} \to {\cal B}$
where ${\cal Q}^*$ is the natural quotient map.

More than one thick subcategory may give the same
quotient category, and that ambiguity is removed if we consider
the corresponding localizing subcategories instead~(\cite{Popescu}).

A composition of localization functors corresponds to {\bf Gabriel
multiplication} $\bullet$ on thick subcategories.
For any two subcategories ${\cal B}$, ${\cal D}$
of an \abelian category ${\cal A}$ one defines $D\bullet B$
to be the full subcategory of ${\cal A}$
consisting of precisely those $A$ in ${\cal A}$ for
which there is an exact sequence $0 \to B\to A \to D\to 0$ with
$B$ in ${\rm Ob}\,{\cal B}$ and $D$ in ${\rm Ob}\,{\cal D}$.
In categories of modules one can redefine Gabriel
multiplication in terms of radical filters, cf.~(\cite{Rosen:book}).

\ppt In this article, we often view exact localizations
(and quotient categories, cf.~\luse{delgab}{\bf b})
as categorical analogues of open spaces.
Their complements should then be the complementary data
to the quotient categories, and such data are localizing subcategories.
A more precise and detailed discussion of those subcategories,
which may be considered as subschemes and closed subschemes,
may be found in~\cite{LuntsRosMP}, Part I
and~\cite{Rosen:book, Smith:3spaces}.
Cf.~the notion of a (co)reflective subcategory in~\luse{idempmonad}.

Thus, in our view, it is geometrically more appealing to
split the data of a category to a localizing subcategory and a quotient
category, than into two {\it sub}categories.
However, the latter point of view is more traditional, under the
name of ``torsion theory' and has
geometrically important analogues for triangulated categories.
A {\bf torsion theory}
~(\cite{Borceux,JaraVerschorenVidal})
in an \abelian category ${\cal A}$ is a pair $({\cal T}, {\cal F})$
of subcategories of ${\cal A}$ closed under isomorphisms and such that
${\cal F}^\perp = {\cal T}$ and ${}^\perp{\cal T} = {\cal F}$.


For any idempotent radical $\sigma$ in ${\cal A}$ (\luse{idemrad}),
the class ${\cal T}_\sigma$ of $\sigma$-{\bf torsion objects}
and the class ${\cal F}_\sigma$ of $\sigma$-{\bf torsion free objects}
are defined by formulas
\[ {\cal T}_\sigma = \{ M \in {\rm Ob}\, {\cal A}\,|\, \sigma(M) = M\},
\,\,\,\,\,\,\,
{\cal F}_\sigma = \{ M \in {\rm Ob}\, {\cal A}\,|\, \sigma(M) = 0\}.\]
This pair $({\cal T}_\sigma,{\cal F}_\sigma)$
is an example of a torsion theory and ${\cal T}_\sigma$
is a thick subcategory of ${\cal A}$.
Not every torsion theory corresponds to a radical, but hereditary
theories do. That means that a subobject of a torsion object is torsion.
The Cohn localization of the next section is not necessarily hereditary,
\ldef{hereditary2}\index{hereditary torsion theory}
but it is always a torsion theory as shown there.

\minsection{Quasideterminants and Cohn localization}

{\bf Notation.} Let $M^n_{m}(R)$ be the set of all $n\times m$ matrices
over a (noncommutative) ring $R$, so that
$M_n(R):= M^n_n(R)$ is a ring as well.
Let $I,J$ be the ordered tuples
of row and column labels of $A = (a^i_j)\in M^n_m(R)$ respectively.
For subtuples $I'\subset I,J'\subset J$ and $A = (a^i_j)\in M^n_m(R)$,
denote by $A^{I'}_{J'}$ the submatrix of $A^I_J := A$ consisting only
of the rows and columns with included labels; e.g.
$A^{\{i\}}_{\{j\}} = a^i_j$ is the entry in $i$-th row
and $j$-th column. When $I$ is known and $K \subset I$,
then $|K|$ is the cardinality of $K$ and $\hat{}$ is the symbol for
omitting, i.e. $\hat{K} = I\backslash K$
is the complementary $(|K|-|I|)$-tuple.

We may consider the $r$-tuple $\tilde{I} = (I_1,\ldots, I_r)$ of
sub-tuples which partitions the $n$-tuple $I =
(i_1,i_2,\ldots,i_n)$, i.e. $I_k$ are disjoint and all labels from
$I$ are included; then $|\tilde{I}|:= r$. Given $\tilde{I},
\tilde{J}$ form the corresponding {\bf block matrix} \index{block
matrix}\index{partitioned matrix} in $M^{\tilde{I}}_{\tilde{J}}$
out of $A$, i.e. the $|\tilde{I}|\times |\tilde{J}|$ matrix
$A^{\tilde{I}}_{\tilde{J}}$ whose entries are matrices
$A^{\tilde{I}_k}_{\tilde{J}_l} := A^{I_k}_{J_l}$ cut-out from $A$
by choosing the selected tuples. Forgetting the partition gives
the canonical bijection of sets
$M^{\tilde{I}}_{\tilde{J}}\rightarrow M^I_J$. The multiplication of
block matrices is defined by the usual matrix multiplication
formula $(AB)^{I_i}_{J_j} = \sum_{l = 1}^r A^{I_i}_{K_l}
A^{K_l}_{J_j}$ if $AB$ {\it and} the sizes of subtuples for
columns of $A$ and rows of $B$ match. One can further nest many
levels of partitions (block-matrices of block-matrices $\dots$).
Some considerations will not depend on whether we consider
matrices in $R$ or block matrices, and then we'll just write
$M^I_J$ etc. skipping the argument. More generally, the labels may
be the objects in some \abelian category ${\cal A}$, and entries
$a^i_j \in {\cal A}(i,j)$; $I$ will be the sum $\oplus_{i \in I}
i$, hence $A : I \rightarrow J$. Ring multiplication is replaced
by the composition, defined whenever the labels match.

{\bf Observation.} {\it Multiplication of block matrices commutes with
forgetting (one level) of block-matrix structure.}
In other words we may multiply in stages (if working in ${\cal A}$
this is the associativity of $\oplus$). Corollaries:

{\it (i)  if $\tilde{I}= \tilde{J}$
then $M_{\tilde{I}}(R):= M^{\tilde{I}}_{\tilde{I}}(R)$
is a ring.

(ii) We can invert matrices in stages as well ('heredity').

(iii) The same for linear equations over noncommutative rings.
}

\vskip .02in
Any pair $(i,j) \in I\times J$
determines partitions $\tilde{I} = (i,\hat{i})$ and$\tilde{J} = (j,\hat{j})$.
For each $A$ in $M^I_J$ it induces a $2 \times 2$ block-matrix
$A^{\tilde{I}}_{\tilde{J}}$.
Reader should do the exercise
of inverting that block matrix (with noncommutative entries),
in terms of the inverses of blocks.
As we will see, the $(i,j)$-quasideterminant of $A$
is the inverse of the $(j,i)$-entry of $A^{-1}$
if the latter is defined; though it may be defined
when the latter is not.

\ppt The $(i,j)$-th quasideterminant $|A|_{ij}$ of $A$ is
\begin{equation}\label{eq:quadetdef}\fbox{$ |A|_{ij} = a^i_j -
\sum_{k \neq i, l\neq j} a^i_l  (A^{\hat{i}}_{\hat{j}})^{-1}_{lk} a^k_j
$}\end{equation}
provided the right-hand side is defined (at least in the sense
of evaluating a rational expression, which will be discussed below).
In alternative notation, the distinguished labels $ij$ may be replaced
by a drawing of a box around the entry $a^i_j$ as in
\[
\left|\begin{array}{ccc} a^1_1 & a^1_2 & a^1_3 \\
a^2_1 & a^2_2 & a^2_3 \\
a^3_1 & \fbox{$a^3_2$} & a^3_3
\end{array}\right| \equiv
\left|\begin{array}{ccc} a^1_1 & a^1_2 & a^1_3 \\
a^2_1 & a^2_2 & a^2_3 \\
a^3_1 & a^3_2 & a^3_3
\end{array}\right|_{32}
\]
At most $n^2$ quasideterminants of a given $A \in M_n(R)$
may be defined.
\ppt {\it If all the $n^2$ quasideterminants $|A|_{ij}$ exist
and are invertible then the inverse $A^{-1}$ of $A$ exist in
$A \in M_n(R)$ and }
\begin{equation}\label{eq:quasiInv}\index{$|A|_{ji}$}
(|A|_{ji})^{-1} = (A^{-1})^i_j.
\end{equation}

Thus we also have
\begin{equation}\label{eq:quasidet-recursive-def}
|A|_{ij} = a^i_j -
\sum_{k \neq i, l\neq j} a^i_l |A^{\hat{i}}_{\hat{j}}|^{-1}_{kl} a^k_j
\end{equation}
\ppt
Sometimes the RHS of~(\ref{eq:quadetdef}) makes sense while
(\ref{eq:quasidet-recursive-def}) does not. So for subtle existence
questions one may want to be careful with alternative formulas
for quasideterminants. Some identities are often proved using alternative
forms, so one has to justify their validity. Different expressions
differ up to~\wind{rational identities} (\cite{Amitsur:rat,Cohn:free}),
and under strong assumptions on the ring $R$ (e.g. a skewfield which is
of $\infty$ dimension over the center which is also infinite),
the rational identities induce a well-behaved equivalence
on the algebra of rational expressions and the results of calculations
extend in an expected way to alternative forms once they are proved
for one form having a nonempty domain of definition 
~(\cite{Amitsur:rat,Cohn:free, Skoda:nflag})

\ppt \ldef{quasiwarning}
On the other hand, the existence of inverse $A^{-1}$
does {\it not} imply the existence of quasideterminants.
For example, the unit $2 \times 2$ matrix ${\bf 1}_{2\times 2}$
over field ${\DDl Q}$ has only 2 quasideterminants, not 4. Or, worse,
matrix $\left(\begin{array}{cc}3 & 2 \\ 2 & 3\end{array}\right)$ over
commutative ring ${\DDl Z}[\frac{1}{5}]$ is invertible, but no single entry is
invertible, and in particular no quasideterminants exist.

\ppt Quasideterminants are {\it invariant under permutation} of rows or columns
of $A$ if we appropriately change the distinguished labels.

\ppt Suppose now we are given an equation of the form
\[ A x = \xi \]
where $A \in M_n(R)$ and $x,\xi$ are $n$-tuples of indeterminates and
free coefficients in $R$ respectively (they are column ``vectors'').
Then one can attempt to solve the system by finding the inverse of
matrix $A$ and multiply the equation by $A^{-1}$ from the left, or one
can generalize the Cramer's rule to the noncommutative setup.

Define thus $A(j,\xi)$ as the $n \times n$ matrix whose entries are
the same as of $A$ except that the $j$-th column is replaced by $\xi$.
Then the noncommutative \wind{left Cramer's rule} says
\[\fbox{$ |A|_{ij} x^j = |A(j,\xi)|_{ij} $} \]
and the right-hand side does not depend on $i$.

To see that consider first $n = 2$ case:
\[\begin{array}{l}
        a^1_1 x^1 + a^1_2 x^2 = \xi^1 \\
        a^2_1 x^1 + a^2_2 x^2 = \xi^2
\end{array}\]
Then
\[\begin{array}{lcl}
 |A|_{11} x^1 &=& a^1_1 x^1 - a^1_2 (a^2_2)^{-1} a^2_1 x^1 \\
&=& (\xi^1 - a^1_2 x^2) - a^1_2 (a^2_2)^{-1}a^2_1 x^1\\
&=& \xi^1 - a^1_2 (a^2_2)^{-1} a^2_2 x^2 - a^1_2 (a^2_2)^{-1}a^2_1 x^1\\
&=& \xi^1 - a^1_2 (a^2_2)^{-1} \xi^2 = |A(1,\xi)|_{11}.
\end{array}\]
The general proof is exactly the same, just one has to understand which
indices are included or omitted in the sums involved:
\[\begin{array}{lcl}
 |A|_{ij} x^j &=& a^i_j x^j - \sum_{k \neq j, l \neq i}
a^i_k (A^{\hat{i}}_{\hat{j}})^{-1}_{kl} a^l_j x^j \\
&=& (\xi^i - \sum_{h \neq j} a^i_h x^h) - \sum_{k \neq j, l \neq i}
a^i_k (A^{\hat{i}}_{\hat{j}})^{-1}_{kl} a^l_j x^j\\
&=& \xi^i - \sum_{h \neq j,k \neq j, l \neq i}
a^i_k (A^{\hat{i}}_{\hat{j}})^{-1}_{kl} a^l_h x^h
- \sum_{k \neq j, l \neq i}
a^i_k (A^{\hat{i}}_{\hat{j}})^{-1}_{kl} a^l_j x^j\\
&=& \xi^i - \sum_{1\leq h\leq n, k \neq j, l \neq i}
a^i_k (A^{\hat{i}}_{\hat{j}})^{-1}_{kl} a^l_h x^h \\
&=& \xi^i - \sum_{1\leq h\leq n, k \neq j, l \neq i}
a^i_k (A^{\hat{i}}_{\hat{j}})^{-1}_{kl} a^l_h x^h \\
&=& |A(j,\xi)|_{ij}.
\end{array}\]

Similarly consider equation $\sum_k y^k B^l_k = \zeta^l$.
Apparently the individual coefficients multiply $y^k$ from the
right, but the combinatorics of matrix labels is organized as if
we multiply $By$ (alas, otherwise the rule of writing upper
indices for rows would force us to write such equations
upside-down!). The canonical antiisomorphism $R \rightarrow R^{\rm
op}$ clearly sends any quasideterminant into the quasideterminant
of the {\it transposed matrix}. Hence the left Cramer's rule
implies the \wind{right Cramer's rule}
\[\fbox{$  y^j |B^T|_{ji}= |(B(j,\zeta))^T|_{ji}. $}\]

\ppt {\bf Row and column operations.}~\index{row and column operations}
Ordinary determinants do not change
if we add a multiple of one row to another, and similarly for the columns.

We have to distinguish between left and right linear combinations.

{\it If $|A|_{ij}$ is defined and $i\neq l$,
then it is unchanged under left-row operation }
\[ A^l \rightarrow A^l + \sum_{s \neq l} \lambda_s A^s \]

{\it Proof.} We may assume $i = 1$. Define the row matrix
\[ \vec{\lambda} = (\lambda_2,\ldots,\lambda_n). \]
Then $\vec\lambda T = \sum_{s\neq k} \lambda_s T^s$
for any matrix $T$ with row-labels $s = 2,\ldots,n$.
Then $\Lambda T = \sum_{s\neq k} \lambda_s T^s$.
Assume the matrix $A$ is in the block-form written as
\[ \mat22{a}{\vec{b}}{\vec{c}^T}{D} \]
with $a$ of size $1\times 1$. Then
\[\begin{array}{lcl}
\left| \begin{array}{cc}
\fbox{$a + \vec{\lambda} \vec{c}^T$}&
{\vec{b} + \vec{\lambda} D}\\
  {\vec{c}^T}&{D}
\end{array}\right|
&=& a + \vec{\lambda} \vec{c}^T
- (\vec{b} + \vec{\lambda} D)\,D^{-1}\,\vec{c}^T\\
&=& a - \vec{b} D^{-1}\vec{c}^T.
\end{array}\]

If we multiply the $l$-th row from the left by an invertible element $\mu$
then the quasideterminant $|A|_{ij}$ won't change for $i \neq l$ and
will be multiplied from the left by $\mu$ if $i=l$.
Actually, more generally, left multiply the $i$-th row by $\mu$ and
the block matrix consisting of other rows by invertible
square matrix $\Lambda$ (i.e. other rows can mix among themselves,
and scale by different factors):
\[
A \rightarrow \left(\begin{array}{cc}\mu& 0\\ 0 &\Lambda \end{array}\right)
\]
Then $|A|_{ij}$ gets left-multiplied by $\mu$:
\[\begin{array}{lcl}
\left| \begin{array}{cc}
\fbox{$\mu a$} & {\mu \vec{b}}\\
  {\Lambda \vec{c}^T} & {\Lambda D}
\end{array}\right|_{11}
&=& \mu a
- \mu \vec{b}\,(\Lambda D)^{-1}\,\Lambda \vec{c}^T\\
&=& \mu\, (a - \vec{b} D^{-1}\vec{c}^T) = \mu\,|A|_{ij}.
\end{array}\]

\ppt \ldef{jacobi-quasi}
\wind{Jacobi's ratio theorem}. (\cite{GelRet:ncr})
{\it Let $A$ be a matrix with possibly
noncommutative entries such that the inverse $B = A^{-1}$ is defined.
Choose some row index $i$ and some column index $j$.
Make a partition of the set of row indices as $I \cup \{ i\} \cup J$
and a partition of the set of column indices as $I' \cup \{j \} \cup J'$,
with the requirements ${\rm card}\, I = {\rm card}\,I'$
and ${\rm card}\, J = {\rm card}\,J'$. Then
}
\[
(| A_{I \cup \{i\}, I' \cup \{j\}} |_{ij})^{-1} =
| B_{J' \cup \{j\}, J \cup \{i\}}|_{ji}.
\]
{\it Proof.} (\cite{KrobLec}) The block decomposition of matrices
does not change the multiplication, i.e. we can multiply
the block matrices and then write out the block entries
in detail, or we can write the block entries of the multiplicands
in detail and then multiply and we get the same result.
In particular, as $A = B^{-1}$, the block-entries of $A$ can be
obtained by block-inversion of $B$.

After possible permutation of labels, we may find
the block-entry of the matrix $A = B^{-1}$ at the intersection of rows
$I \cup \{i\}$ and columns $I' \cup \{j\}$ by means of block-inverting
the block matrix
\[
 A = \left(\begin{array}{cc}
A_{I \cup \{i\}, I' \cup \{j\}}
& A_{I \cup \{i\}, J' } \\
A_{J, I' \cup \{j\}}
& A_{J, J'} \\
\end{array}\right)
\]
Then $A_{I \cup \{i\}, I' \cup \{j\}} = (B_{I \cup \{i\}, I' \cup \{j\}}
- B_{I \cup \{i\}, J'}(B_{JJ'})^{-1} B_{J,I \cup \{i\}})^{-1} $
or, equivalently,
\[
(A_{I \cup \{i\}, I' \cup \{j\}})^{-1} =
B_{I \cup \{i\}, I' \cup \{j\}}
- B_{I \cup \{i\}, J'}(B_{JJ'})^{-1} B_{J,I' \cup \{j\}}
\]
This is a matrix equality, and therefore it implies the equality
of the $(i,j)$-th entry of both sides of the equation.
We obtain
\[((A_{I \cup \{i\}, I' \cup \{j\}})^{-1})_{ij} =
b_{ij} - \sum_{k \in J', l \in J}
b_{i, k}(B_{JJ'})^{-1}_{kl} b_{lj}.
\]
Finish by applying the formula $(|C|_{ji})^{-1} = (C^{-1})_{ij}$
at LHS.

\ppt \wind{Muir's law of extensionality}.
(\cite{GelRet:ncr,GelRet:quasiI, GelRet:quasi})
Let an identity ${\cal I}$ between
quasiminors of a submatrix $A^I_J$ of a generic matrix $A$ be given.
Let $K \cap I = \emptyset$, $L \cap J = \emptyset$
and $K = L$. If every quasiminor $|A^U_V|_{uv}$ of $A^I_J$
in the identity ${\cal I}$ is replaced
by the quasiminor $|A^{U \cup K}_{V \cup L}|_{uv}$
of $A^{I \cup K}_{J\cup L}$ then we obtain a new identity ${\cal I}'$
called the extensional to ${\cal I}$.

\ppt {\bf Quasitelescoping sum.}~\index{quasitelescoping sum}
 Let $A = (a^i_j)$ be a generic $n \times n$
matrix. For any $k > 2$, and $i, j\in \{1, k-1\}$ consider the quasiminor
\[ |A^{i,k,k+1,\ldots,n}_{j,k,k+1,\ldots,n}|_{ij}. \]
The quasitelescoping sum involves such minors:
\[ QT(A^{1,\ldots,n}_{1,\ldots,n}) =
\sum_{k = 3}^{n} |A^{1, k, \ldots, n}_{k-1,k, \ldots, n}|_{1, k-1}
|A^{k-1, k, \ldots, n}_{k-1,k, \ldots, n}|^{-1}_{k-1, k-1}
|A^{k-1, k, \ldots, n}_{1,k, \ldots, n}|_{k-1, 1} \]
Then, by Muir's law and induction on $n$, we obtain
\begin{equation}\label{eq:quasitel-eq}
 QT(A^{1,\ldots,n}_{1,\ldots,n}) = a^1_1 - |A|_{11}.
\end{equation}
For $n=3$ this is simply the identity
obtained by extending by the third
row and column the identity
expressing the expansion of the $2\times 2$ upper left quasiminor.
Suppose now we have proved~(\ref{eq:quasitel-eq}) for $n$.
Take an $(n+1)\times(n+1)$-matrix $A$. Then, by induction,
this is true for the submatrix
\[ A^{\hat{2}}_{\hat{2}} = A^{1,3,\ldots,n}_{1,3,\ldots,n}. \]
But
\[\begin{array}{l}
 QT(A^{1,\ldots,n}_{1,\ldots,n}) \,=\, QT(A^{1,3, \ldots,n}_{1,3, \ldots,n})
+ |A^{1, 3, \ldots, n}_{2,3, \ldots, n}|_{1, 2}
|A^{2, 3, \ldots, n}_{2,3, \ldots, n}|^{-1}_{2, 2}
|A^{2, 3, \ldots, n}_{1,3, \ldots, n}|_{2, 1}\\
\,\,\,\,\,\,\,\,\,
=\, a^1_1 - |A^{1,3,\ldots,n}_{1,3,\ldots,n}|_{11} +
|A^{1, 3, \ldots, n}_{2,3, \ldots, n}|_{1, k-1}
|A^{2, 3, \ldots, n}_{2,3, \ldots, n}|^{-1}_{k-1, k-1}
|A^{2, 3, \ldots, n}_{1,3, \ldots, n}|_{k-1, 1}\\
\,\,\,\,\,\,\,\,\,
=\, a^1_1 - |A^{1,2, 3,\ldots,n}_{1,2,3,\ldots,n}|_{11}
\end{array}\]
where the last two summands were added up, using the identity
which represents the expansion of $2 \times 2$ upper left corner of $A$
and extending the identity by rows and columns having labels $3,\ldots,n$.

\ppt {\bf Homological relations.}~\index{homological relations}
Start with the identity
\[ (a^1_1 - a^1_2 (a^2_2)^{-1} a^2_1 )(a^2_1)^{-1}
= - (a^1_2 - a^1_1 (a^2_1)^{-1} a^2_2 )(a^2_2)^{-1}. \]
which in the quasideterminant language reads
\[ \left|\begin{array}{cc}\fbox{$a^1_1$}&{a^1_2}\\{a^2_1}&{a^2_2}
\end{array}\right|
(a^2_1)^{-1}
= -\left|\begin{array}{cc}\fbox{$a^1_1$}&{a^1_2}\\{a^2_1}&{a^2_2}
\end{array}\right| (a^2_1)^{-1} \]
and extend the latter applying Muir's law, adding the rows $3,\ldots,n$
of $A$ to each minor in the expression. Renaming the indices
arbitrarily we obtain the \wind{row homological relations}:
\begin{equation}\label{eq:row-hom-rel}\fbox{$
 |A|_{ij} |A^{\hat{i}}_{\hat{j'}}|_{i'j}^{-1} =
- |A|_{ij'} |A^{\hat{i}}_{\hat{j}}|_{i'j'}^{-1}$}
\end{equation}
for $j\neq j'$. Similarly, starting with the identity
\[
(a^1_2)^{-1}(a^1_1 - a^1_2 (a^2_2)^{-1} a^2_1 )
-(a^2_2)^{-1}(a^2_1 - a^2_2 (a^2_2)^{-1} a^2_1 ),
\]
we obtain the \wind{column homological relations}
\begin{equation}\label{eq:col-hom-rel}\fbox{$
 |A^{\hat{i'}}_{\hat{j}}|_{ij'}^{-1}|A|_{ij}  =
- |A^{\hat{i}}_{\hat{j}}|_{i'j'}^{-1}|A|_{i'j}.$}
\end{equation}

\ppt \wind{Laplace expansion for quasideterminants}. Start with the identity
\[ \sum_{j} a^i_j (A^{-1})^j_k = \delta^i_k. \]
If $i\neq k$ and $A^{-1}$ exists, then substituting
$(A^{-1})^j_k = |A|_{kj}^{-1}$ this becomes
\[ \sum_{j} a^i_j |A|_{ij}^{-1} = 1. \]
Multiply this equation from the right by $|A|_{il}$
and split the sum into the part with $j\neq l$ and the remaining term:
\[ a^i_{l} + \sum_{j \neq l} a^i_j |A|^{-1}_{ij}|A|_{il} = |A|_{il} \]
and apply the row homological relations~(\ref{eq:row-hom-rel})
to obtain the following Laplace expansion
for the $(i,j)$-th quasideterminant by the $k$-th row:
\begin{equation}\label{eq:row-Laplace-exp-quasid}\fbox{$
a^i_{l} - \sum_{j \neq l} a^i_j |A^{\hat{i}}_{\hat{l}}|^{-1}_{kj}
|A^{\hat{i}}_{\hat{j}}|_{kl} = |A|_{il}
$}
\end{equation}
Similarly, multiplying from the left the equation $\sum_{i}
|A|_{ij}^{-1} a^{i}_{j} = 1$ by $|A|_{lj}$ and splitting the sum
into two terms we obtain
\[ a^l_{j} + \sum_{i \neq l} |A|_{lj}|A|_{ij}^{-1} a^i_j = |A|_{lj}, \]
which after the application of the column homological relations
~(\ref{eq:col-hom-rel}) gives the following Laplace expansion
for the $(i,j)$-th quasideterminant by the $k$-th column:
\begin{equation}\label{eq:col-Laplace-exp-quasid}\fbox{$
a^l_{j} - \sum_{j \neq l}
|A^{\hat{i}}_{\hat{j}}|_{lk} |A^{\hat{l}}_{\hat{j}}|^{-1}_{ik}
a^{i}_l  = |A|_{lj}
$}
\end{equation}
Notice that the summation sign involves $(n-1)$ summands whereas
the similar summation
in the recursive formula~(\ref{eq:quasidet-recursive-def})
for quasideterminants involves $(n-1)^2$ summands.


\ppt ~(\cite{Cohn:free, Cohn:invloc}) 
Let $R$ be an associative unital ring, and
$\Sigma$ a given set of square {\it matrices} of
possibly different (finite) sizes with entries in $R$.
Map $f: R \rightarrow S$ of rings is $\Sigma$-inverting
if each matrix in $\Sigma$ is mapped to an invertible
matrix over $S$.  A $\Sigma$-inverting ring map
$i_\Sigma : R \rightarrow R_\Sigma$
is called {\bf Cohn localization} (or universal
$\Sigma$-inverting localization) if for
every $\Sigma$-inverting ring map $f: R \rightarrow S$
there exist a unique ring map $\tilde{f} : R_\Sigma \rightarrow S$
such that $f = \tilde{f}\circ i_\Sigma$.

A set $\Sigma$ of matrices is called (upper) {\bf multiplicative}
\index{multiplicative set!of matrices} if $1 \in \Sigma$ and, for
any $A,B \in \Sigma$ and $C$ of right size over $R$,
$\left(\begin{array}{cc}A & C \\ 0 & B \end{array}\right)$ is in
$\Sigma$. If $\Sigma$ is the smallest multiplicative set of
matrices containing $\Sigma_0$, then a map is $\Sigma_0$-inverting
iff it is $\Sigma$-inverting. Inclusion $\Sigma_0 \subset \Sigma$
makes every $\Sigma$-inverting map $f: R \rightarrow S$ also
$\Sigma_0$-inverting. Conversely, if each of the diagonal blocks
can be inverted, a block-triangular matrix can be inverted, hence
$\Sigma_0$-inverting maps are $\Sigma$-inverting.

The universal $\Sigma$-inverting localization can be constructed by
``invertive method'', as follows. Represent $R$ as a free
algebra $F$ on a generating set $\bf f$
modulo a set of relations $I$. For each quadratic matrix $A \in \Sigma$
of size $n\times n$, add $n^2$ generators $(A,i,j)$ to $\bf f$.
This way we obtain a free algebra $F'$ over some generating set $\bf f'$.
All $(A,i,j)$ for fixed $A$ clearly form a $n\times n$-matrix $A'$ over $F'$.
Then $\Sigma^{-1}R = F'/I'$ where $I'$ is the ideal generated by
$I$ and by all elements of matrices $AA'-I$ and $A'A-I$ for all
$A\in \Sigma$. Then $i_\Sigma : R\to \Sigma^{-1}R$ is the unique
map which lifts to the embedding $F \hookrightarrow F'$.

\ppt {\bf Warning.}\ldef{rationalclos}
 A naive approach to quotient rings,
would be just adding new generators $a'$ and relations $aa' = a'a
= 1$ for each $a\in R\backslash \{0\}$
which needs to be inverted in the first place.
In geometrical applications this could induce pretty unpredictable
behaviour on modules etc. But suppose we just want to do this in an
extreme special case: constructing a quotient skewfield.
After inverting all the nonzero elements,
we try inverting all their nonzero sums and so on.
The problem arises that one may not know which elements
from $m$-th step will be forced to zero by
new relations added a few steps later. So one should skip inverting some
new elements, as they will become zeros after a few more steps of inverting
other elements. There is no recipe which elements to leave out at
each step. For a given ring $R$, there may be none (no
quotient field) or multiple possibilities for such a recipe. More
precisely, given two embeddings $R\hookrightarrow K_i$ into
skewfields $K_1 \neq K_2$, there may be {\it different} smallest
subskewfields $L_i \hookrightarrow K_i$ containing $R$.

\ppt {\bf Proposition.} {\it
Let $\Sigma$ be multiplicative set of square matrices
over $R$ and $f : R \rightarrow S$ a $\Sigma$-inverting map.
Let $S(i,\Sigma)\subset R$ consists of all
components of solutions over $S$ of all equations $f(A)x = f(b)$
where $A \in \Sigma$, $b$ is a column-vector over $R$
and $x$ a column of unknowns.

(i) $S(i,\Sigma)$ is a subring of $S$.

(ii) $S(i,\Sigma)$ coincides with the image of $R_\Sigma$ under
the unique map $\tilde{f} : R_\Sigma\to S$ for which
$f = \tilde{f}\circ i_\Sigma$.

In particular, if $f$ is 1-1 then $\tilde{f}$
is isomorphism and $i_\Sigma$ is 1-1.
}

(i) If components $x_i$ and $y_j$ of column vectors $x$ and $y$
over $S$ are in $S(i,\Sigma)$, with $f(A)x = b$ and $f(B)y = c$,
then, possibly after enlarging $x,y,b,c$ by zeroes and $A$ and $B$ by
diagonal unit blocks, we may always make $i=j$ and $b$ and $c$ of
the same length. Then
$f\left(\begin{array}{cc}A&-A+B\\0&B\end{array} \right)
\left(\begin{array}{c}x+y\\y \end{array}\right) =
f\left(\begin{array}{c}b+c\\ c\end{array}\right)$ and as the
left-hand side matrix is in $f(\Sigma)$ by multiplicativity, then
$x_i + y_i\in S(i,\Sigma)$, as claimed. For $z$ a (row or column)
vector consider the diagonal square matrix ${\rm diag}(z)$ with
diagonal $z$. Then  ${\rm diag}(z)(1,1,\ldots, 1)^T = (z_1,\ldots,
z_n)$. For a fixed $i$, there is a matrix $P_i$ such that $P_i
(y_1,\ldots,y_n)^T = (y_i,\ldots,y_i)^T$. Hence,
 $f\left(\begin{array}{cc}B& -{\rm diag}(c)P_i \\0&A\end{array} \right)
\left(\begin{array}{c}y\\x \end{array}\right) =
f\left(\begin{array}{c}0\\b \end{array}\right)$ has as the $j$-th
component $(f(B)^{-1}f(c))_j (f(A)^{-1}f(b))_i$. But our
block-triangular matrix is in $\Sigma$, hence $x_i y_j$ is in
$S(i,\Sigma)$. Similarly, had we worked with algebras over $\genfd$,
we could have considered all possible weights on the diagonal
instead of just using the non-weighted diagonal ${\rm diag}(c)$
to obtain any possible $\genfd$-linear combination of
such.

(ii) The corestriction of $i$ onto $S(i,\Sigma)$ is also
$\Sigma$-inverting. Hence there is a unique map form $R_\Sigma$.
But, by construction, there is no smaller ring than $S(i,\Sigma)$
containing $f(R_\Sigma)$. As $i(R_\Sigma)$ is a ring they must
coincide. If the map is 1-1 it has no kernel hence $\tilde{f}$ is
an isomorphism.

\ppt {\bf Proposition.}
(left-module variant of {\sc P.~M. Cohn}~\cite{Cohn:invloc}, 2.1)
{\it If $\Sigma$ is multiplicative, then $\exists!$ subfunctor
$\sigma_\Sigma : R-{\rm Mod} \to R-{\rm Mod}$ of identity such that,
as a subset, $\sigma_\Sigma(M)$ equals
\[ \{ m \in M \,|\, \exists u = (u_1,\ldots, u_n)^T \in M^{\times n},
\, \exists i, \, m = u_i\mbox{ and } \exists A \in \Sigma,
\,\, A u = 0\} \]
for every $M \in R-{\rm Mod}$.
Moreover, $\sigma_\Sigma$ is an idempotent preradical.
}
\vskip .017in
{\it Proof.} 1. {\it $\sigma_\Sigma(M)$ is an $R$-submodule of $M$.}
It is sufficient to show that for any $r \in R$,
$m, m' \in \sigma_\Sigma(M)$ the left linear combination
$m + rm' \in \sigma_\Sigma(M)$. Choose
$A,B \in \Sigma$, $Au = 0, Bv = 0$, $u \in M^{\times k},
 v \in M^{\times l}$, $m = u_i, m' = v_j$.
 We may assume $k = l$, $i = j$, hence $m + rm' = (u + rv)_i$,
 by adjusting $A,B,u,v$.
For example,
 $\tilde{A} := {\rm diag}\,(I_s,A,I_t) \in \Sigma$
 and $\tilde{u} := (0_s, u_1, \ldots, u_k, 0_t)^T \in M^{s+k+t}$
 satisfy $\tilde{A} \tilde{u} = 0$ with $m =\tilde{u}_{i+s}$.
\newline Then
$\left(\begin{array}{cc} A & -Ar \\ 0 & B \end{array}\right)\in\Sigma$
and $\left(\begin{array}{cc} A & -Ar \\ 0 & B \end{array}\right)
\left(\begin{array}{c} u + rv \\ v \end{array}\right) =
\left(\begin{array}{c} 0 \\ 0 \end{array}\right)$.

2. {\it $M \mapsto \sigma_\Sigma(M)$ extends to a unique subfunctor
of identity.} If $m  = u_i\in \Sigma$
for some $i$ and $A(u_1,\ldots,u_k)^T = 0$ then
$A(f(u_1),\ldots, f(u_k))^T = 0$ whenever $f : M \to N$ is $R$-module map.
As $f(m) = f(u_i)$ this proves that
$f(\sigma_\Sigma(M))\subset \sigma_\Sigma(f(M))$ as required.

3. {\it $\sigma_\Sigma(M)$ is a preradical:
$\sigma_\Sigma(M/\sigma_\Sigma(M)) = 0$.}
If $m \in \sigma_\Sigma(M)$, then
$\exists u_1,\ldots, u_k \in M$, $\exists A \in \Sigma$, $\exists p \leq k$,
such that $A \vec{u}^T = 0 {\rm mod}\, \sigma_\Sigma(M)$
and $m  = u_{p}$
where $\vec{u} := (u_1,\ldots, u_k)$.
Hence $\exists v_1,\ldots,v_k \in \sigma_\Sigma(M)$
such that
$A\vec{u}^T = (v_1,\ldots, v_k)^T$ and
there are matrices $B_1,\ldots, B_k$, where $B_s$ is of size
$h_s \times h_s$, and vectors
$(w_{1s},\ldots,w_{h_s s})$ of size $h_s$, such that
$B_s (w_{1s},\ldots,w_{hs})^T = 0$ for all $s$; and we have
that $v_i = w_{s_i i}$ for some correspondence $i \mapsto s_i$.
Let $\vec{w} = (w_{11}, \ldots, w_{h_1 1}, w_{12},\ldots, w_{k s_k})$.
Let matrix $J = (J^i_j)$ be defined by
$J^i_{s_i} = 1$ for each $i$ and all other entries are $0$.
This matrix by construction satisfies $J\vec{w}^T = \vec{v}$.
Define also the block matrix $B := {\rm diag}\,\{B_1,\ldots, B_k\}$.
Clearly $B\vec{w}^T = 0$ by construction
and $B\in \Sigma$ by multiplicativity of $\Sigma$.
In this notation the summary of just said is encoded in this block identity
\[ \left(\begin{array}{cc} A & -J \\ 0 & B \end{array}\right)
(\vec{u}, \vec{v})^T = 0,\,\,\,\,\,\,
\left(\begin{array}{cc} A & -J \\ 0 & B \end{array}\right)\in \Sigma,\,\,\,
m = u_{p} .\]

4.  {\it $\sigma_\Sigma(\sigma_\Sigma(M)) = \sigma_\Sigma(M)$.}
If $m = u_i$ for some $i$ and $A (u_1,\ldots, u_k)^T = 0$
for some $A \in \Sigma$ with all $u_j \in \sigma_\Sigma(M)$, then
in particular, all  $u_j \in M$.

{\bf Exercise.}
{\it Let $\Sigma$, $\Sigma'$ be multiplicative sets of matrices
over $R$.

 If for every $A \in \Sigma$ there are permutation
matrices $w, w'\in GL(k,{\DDl Z})$ such that $wAw' \in \Sigma'$ then
$\sigma_{\Sigma'}(M) \subset\sigma_\Sigma(M)$ for all $M$.
}

\vskip .03in \ppt {\bf Warning.} $\sigma_\Sigma$ is not
necessarily left exact. Equivalently, the associated torsion
theory is not always hereditary (i.e. a submodule of a
$\sigma_\Sigma$-torsion module is not necessarily
$\sigma_\Sigma$-torsion). Hereditary torsion theories correspond
to Gabriel localizations.\ldef{hereditary3}

\ppt {\bf Quasideterminants vs. Cohn localization.}
Quasideterminants are given by explicit formulas. It is sometimes
more algorithmically manageable to invert them,
than the matrices (if the inverse can not be
expressed in terms of them anyway). The two procedures often
disagree, as some simple (e.g. diagonal) matrices do not
possess some among $n^2$ possible quasideterminants.
One may combine the process, by first inverting quasideterminants
which exist, and then performing the Cohn localization for the
simpler matrix so obtained.
The combination is not necessarily a Cohn localization.

Thus, let $\Sigma$ be as before. For each $A \in \Sigma$ and pair
$(ij)$ such that $|A|_{ij}$ exists and is nonzero add a variable
$B_{ij}$ and require $B_{ij}|A|_{ij} = |A|_{ij}B_{ij} = 1$. One
obtains a localization  $j : R \to R^{q0}_\Sigma$. Then one
inverts $j(\Sigma)$ by the Cohn method, which amounts to adding
formal variables just for those entries which are not added before
as quasideterminants, and adding relations for them.
The result is some localization
$i^q_\Sigma : R \to R^{q}_\Sigma$ which is $\Sigma$-inverting and
clearly a quotient ring of the Cohn localization. If $i^q_\Sigma$
is injective, it is just the Cohn localization.

There are obvious variants of this method,
(cf. reasoning in~\luse{rationalclos}, and recall that
quasideterminants may be defined inductively by size). Some rings
may be quotiented by ideals to get commutative or Ore domains.
A quasideterminant may be proven to be nonzero, as its
image in the quotient is nonzero, which is a good procedure for
some concrete $\Sigma$ (cf.\cite{Skoda:ban}, Th.7).

We have seen in Ch.8 that for the usual descent of quasicoherent sheaves
one needs flatness, which is often lacking for Cohn localization.
In the special case of Cohn localization at a 2-sided ideal,
the flatness of the localization map $i_\Sigma$ as a map of left modules
is equivalent to the right Ore condition. 
Though in geometrical situations one inverts sets of matrices
for which this theorem does not apply, flatness is not
expected for useful non-Ore universal localizations. Less essential,
but practically difficult, is to find the kernel of the localization
map $i_\Sigma$ (there is a criterium using the normal form mentioned below).

We would still like, in the spirit of an example~\cite{Skoda:ban},
Th.7 (more accurately described in~\cite{Skoda:nflag}), to be able
to consider some global noncommutative spaces where the local
coordinates are compared using nonflat Cohn localizations.
Localizations between full categories of modules ('perfect
localizations') \index{perfect localization} are described by
their underlying rings (the forgetful functor from the localized
category to the category of modules over localized ring is an
equivalence). Similarly, the knowledge of the restriction of the
Cohn localization functor to the category of finitely generated
(f.g.) projective modules is equivalent to the knowledge of the
localization morphism on the level of rings. Of course, the theory
is here not any more complicated if one inverts any multiplicative
set of morphisms between f.g. projectives than only matrices. The
localization functor for f.g. projectives has the explicit
description ({\sc Gerasimov-Malcolmson} normal
form~(\cite{Geras:loc,Malcolmson})) ~\index{normal form}
(analogous to the description of Ore localization as $S\times
R/\sim$ where $\sim$ is from ~\luse{relFrac}) and has an
interesting homological interpretation and
properties~(\cite{NR:nclchain}). Thus while the torsion theory
corresponding to Cohn localization is bad (nonhereditary,
cf.~\luse{hereditary},\luse{hereditary2},\luse{hereditary3}),
\index{hereditary torsion theory} other aspects are close to
perfect localizations (thus better than arbitrary hereditary
torsion theory). This suggests a hope in a geometry of ``covers by
Cohn localizations'' if we learn a way beyond the case of flat
descent for full categories of modules.


\vskip  .1in

\footnotesize{
{\bf Acknowledgements.} My biggest thanks go to Prof. {\sc A.~Ranicki},
who patiently encouraged the completion of this article
and gave me numerous small suggestions.

Many thanks to Prof. {\sc J.~Robbin} for educating me to
improve my writing style, to be straight,
explicit and to avoid a too frequent use of 'it is well-known that'.
Thanks to Indiana University, Bloomington
where a piece of my Wisconsin student 'preprint' started
growing into this paper, and Prof. {\sc V.~Lunts} of Indiana
for insights into mathematics and his manuscripts. I thank
the  Rudjer Bo\v{s}kovi\'{c} Institute, Zagreb, Max Planck
Institute, Bonn, and IH\'{E}S, Bures-sur-Yvette,
for allowing me to use a big chunk of my research time there
to complete this survey.
Finally, my apology to the reader that I needed
much more space and time for this article
to be polished and balanced than I could deliver.
}

\footnotesize{ {\bf Bibliography.} We try to give a few of the most
useful references for a reader who has similar geometric needs to
the author. The bibliography is intentionally incomplete. For
the benefit of the geometrically oriented reader, we used the
following preference criteria: geometrically motivated,
historically important, readable (for the author at least),
irreplacable. The literature which is obscure to me in a major way
is naturally not in the list. However I mention here some undoubtedly
important alternative works by listing only Math. Reviews code. For ring
theorists there is a monograph on torsion theories by {\sc Golan}
MR88c:16934 and on localization by {\sc Jatageonkar} MR88c:16005 and by
{\sc Golan} MR51:3207.
Many equivalent approaches to Gabriel localization
have been multiply discovered
({\sc Goldman} (1969)~\cite{Goldman},
{\sc Silver} (1967) MR36:205, 
{\sc Maranda} (1964) MR29:1236  
etc.) in various formalisms, e.g. torsion theories
(the term is basically from {\sc Dickson} (1965) MR32:2472).

Despite their historical importance, we ignore these, and
recommend the systematic treatment in Gabriel's
thesis~\cite{Gab:catab} as well as the
books~\cite{JaraVerschorenVidal,Popescu, stenstrom} and Ch.~6
of~\cite{BD:cats}. For \abelian categories see~\cite{Borceux,
Faith, Gab:catab, GM:hom, GM:met, JaraVerschorenVidal, Popescu,
Smith:nag, Tohoku, Weibel}; for localization in \abelian
categories see
books~\cite{Borceux,Faith,JaraVerschorenVidal,Popescu,Smith:nag}.
Other longer bibliographies on localization are
in~\cite{Faith,JaraVerschorenVidal,Lambek:bk,Popescu}.
Neither the present article nor the bibliography survey
noncommutative geometry beyond the localization aspects;
rather consult~\cite{Cart:madday,Connes:book,Var:el, Landsman:bookCQM,
Laudal:MPI,Manin:TNG,Rosen:book,SilvaWein,
Smith:nag,Staf:ICM,vdB:blow,vOyst:assalg,vOV:LNM887}
and bibliographies therein; quantum group literature\index{quantum groups}
(e.g.~\cite{Lusztig:qgrbk,ManinQGNG}; {\sc Chari, Pressley} MR95j:17010;
{\sc Majid} MR90k:16008, MR97g:17016, MR2003f:17014;
{\sc Klimyk-Schm\"udgen} MR99f:17017; {\sc V\'arilly} MR2004g:58006);
and, for the physics, also~\cite{dougnek}.
Abbrev.:
LMS = London Mathematical Society, MPI = Max Planck Inst. preprint
(Bonn). \& for Springer series: GTM (Graduate Texts in Math.), LNM
(Lecture Notes in Math.), Grundl.MW (Grundlehren Math.Wiss.).}

\noindent
Institute Rudjer Bo\v{s}kovi\'{c}, P.O.Box 180, HR-10002 Zagreb, Croatia
\bigskip

\noindent Institut des Hautes \'{E}tudes Scientifiques,\\
Le Bois-Marie, 35 route de Chartres, F-91440 Bures-sur-Yvette, France
\vskip .1in

\noindent Email: zskoda@irb.hr{\footnotesize
}

\printindex

\end{document}